\documentclass[11pt]{article} 

\usepackage{amsfonts} 
\usepackage{amsmath} 
\usepackage{amsthm} 
\usepackage{mathrsfs} 
\usepackage{dsfont} 
\usepackage{bbold}  
\usepackage{caption} 
\usepackage{float} 
\usepackage{tensor}

\DeclareMathOperator{\meas}{meas}
\DeclareMathOperator{\diag}{diag} 
\DeclareMathOperator{\dist}{dist}

\newtheorem{theorem}{Theorem} 
\newtheorem*{theorem*}{Theorem} 
\newtheorem*{prop*}{Theorem} 
\newtheorem{theo}[theorem]{Theorem} 
\newtheorem{coro}[theorem]{Corollary} 
\newtheorem{defi}[theorem]{Definition} 
\newtheorem{lemma}[theorem]{Lemma} 
\newtheorem{prop}[theorem]{Proposition} 
\newtheorem{rmk}[theorem]{Remark}

\newtheorem{hyp}{Hypothesis} 
\newtheorem*{hyp*}{Hypothesis}
 
 \newtheorem{example}{Example}

\newcommand{\zerarcounters}{\setcounter{equation}{0}\setcounter{theorem}{0}} 
 
\newcommand{\ZZZ}{\mathds{Z}} 
\newcommand{\CCC}{\mathds{C}} 
\newcommand{\NNN}{\mathds{N}} 
 
\newcommand{\RRR}{\mathds{R}} 
\newcommand{\TTT}{\mathds{T}} 
\newcommand{\uno}{\mathds{1}} 
 
\newcommand{\BB}{{\mathcal B}} 
\newcommand{\SSSSG}{\Lambda_{+}(G)}
\newcommand{\SSSSGa}{\Lambda_{+}(\Gacca)}
\newcommand{\CCCC}{{\mathcal C}} 
\newcommand{\DD}{{\mathcal D}} 
 
\newcommand{\calF}{{\mathcal F}} 
\newcommand{\calG}{{\mathcal G}} 
\newcommand{\calH}{{\mathcal H}} 
\newcommand{\calI}{{\mathcal I}}

\newcommand{\MM}{{\mathcal M}} 
\newcommand{\NN}{{\mathcal N}}

\newcommand{\RR}{{\mathcal R}} 
\newcommand{\SSSS}{{\Lambda_+}}

\newcommand{\calmM}{{\mathscr M}}

\newcommand{\gota}{{\mathfrak a}}

\newcommand{\gotd}{{\mathfrak d}} 
\newcommand{\gote}{{\mathfrak e}} 
\newcommand{\gotf}{{\mathfrak f}}

\newcommand{\gotn}{{\mathfrak n}}

\newcommand{\gots}{{\mathfrak s}}

\newcommand{\gotA}{{\mathfrak A}} 
\newcommand{\gotB}{{\mathfrak B}} 
\newcommand{\gotC}{{\mathfrak C}} 
\newcommand{\gotD}{{\mathfrak D}} 
 
\newcommand{\gotF}{{\mathfrak F}} 
\newcommand{\gotG}{{\mathfrak G}} 
\newcommand{\gotH}{{\mathfrak H}}

\newcommand{\gotK}{{\mathfrak K}}

\newcommand{\gotM}{{\mathfrak M}}

\newcommand{\gotS}{{\mathfrak S}} 
\newcommand{\gotT}{{\mathfrak T}}

\newcommand{\gotZ}{{\mathfrak Z}}

\newcommand{\pD}{\zeta}
 \newcommand{\pippa}{\mathtt c}
 \newcommand{\Pippa}{{\mathtt C}} 

\newcommand{\ol}{\overline} 
 
\newcommand{\Gacca}{\mathtt G}
\newcommand{\Fullbox}{{\rule{2.0mm}{2.0mm}}} 
 
\newcommand{\EP}{\hfill\Fullbox\vspace{0.2cm}} 
\newcommand{\prova}{\noindent{\it Proof. }} 
\newcommand{\io}{\infty} 
\newcommand{\e}{\varepsilon} 
\newcommand{\al}{\alpha} 
\newcommand{\de}{\delta} 
\newcommand{\be}{\beta} 
 
\newcommand{\m}{\mu} 
\newcommand{\x}{\xi}

\newcommand{\ka}{\kappa} 
\newcommand{\g}{\gamma} 
\newcommand{\om}{\omega} 
\newcommand{\h}{\eta} 
 
\newcommand{\la}{\lambda} 
\newcommand{\f}{\varphi} 
\newcommand{\s}{\sigma} 
 
\newcommand{\del}{\partial} 
\newcommand{\laa}{\langle} 
\newcommand{\raa}{\rangle}

\newcommand{\av}[1]{\langle #1 \rangle}

\newcommand{\oo}{{\omega}}

\newcommand{\ff}{\boldsymbol{f}}

\newcommand{\BBB}{\boldsymbol{B}} 
\newcommand{\bvert}{\boldsymbol{\vert}}

\newcommand{\ii}{{\rm i}}

 \newcommand{\se}{\mathfrak a}

\def\tilde#1{\widetilde{#1}}

\def\leftinv#1{\tensor*[^{[-1]}]{#1}{}}
 
\def\ins#1#2#3{\vbox to0pt{\kern-#2 \hbox{\kern#1 #3}\vss}\nointerlineskip}

\usepackage{hyperref}

  \makeindex

\begin{document} 
 
\title{\bf 
An abstract Nash-Moser theorem
and quasi-periodic solutions for NLW and NLS on 
compact \\ Lie groups and homogeneous manifolds} 
 
\author 
{\bf Massimiliano Berti$^{1}$, Livia Corsi$^{1}$, Michela Procesi$^{2}$
\vspace{2mm} 
\\ \small 
$^{1}$ Dipartimento di Matematica, Universit\`a di
Napoli ``Federico II'', Napoli, I-80126, Italy
\\ \small 
$^2$ Dipartimento di Matematica, Universit\`a di Roma ``La Sapienza", Roma, I-00185, Italy
\\ \small 
E-mail:  m.berti@unina.it, livia.corsi@unina.it, mprocesi@mat.uniroma1.it}
 
\date{} 
 
\maketitle 
 
\begin{abstract} 
We prove an abstract Implicit Function Theorem with parameters 
for smooth operators defined on sequence scales, modeled for the
search of  quasi-periodic solutions of PDEs. 
The tame estimates required for the inverse linearised operators at each step of the iterative scheme 
are deduced 
 via a multiscale inductive argument.  The Cantor like set of parameters 
where the solution exists is defined in a non inductive way. 
This formulation completely decouples the iterative scheme from the  measure theoretical analysis of the parameters
where the small divisors non-resonance conditions are verified. 
As an application, we deduce the existence of quasi-periodic solutions for forced NLW and NLS equations 
on {\it any} compact Lie group or manifold which is homogeneous
with respect to a compact Lie group, extending previous results valid only for tori.
A basic tool of harmonic analysis is the highest weight  theory  for the irreducible representations of compact
Lie groups.
\smallskip


\smallskip

\end{abstract} 

  \tableofcontents
 
\section{Introduction} 
\label{sec.intro} 

In the last years several works have been devoted to the search of quasi-periodic solutions of Hamiltonian  PDEs
in higher space dimensions, like analytic 
nonlinear Schr\"odinger  (NLS) and  nonlinear wave (NLW) equations  on $ \TTT^n $. 
A major difficulty concerns the verification of the so-called Melnikov non-resonance
conditions. 
The first successful approach, due to Bourgain \cite{B3}, \cite{B5},  
used a Newton iterative scheme  
which requires only the minimal (first-order) Melnikov conditions, which are verified inductively 
at each step of the iteration. 

In this paper we prove an abstract, differentiable Nash-Moser Implicit Function Theorem  with parameters
for smooth operators defined on Hilbert sequence scales.
 As applications, we prove the existence of quasi-periodic solutions with Sobolev regularity 
 of  the forced nonlinear wave equation 
\begin{equation}\label{nlw}
u_{tt}-\Delta u+mu = \e f(\oo t,x,u),
\qquad x\in {\mathtt M} \, , 
\end{equation}
and the nonlinear Schr\"odinger equation 
\begin{equation}\label{nls}
\ii u_{t}-\Delta u+mu = \e {\mathtt f}(\oo t,x,u),
\qquad x\in {\mathtt M} \, , 
\end{equation}
where $ {\mathtt M} $ is {\it any} compact Lie group or manifold which is homogeneous 
with respect to a compact Lie group, namely there exists a compact Lie group
which acts on ${\mathtt M}$ transitively and differentiably. 
In \eqref{nlw}-\eqref{nls} we denote by $ \Delta $ the Laplace-Beltrami operator,  
the ``mass" $ m>0 $, 
the parameter $\e >  0 $ is  small, and  the frequency vector 
$\oo\in\RRR^{d}$ is non-resonant, see \eqref{omegabar}-\eqref{diophquad} below. 

Examples of compact connected Lie groups are the standard torus $ \TTT^n $, the special
orthogonal group $ SO(n) $, the special unitary group $ SU(n) $, and so on. Examples of
(compact) manifolds homogeneous with respect to a compact Lie group are the spheres $ S^n $, the real
and complex Grassmannians,
and the moving frames, namely, the manifold of the $k$-ples of orthonormal vectors
in $ \RRR^n $ with the natural action of the orthogonal group $ O(n) $ and many
others; see for instance \cite{BTD}. 

The study of \eqref{nlw}-\eqref{nls} on a
manifold $ \mathtt M $   which is homogeneous with respect to a compact Lie 
group $ G $ is reduced to that  of \eqref{nlw}-\eqref{nls}  on the Lie group $ G $ itself. Indeed
$ \mathtt M $ 
is 
 diffeomorphic to 
$ \mathtt M =  G/ N $ where 
$ N $ is a closed  subgroup of $  G  $ and 
the Laplace-Beltrami operator on  $ \mathtt M $ can be identified with the Laplace-Beltrami operator 
on  $ G  $, acting on functions invariant under $ N $ (see \cite{BP}-Theorem 2.7 and \cite{H1,H2,Pro}).  

Concerning regularity we assume that the nonlinearity $ f \in C^{q}(\TTT^{d}\times{\mathtt M}\times\RRR; \RRR) $, resp.
 ${\mathtt f}(\varphi, x, u ) \in C^q(\TTT^d\times{\mathtt M}\times \CCC ;\CCC)$ in the real sense
 (namely as a function of ${\rm Re}(u), {\rm Im}(u)$),  for some $ q $ large enough. 
We also require that 
\begin{equation}\label{HamNLS}
{\mathtt f}(\oo t,x,u) = \partial_{\ol{u}} H(\oo t,x,u) \, , 
\quad H( \varphi,x,u) \in \RRR \, , \ \forall u \in \CCC \, , 
\end{equation}
  so that the NLS equation
 \eqref{nls} is Hamiltonian.

We assume that the frequency $\oo$ has a fixed direction,  namely 
\begin{equation}\label{omegabar}
\oo=\la\ol{\oo},\qquad \la\in\mathcal I:=[1/2,3/2],\qquad
|\ol{\oo}|_1:= {\mathop\sum}_{p=1}^d|\om_p| \le1,
\end{equation}
for some fixed diophantine vector $\ol{\oo}$, i.e. $\ol{\oo}$
satisfies 
\begin{equation}\label{dioph}
|\ol{\oo}\cdot l|\ge 2\g_{0} |l|^{-d}, \quad \forall\,l\in\ZZZ^{d}\setminus\{0\}.
\end{equation}
For the  NLW  equation \eqref{nlw} we assume also the quadratic diophantine condition
\begin{equation}\label{diophquad}
\Big|\sum_{1\le i,j\le d}\ol{\om}_{i}\ol{\om}_{j}p_{ij} \Big|
\ge\frac{\g_{0}}{|p|^{d(d+1)}},
\quad \forall\,p\in\ZZZ^{d(d+1)/2}\setminus\{0\}  
\end{equation}
which is  
satisfied for all  $|\overline \oo |_1 \le 1$ 
except a set of measure $O(\g_{0}^{1/2})$, see Lemma 6.1 in \cite{BB2}.

The search of quasi-periodic solutions of \eqref{nlw}-\eqref{nls} reduces to finding solutions $ u(\varphi, x) $  of 
\begin{equation}\label{NLSNLW}
(\omega \cdot \partial_\varphi)^2 u - \Delta u+mu = \e f(\varphi ,x,u) \, , 
\quad
\ii \om \cdot \partial_\varphi u -\Delta u+mu = \e {\mathtt f}(\varphi ,x,u,\ol{u}) \, ,  
\end{equation}
in some Sobolev space $ H^s $ of both the variables $ ( \varphi , x ) $, see Section \ref{hilbertspaces}.

\begin{theo}\label{teoremone}
Let  ${\mathtt M}$ be  {\rm any}  compact Lie group or manifold which is homogeneous with respect to a
compact Lie group. 
Consider the NLW equation \eqref{nlw}, 
and assume \eqref{omegabar}-\eqref{diophquad}; for the  
NLS equation \eqref{nls}-\eqref{HamNLS} assume only  \eqref{omegabar}-\eqref{dioph}.
Then there are $s,  q \in \RRR $ 
such that, for any $ f  $, $ {\mathtt f} \in C^{q}$ and
for all $\e\in[0,\e_{0})$ with $\e_{0}>0$ small enough, there is a map
$$
u_\e\in C^{1}([1/2,3/2],H^{s}),\qquad
\sup_{\la\in[1/2,3/2]}\|u_\e (\la)\|_{s}\to0,\mbox{ as }\e\to0,
$$
and a Cantor-like set $\CCCC_{\e}\subset[1/2,3/2]$, satisfying $ {\meas(\CCCC_{\e})}\to1\mbox{ as }\e\to 0 $, 
such that, for any $\la\in\CCCC_{\e}$, $u_\e (\la)$ is a solution of
\eqref{NLSNLW}, with $\oo=\la\ol{\oo}$. Moreover if 
 $f  $, $ {\mathtt f} \in C^{\io}$ then
the solution $u_\e (\la) $ is of class $ C^{\io} $ both in time and space. 
\end{theo}

Actually Theorem  \ref{teoremone} is deduced by  the abstract {\it Implicit Function Theorems \ref{principe}, \ref{mamma}}
(and {\it Corollary \ref{coromerd}}) 
on scales of Hilbert sequence spaces. 
We postpone their precise formulations   to Section \ref{see:abtheo}, since  some preparation is required.

Theorem \ref{teoremone} is  a first step in the direction of tackling the very hard problem of finding quasi-periodic solutions
for NLW and NLS on any compact Riemannian manifold, if ever true. This is an open problem also for periodic solutions. 
In the particular case that the manifold $ {\mathtt M} = \TTT^n $ is a $ n $-dimensional torus, Theorem \ref{teoremone} is proved in
\cite{BB1} for NLS, and, in \cite{BB2}, for NLW.

So far,  the  literature about quasi-periodic
solutions is restricted to NLS and NLW on 
tori 
(which are compact commutative Lie groups). 
The first results were proved for the interval $[0,\pi]$ by Kuksin \cite{K1, K2}, Wayne \cite{W1},  
P\"oschel \cite{Po2, KP}, 
or for the 1-dimensional circle $\TTT $ by Craig-Wayne \cite{CW}, Bourgain \cite{B1}, 
and Chierchia You \cite{CY}. For  higher dimensional tori $ \TTT^n $, $ n \geq 2 $, the first existence results have been obtained by 
Bourgain \cite{B3, B5} for NLS and NLW with Fourier multipliers via a multiscale analysis, recently applied by  
Wang \cite{W2} for completely resonant  NLS.  
Using KAM techniques,
Eliasson-Kuksin \cite{EK, EK1} proved existence and stability of quasi-periodic solutions for NLS with 
Fourier multipliers,  see also Procesi-Xu \cite{PX}. 
Then 
Geng-Xu-You \cite{GXY} proved KAM results for the cubic NLS in dimension $ 2 $ and 
Procesi-Procesi \cite{PP} in any dimension $ n $ and polynomial nonlinearity. 

The reason why these results are confined to tori 
 is that these proofs  require specific properties  of the eigenvalues and the eigenfunctions 
must be the exponentials or, at least,  strongly ``localized close to exponentials".
Recently,  Berti-Bolle \cite{BB1, BB2}
have extended the multiscale analysis to deal with  
NLS and NLW on $ \TTT^d $ with a multiplicative potential. In such a  case the 
eigenfunctions may not be  localized close to the exponentials. 

In the previous paper \cite{BP}, Berti-Procesi proved existence of periodic solutions for NLW and NLS on any compact  Lie group
or manifold homogenous  with respect to a compact Lie group. Main difficulties concern the  
eigenvalues of the Laplace-Beltrami operator,  with their unbounded multiplicity, 
and the rule of multiplications of the eigenfunctions. 
A key property which is exploited is that, for a Lie group, 
the product of two eigenfunctions is a finite linear combinations of them (as for the exponentials or the spherical harmonics).
From a dynamical point of view, it implies, roughly speaking, that 
 only finitely many normal modes are strongly coupled. 

Theorem \ref{teoremone}   extends the result in \cite{BP} to the harder quasi-periodic setting. As already said, 
it is deduced by  the abstract {\it Implicit Function Theorems \ref{principe}-\ref{mamma}}.
These results rely on the Nash-Moser iterative Theorem 
\ref{thm:nm1} and a multiscale inductive scheme for deducing 
tame estimates for the inverse linearised operators at each step of the iteration, see  Section \ref{section 5}. 
A main advantage of Theorem \ref{principe} is that 
  the Cantor-like set of parameters $ \CCCC_\e $ in \eqref{Cantorfinale} for which a solution exists
  is defined in terms of the  ``solution" $ u_\e $,
   and it is 
{\it not} inductively defined as in previous approaches.
This formulation completely decouples the Nash-Moser 
iteration  from the discussion about the measure of the  parameters where all the required
``non-resonance" conditions are verified. The possibility to impose the non-resonance conditions through 
the ``final solution"  was yet observed in \cite{BB06} (in a Lyapunov-Schmidt context) and in \cite{BBi10}
for a KAM theorem. In the present case the Cantor set $ \CCCC_\e $ is 
rather involved. 
Nevertheless we are able to provide efficient measure estimates in the applications.
This simplifies considerably the presentation because the measure estimates are not required at each step.
In conclusion, in order to apply Theorems \ref{principe}-\ref{mamma}, one does not need to know the multiscale
techniques nor the Nash-Moser approach: they can be used as a {\it black box}.

We believe that  Theorems \ref{principe}-\ref{mamma} can be applied  to several other cases.
The  abstract hypotheses can be verified by informations of the harmonic analysis on the manifold. 
In the case of compact Lie groups and homogeneous manifolds,  
Theorem \ref{teoremone}  follows by using only the harmonic analysis  in \cite{BP}
(see Section \ref{section 3}),
which stems from the informations on the eigenvalues and eigenspaces of the Laplace-Beltrami operator  
provided  by  the highest weight theory, see \cite{Pro}. 

We find it convenient to use a Nash-Moser scheme because the eigenvalues of the Laplacian are highly degenerate
and the second order Melnikov non resonance conditions required for the KAM reducibility scheme
might not be satisfied. 
\\[2mm]
{\it Informal presentation of the ideas and techniques.}
Many nonlinear PDEs (such as \eqref{NLSNLW}) can be seen as implicit function equations of the form
$$
F(\e,\la,u)=0,
$$
having for $\e=0$ the trivial solution $u(t, x) = 0 $, for each parameter $ \la \in {\cal I} $. Clearly, due to the small divisors,
 the standard Implicit Function Theorem fails, and one must rely on some Nash-Moser or KAM 
 quadratic scheme.
They are rapidly convergent iterative algorithms
based on the Newton method 
and hence  need some informations about the invertibility of the linearisation $L(\e,\la,u)$ of $F$ at any 
function $u$ close to zero. 

Due to the Hamiltonian structure, the operator $L(\e,\la,u)$ is self-adjoint and it
is  easy to obtain informations on its eigenvalues,  
implying the invertibility of $L(\e,\la,u)$ with bounds
of the $L^2$-norm of $L^{-1}(\e,\la,u) $ for ``most" parameters $ \la $. However  these informations
are not enough to prove the convergence of the algorithm: one needs
estimates on the high Sobolev norm of the inverse which  do not follow only
 from bounds on the eigenvalues. Usually this property is implied by a
 sufficiently fast polynomial off-diagonal decay of the matrices which represent the inverse operators.  

In the case of the interval $[0, \pi] $, the eigenvalues of $L(\e,\la,u)$ are often 
distinct, a property which enables to diagonalise $L(\e,\la,u)$ 
via a smooth change of variables (reducibility) implying very strong estimates of the 
inverse operator in high Sobolev norm. This method
automatically implies also the stability of the solution.
Unfortunately, the eigenvalues of $\Delta$ are not simple already on $\TTT$
(a fortiori neither on $\TTT^n $, $n\ge2$), so that generalising these {\it reducibility}
methods is complicated
and strongly depends on the equation. For NLS it is obtained in \cite{EK}.

However, 
the convergence of the Nash-Moser scheme only requires ``tame'' 
estimates of the inverse in high Sobolev norm (for instance like the one in \eqref{buonini}) which may
be obtained under weaker spectral hypotheses.  In the case of NLS and NLW 
on  $\TTT^n $ these estimates have been obtained in \cite{BB1, BB2} 
via a {\it multiscale analysis} on Sobolev spaces (see \cite{B5} in an analytic setting).
Informally, the multiscale method is a way to prove an off-diagonal decay 
for the inverse of a finite-dimensional invertible matrix with off-diagonal decay, by using informations on the
invertibility (in high norm) of a great number of principal minors of order $N$ much
smaller than the dimension of the matrix. The polynomial off-diagonal decay of a matrix implies that it
defines a ``tame'' operator between Sobolev spaces. 

In this paper we  extend these techniques also to the case
of compact Lie groups and  manifolds  which are homogeneous with respect to a compact Lie group
(in the latter case we ``lift-up"  the equation to the Lie group).
Two key points concern
\begin{enumerate}
\item  the matrix representation of a multiplication operator
$ u \mapsto b u $,
\item the properties of the eigenvalues of the Laplace-Beltrami operator. 
\end{enumerate}
The multiplication rules for the eigenfunctions, 
together with the numeration of the eigenspaces provided the highest weight theory, implies that 
the multiplication operator by a Sobolev function $ b \in H^s ({\mathtt M})$ is
represented in the eigenfunction basis as a block matrix with off-diagonal
decay, as  stated precisely in Lemmas \ref{mult.oper}, \ref{moltiplicazione.matrici} (proved in \cite{BP}). 
 The block structure of this matrix takes into account the (large) multiplicity of the
degenerate eigenvalues of $ \Delta $ on ${\mathtt M} $ (several blocks could correspond
to the same eigenvalue). 
This in principle could be a problem because one can not
hope to achieve any off-diagonal decay property for the matrices
 restricted to such blocks. However, as in \cite{BP}, 
 we do not need such a decay, being sufficient  to control only the $L^2$-operator norm on these blocks. 
Interestingly, properties of this type have been used by Bambusi, Delort, Gr\'ebert, Szeftel \cite{BDGS} 
 for Birkhoff normal form results of Klein-Gordon equations on Zoll manifolds
 (a main difficulty in \cite{BDGS}  is to verify the 
Birkhoff normal form non-resonance conditions). 

Concerning item 2, the eigenvalues of the Laplace-Beltrami operator on a Lie group
are very  similar to those on a torus, as stated in  \eqref{stok}. This enables to prove 
``separation properties'' of clusters of singular/bad sites (i.e. Fourier indices with a corresponding small divisor) 
\`a la Bourgain \cite{B3}, \cite{B5}. 
Thanks to the off-diagonal decay property proved in item 1 such ``resonant'' clusters 
 interact only weakly.   As in the case of $ \TTT^n $ (where the eigenvalues of $ - \Delta $ are  $|j|^2 $, $ j \in \ZZZ^n $) 
 one does not gather into the same cluster 
 all the indexes corresponding to the same eigenvalue. The reason is that such clusters would not satisfy
 the needed separation properties. 

\smallskip


We now give some more detail about the proof. 
In the usual PDE applications the function spaces decompose as a direct sum of eigenspaces of $ L ( 0 , \la  , 0 ) $, 
each of them being
a direct product of the exponentials
(for the time-direction) and the eigenspaces of $-\Delta+m$ (space direction). Hence we decompose
$u=\sum_{k}u_k$ with $k=(l,j)\in\ZZZ^d\times\Lambda_+ $ ($ l \in \ZZZ^d$ is the time-Fourier 
component and $ j \in \Lambda_+ $ the space-Fourier component). 
In particular $L(0,\la,0)$ is a diagonal operator which is proportional
to the identity on each eigenspace. Moreover the dependence on $l$  appears only through a
scalar function $\gotD_j( \om \cdot l)$, see \eqref{hyp:diag}.

In Theorems  \ref{principe}-\ref{mamma} 
we  revisit in a more abstract way the strategy of \cite{BB1,BB2}, obtaining a unified and more general result
for smooth operators $ F (\e, \la, u) $ acting on a  Hilbert scale of sequences spaces.
These results are based on  three hypotheses which allow the possibility of passing from $L^2$-norm estimates 
to high Sobolev norm bounds for the inverse linearized operators. We try to explain the meaning of this assumptions:

\begin{itemize}

\item[i.] the linearised operator can be written as the sum of a diagonal part $D=D(\la)$ (which is
the linearised operator at $\e=0,u=0$) plus a perturbation which has off-diagonal decay and is T\"oplitz
in the time indices (see Hypothesis \ref{matricione}),

\item[ii.] a uniform lower bound for the derivative of $\gotD_j(y)$ on the set where $ |\gotD_j(y)| $ is small 
(see Hypothesis \ref{hyp.misura.diag}), 

\item[iii.] an assumption on the length of chains of  ``singular sites''
(see Hypothesis \ref{hyp.separacatene}).

\end{itemize}

Under these hypotheses Theorem \ref{principe} implies the existence, for 
$\e$ sufficiently small, 
of a  function $u=u_\e(\la)$ which is a solution of
the equation $F(\e,\la,u)=0$ for all $\la$ in the Cantor-like set $\CCCC_\e$ in \eqref{Cantorfinale}, which is
defined only in terms of the eigenvalues of submatrices of $L(\e,\la,u_\e(\la))$.
Roughly speaking the set $\CCCC_\e$ is defined as the intersection of two families of sets:

\begin{itemize}

\item[1.] the sets $\ol{\gotG}_N$ of parameters $\la$ for which the $N$-truncation
of $L(\e,\la,u_\e)$  
 is invertible in  $ L^2 $ with good bounds of the $ L^2 $-norm of the inverse (see \eqref{buoninormal2finale}),

\item[2.] the sets $\ol{\calG}_N^0$ of parameters $\la$ for which the principal minors
of order $N$ having a small eigenvalue are {\it  separated} (see \eqref{buoniautovaloriinf}).

\end{itemize}

Technically the sets of type 2 are defined by exploiting the time-covariance property \eqref{topliz}
and analysing the complexity of the real parameter $\theta$ (see \eqref{innomin})
for which the $N$-truncation of the time-traslated matrix $L(\e,\la,\theta,u_\e)$ have a small eigenvalue:
since $\om = \la \overline \om $ is diophantine these two definitions are equivalent.

Finally we underline two main differences 
 with respect to the abstract Nash-Moser theorem in \cite{BBP}. 
  The first is that the tame estimates \eqref{buonini}  required for the inverse linearized operators 
 are  much weaker than in \cite{BBP}. Note, in particular, 
 that the tame exponent in  \eqref{buonini} grows like $ \sim \delta s $ (this corresponds  to an unbounded  loss of derivatives
as $ s $ increases).
 This improvement is necessary to deal with quasi-periodic solutions. 
 The second difference is that in \cite{BBP} the measure 
 issue was not yet completely decoupled from the Nash-Moser iteration, as, on the contrary, it is achieved in this paper
 thanks to the introduction of the set ${\cal C}_\varepsilon $ in \eqref{Cantorfinale}.

\smallskip

The paper is essentially self-contained.  The Appendix \ref{prova.multiscala} contains
the proof of the multiscale  proposition \ref{multiscala} which follows verbatim as in \cite{BB1}. 
We have added it for the convenience of the reader.

\smallskip

\noindent
{\it Acknowledgements}. We  thank L. Biasco, P. Bolle, C. Procesi for many useful comments. 

\zerarcounters 
\section{An implicit function theorem with parameters on sequence spaces}
\label{hilbertspaces}

We work on a scale of Hilbert sequence spaces defined as follows. We start from an index set
\begin{equation}\label{frakK}
\mathfrak K := \mathfrak I\times\gotA=\ZZZ^d\times \SSSS\times \gotA 
\end{equation}
where $\SSSS \subset \Lambda $ is contained in a $r$-dimensional lattice 
(in general not orthogonal)
\begin{equation}\label{cono}
\Lambda:= \Big\{j\in \RRR^r\,:\quad j= \sum_{p=1}^r j_p\mathtt w_p\,,\; j_p\in \ZZZ \Big\} 
\end{equation}
generated by independent vectors $\mathtt w_1,\ldots\mathtt w_r\in \RRR^r $.  
The set  $\gotA$ is finite, and in the applications will be  either $\gotA=\{1\}$ (for NLW) or 
$\gotA=\{1,-1\}$ (for NLS). Given 
$k\in \mathfrak K$ we denote
\begin{equation}\label{moduliki}
k=(i,\se)= (l,j,{\se}) \in \ZZZ^d\times \SSSS\times \gotA  \,,\quad |k|=|i|:= \max(|l|,|j|)\,,\quad |j|:=|j|_\io= \max_{p} |j_p|\in \NNN.
\end{equation}
If $\gotA =\{1\} $ we  simply write $k=(l,j)$.   

We require that
$\Lambda_+$ has a {\em product} structure, namely that
\begin{equation}\label{nonsodimeglio}
j=\sum_{p=1}^r j_p \mathtt w_p\,,j'=\sum_{p=1}^r j'_p \mathtt w_p\in \SSSS \quad \Rightarrow \; 
j''=\sum_{p=1}^r j''_p\mathtt w_p \in \SSSS\ \mbox{if } \forall\, p \ \min(j_p,j'_p)\le j''_p\le \max(j_p,j'_p).  
\end{equation} 
Condition \eqref{nonsodimeglio} will be used only in order to prove Lemma \ref{stronglygoodgood}.
It could be probably weakened. In the applications it is satisfied.  

To each $j\in \SSSS$ we associate a
``multiplicity" $ d_j \in \NNN $.  
Then, for $ s \ge 0 $, we define the
(Sobolev) scale of Hilbert sequence spaces
\begin{equation}\label{spazi}
H^{s}:=H^{s}({\mathfrak K} ):=\Big\{u=\!\!
\{u_{k}\}_{k\in {\mathfrak K}} \, , \ u_k  \in \CCC^{d_j}   \,:\, 
\|u\|^2_{s}:= 
\sum_{k\in{\mathfrak K}} \langle w_k\rangle ^{2s}
\|u_{k}\|_{0}^{2}<\infty \Big\}
\end{equation}
where $\| \ \|_{0}$ denotes the  $L^{2}$-norm in $\CCC^{d_j} $ and
the weights $\langle w_k \rangle := \max(c, 1, w_k  )$  satisfy
 $$ 
c |k| \leq w_k \leq C |k| \, , \quad \forall k \in{\mathfrak K} \, , 
 $$
 for suitable constants $ 0 < c \leq C $. In the applications 
 the weights $w_k$ are related to the eigenvalues of the Laplacian, see Examples \ref{zoppa}, \ref{zoppa2} below.
 
 \begin{rmk}\label{phil1}
The abstract Theorem \ref{principe} does not require any bound on the multiplicity $d_j$. 
In the applications we use the polynomial bound \eqref{pol-bound} 
for Lemmas \ref{mult.oper}, \ref{moltiplicazione.matrici} and for the measure  estimates.
 \end{rmk}
 
 
For any $B\subseteq{\mathfrak K}$ we define  the subspaces
\begin{equation}\label{HB subspaces}
H_{B}^{s}:=\left\{u\in H^{s}\,:\, 
u_k=0 \mbox{ for }k\notin B\right\} \, .
\end{equation}
If $B$ is a finite set the space $H_{B}^{s} = H_B $ does not depend
on $s$ and it is included in $ \cap_{s\geq 0}H^s $.

Finally, for $ k = ( i , \se ), k' = (i' , \se' ) \in{\mathfrak K}$ we denote
\begin{equation}\label{modulo}
\dist(k,k'):=\left\{
\begin{aligned}
&1,\qquad\qquad i=i',\,\se\ne\se',\\
& |i-i'|,\qquad\mbox{otherwise} \, , 
\end{aligned}
\right.
\end{equation}
where $ |i|$ is defined in \eqref{moduliki}.

\begin{rmk}\label{vday}
In principle $i-i'$ { may not be in }$ \ZZZ^d \times \Lambda_+ $ { because }$\SSSS$
{ is not a lattice}. {However, since }$\SSSS\subset\Lambda$ { we can always
compute }$|i-i'|$ { by considering }$i-i'\in\ZZZ^d \times\Lambda$.
In order to avoid this problem we will extend our vectors by setting them to zero on 
$ (\ZZZ^d\times \Lambda\times \gotA) \setminus \gotK$. 
\end{rmk}

All the constants that will appear in  the sequel may depend on the index set $\gotK$,
the weights $w_k$ and on $s$. We will  evidence only the dependence on $s$.

\subsection{Linear operators on \texorpdfstring{$H^s$}{Hs} and matrices}
\label{sub.linearop}

Let $B,C\subseteq{\mathfrak K}$. 
A bounded linear operator $ L : H_B^s \to H_C^s $ is represented, as usual,
 by a matrix in
\begin{equation}\label{lematrici}
\MM^B_C := 
\Big\{\big(M_{k}^{k'}\big)_{k\in C, k'\in B},\,
M_{k}^{k'}\in {\rm Mat}(d_{j}\times d_{j'},\CCC)\Big\}.
\end{equation}
It is useful to evidence a bigger block structure. 
We decompose
$$
B= {\overline B} \times \gotB\,,\quad  \overline{B}:={\rm Proj}_{\ZZZ^{d}\times\SSSS} B\,,\quad 
\gotB:={\rm Proj}_{\gotA} B\,
$$
and  $C= \overline C\times \gotC$, defined in the same way. Now, for $i=(l,j)\in \overline C$, $i'=(l',j')\in \overline B $, we 
consider the matrix
\begin{equation}\nonumber
M_{\{i\}}^{\{i'\}}:= \{M_{i,\gota}^{i',\gota'}\}_{\gota\in \gotC,\gota'\in \gotB} \,,\quad 
M_{\{i\}}^{\{i'\}}\in {\rm Mat}( |\gotC | d_{j}\times|\gotB| d_{j'},\CCC) \, , 
\end{equation}
where $ |\gotB|$, $ |\gotC | $ denote the cardinality of $ \gotB, \gotC \subseteq \gotA $ respectively. 
In the same way, given a vector $v:=\{v_k\}_{k\in \overline C\times \gotC}$, for $i=(l,j)\in \overline C$,  we set
$ v_{\{i\}} :=\{v_{i,\gota}\}_{\gota\in \gotC}$.

\begin{rmk}\label{dimblo}
The difference with respect to \cite{BB1},  \cite{BB2} 
is that the dimension of  the {\it matrix blocks } $M_{\{i\}}^{\{i'\}}$ 
may not be uniformly bounded. 
They are  scalars for the NLW equation in \cite{BB2} and, in \cite{BB1}, for NLS, at most  
$ 1 \times 2 $, $ 2 \times 1 $ or $ 2 \times 2 $ matrices, 
because $ d_j = d_{j'} = 1 $ and  $ 1 \leq |\gotB|, |\gotC| \leq 2$.
\end{rmk}

We endow ${\rm Mat}(|\gotC|d_{j}\times|\gotB| d_{j'},\CCC)$  with the  $L^{2}$-operator norm, which we denote $\|\cdot\|_0$.  
Note that whenever a multiplication is possible one has the algebra property.

\begin{defi}\label{def:Ms} {\bf ($s$-decay norm)}
For any $M\in \MM^{B}_{C}$ 
we define its $s$-norm   
\begin{equation}\label{s-norm}
\bvert M\bvert_{s}^{2}:=K_{1}\sum_{i\in{\mathfrak \ZZZ^{d}\times\Lambda}}[M(i)]^{2}
\langle i\rangle^{2s}
\end{equation}
where  $\av{i} := \max(1,|i| )$, 
\begin{equation}\label{maggiorante}
[M(i)]:=\left\{
\begin{aligned}
&\sup_{\substack{h-h'=i, h\in \ol{C},\,h'\in \ol{B}}} \big\|M^{\{h'\}}_{\{h\}} \big\|_0,
 &i\in \ol{C}-\ol{B}, \\
&\qquad 0, & i\notin \ol{C}-\ol{B}\,,
\end{aligned}
\right.
\end{equation}
and  $K_{1}> 4 \sum_{i\in \ZZZ^d\times \Lambda} \av{i}^{-2s_0} $.

We denote by $ (\MM^s)_{C}^{B}\subset \MM_C^B $ the set of matrices with finite $s$--norm
$\bvert \cdot\bvert_{s}$.
If $B, C $ are finite sets  then $(\MM^s)_{C}^{B} = \MM_{C}^{B} $ does not  depend on $s$, and, for simplicity, 
we drop the apex $ s $.  
\end{defi}

Note that the norm $\bvert\cdot\bvert_{s}\le\bvert\cdot
\bvert_{s'}$ for $s\le s'$.

The norm defined in \eqref{s-norm} is a variation of that introduced in Definition 3.2 of \cite{BB1}. 
The only difference concerns the dimensions of the blocks $M_{\{i\}}^{\{i'\}}$ as noted in Remark \ref{dimblo}.
However, since  the matrices  $M_{\{i\}}^{\{i'\}}$ are measured with the operator norm $ \| \cdot \|_0 $ 
the algebra and interpolation properties  of the norm $\bvert \cdot \bvert_s $ follow similarly to 
\cite{BB1}, as well as all the properties in section 3-\cite{BB1}. 
Indeed, given $  M \in \MM^{B}_{C} $ we introduce the T\"opliz matrix
\begin{equation}\label{toplizzizata}
 \calmM := (  \calmM^{\{i'\}}_{\{i\}} ) \in \MM^{B}_{C}  \, , \quad  \calmM^{\{i'\}}_{\{i\}} := [ M (i-i')] \uno_{|\gotC| d_j \times |\gotB| d_{j'}}
\end{equation}
which has the same decay norm
\begin{equation}\label{normeuguali}
\bvert M \bvert_s = \bvert \calmM \bvert_s \, . 
\end{equation}
\begin{lemma}\label{lem:mag}
 Let $ M_1 \in {\cal M}^C_D $ and $ M_2 \in {\cal M}^B_C $. Then $ M_1 M_2 \in {\cal M}^B_D $ satisfies
$\bvert M_1 M_2 \bvert_s \leq \bvert \calmM_1 \calmM_2 \bvert_s $.
\end{lemma}

\prova
For $ i' \in \ol{B} $, $ i \in \ol{D} $, we have
\begin{eqnarray}\nonumber 
\Big\| (M_1 M_2)_{\{ i \}}^{\{ i' \}} \Big\|_0 &\leq &  
\sum_{q \in \ol{C}} \Big\| (M_1)_{\{ i \}}^{\{ q \}} \Big\|_0 \Big\| (M_2)_{\{ q \}}^{\{ i' \}} \Big\|_0 
 \stackrel{\eqref{maggiorante}} \leq  \sum_{q \in \ol{C}}  [(M_1)(i-q)] [(M_2)(q-i')] \, 
 \\
 & = &  \!\!\sum_{q \in \ol{C}}   [(M_1)(i-q)] [(M_2)(q-i')]  \Big\| \uno_{|\gotD|d_{j} \times |\gotB| d_{j'}}\Big\|_0 \nonumber
 \\ &=& \Big\| \sum_{q \in \ol{C}}  [M_1(i-q)] \uno_{|\gotD |d_j \times |\gotC|d_{j_q}} 
  [M_2 (q-i')] \uno_{|\gotC|d_{j_q} \times |\gotB| d_{j'}} \Big\|_0 \stackrel{\eqref{toplizzizata}}  = \Big\| (\calmM_1 \calmM_2)_{\{ i \}}^{\{ i' \}} \Big\|_0 \nonumber 
 .
\end{eqnarray}
Therefore
$ [ (M_1 M_2)(i-i')] \leq \| (\calmM_1 \calmM_2)_{\{ i \}}^{\{ i' \}} \|_0 $ and the lemma follows. 
\EP

In what follows we fix $s\ge s_0 >(d+r)/2$. 

\begin{lemma}\label{lem.algebra.bvert} {\bf (Interpolation)}
For all $s\ge s_{0} $ there is $C(s)>1$ with $C(s_{0})=1$
such that, for any  subset $B,C,D\subseteq{\mathfrak K}$
and for all $M_{1}\in\MM^{C}_{D}$, $M_{2}\in\MM^{B}_{C}$, one has
\begin{equation}\label{nanbnbna.bvert}
\bvert M_{1}M_{2}\bvert_{s}\le \frac{1}{2}\bvert M_{1}\bvert_{s_{0}}
\bvert M_{2}\bvert_{s}
+\frac{C(s)}{2}\bvert M_{1}\bvert_{s}\bvert M_{2}\bvert_{s_{0}}.
\end{equation}
In particular, one has the algebra property 
$ \bvert M_{1}M_{2}\bvert_{s}\le C(s) \bvert M_{1}\bvert_{s} \bvert M_{2}\bvert_{s} $.
\end{lemma}

\prova For the T\"opliz matrices $ \calmM_1 $, $ \calmM_2 $ the interpolation inequality \eqref{nanbnbna.bvert} 
follows as usual
(with $ C(s_0) \leq 1 $ possibly taking $ K_1  $ larger).
Hence Lemma \ref{lem:mag} and \eqref{normeuguali} imply \eqref{nanbnbna.bvert}. 
\EP

The $ s -$norm of a matrix also controls the $ \| \  \|_s $ norm (see \cite{BB1}-Lemma 3.5).

\begin{lemma}\label{lem.questo}
For any $B,C\subseteq\gotK$,
let $M\in \mathcal M^B_C $.  Then 
\begin{equation}\label{soboh}
\| M h\|_s \leq C(s)\bvert M\bvert_{s_0}\|h\|_{s}+ C(s)\bvert M\bvert_{s}\|h\|_{s_0} \, , \quad \forall h\in H_B^s \, . 
\end{equation}
\end{lemma}

\prova
Regarding a vector $ h =\{ h_k\}_{k\in \overline B \times \gotB}$ as a column matrix, 
its $s$-decay norm  is 
$ \bvert h \bvert_{s}^2 =K_1\sum_{i\in\overline B} \av{i}^{2s}\| h_{\{i\}}\|^{2}_0$.
Hence \eqref{soboh} follows by Lemma \ref{lem.algebra.bvert} because 
$ c(s)\|h\|_s \leq \bvert h\bvert_s \leq c'(s)\|h\|_s$. 
\EP

We conclude this section stating further properties of the $s$-norm: such
lemmata are proved word by word as Lemmas 3.6, 3.7 3.8 and 3.9 of \cite{BB1} respectively.

\begin{lemma}\label{lem.decadimento} {\bf (Smoothing)}
Let $M\in\MM^{B}_{C}$ and $N\ge2$. For all $s'\ge s\ge0$ the following hold.

\noindent
(i) If $M^{k'}_{k} = 0 $ for all $ \dist(k', k) <N$
(recall the definition \eqref{modulo}), then
\begin{equation}\label{fuoridiag}
\bvert M\bvert_{s}\le N^{-(s'-s)}\bvert M\bvert_{s'}.
\end{equation}
\noindent
(ii) If $M^{k'}_{k} = 0 $  for all $ \dist(k', k) > N$, then
\begin{equation}\label{dentrodiag.a} 
\bvert M\bvert_{s'}\le N^{s'-s}\bvert M\bvert_{s} \, ,
\qquad 
\bvert M\bvert_{s}\le N^{s+d+r}\| M\|_{0} \, .  
\end{equation}
\end{lemma}

\begin{lemma}\label{decadonolerighe} {\bf (Decay along lines)}
Let  $M\in\MM^{B}_{C}$ and denote by $ M_k $, $ k\in C $,  its  $ k $-th line. Then 
\begin{equation}\label{eq.decadonolerighe}
\bvert M\bvert_{s} 
\le |\gotC| K_2\max_{k\in C}\bvert M_{k}\bvert_{s+d+r}, \quad \forall s\ge0 \, .    
\end{equation}
\end{lemma}

\begin{lemma}\label{controlloautov}
Let $M\in\MM^{B}_{C}$. Then $ \|M\|_{0}\le\bvert M\bvert_{s_{0}} $.
\end{lemma}

\begin{defi}\label{def.leftinv}
We say that a matrix $M\in\MM^{B}_{C}$ is \emph{left invertible} if there
exists $N\in\MM^{C}_{B}$ such that $NM=\uno_{B}$. In such a case $N$ is
called a \emph{left inverse} of $M$.
\end{defi}

A matrix $M$ is left-invertible if and only if it is
injective. The left inverse is, in general, not unique.
In what follows we shall denote by $\leftinv{M}$ any left inverse of $M$
when this does not causes ambiguity.

\begin{lemma}\label{lem.inversasinistra}{\bf (Perturbation of left-invertible matrices)}
Let $M\in\MM^{B}_{C}$ be a left invertible matrix. Then for any
$P\in\MM^{B}_{C}$ such that $\bvert \leftinv{M}\bvert_{s_{0}}\bvert P\bvert_{s_{0}}
\le1/2$ there exists a left inverse of $M+P$ such that
\begin{equation}\label{normaltaleftinv}
\bvert \leftinv{(M+P)}\bvert_{s_{0}}\le 2\bvert\leftinv{M}\bvert_{s_{0}},
\quad \bvert \leftinv{(M+P)}\bvert_{s} \le 
C(s)\big(\bvert \leftinv{M}\bvert_{s}+\bvert\leftinv{M}\bvert_{s_{0}}^{2}
\bvert P\bvert_{s}\big) , 
\end{equation}
 for any $s\ge s_{0} $.  Moreover, if 
 $\|\leftinv{M}\|_{0}\|P\|_{0}
\le1/2$, then 
there is a left inverse of $M+P$ which satisfies
\begin{equation}\label{norma0leftinv}
\|\leftinv{(M+P)}\|_{0}\le 2\| \tensor*[^{[-1]}]{M}{}\|_{0}.
\end{equation}
\end{lemma}

\subsection{Main abstract results} 
\label{see:abtheo}

We consider a non-linear operator
\begin{equation}\label{eq.vera}
F(\e,\la,u)=D(\la)u+\e f(u) 
\end{equation}    
where $\e>0$ is small, the parameter
$\la\in{\mathcal I}\subset [1/2, + \infty) $, and $ D(\la)$ is a diagonal linear operator 
$D(\la):H^{s+\nu}\to H^{s}$ such that 
\begin{equation}\label{Hdiag}
\| D(\la) h \|_s, \| \partial_\la D(\la) h \|_s  \leq C(s) \| h \|_{s+ \nu}
\end{equation}
(in the applications $D(\la)=\ii\la\ol{\om}\cdot\del_{\f}-\Delta + m $ or $ (\la\ol{\om}\cdot\del_{\f})^{2}-\Delta + m $) whose action on the subspace associated to a fixed index $k$ is scalar, namely 
\begin{equation}\label{Ddiag}
D(\la)=\diag(D_k(\la)\uno_{d_j})_{k\in \mathfrak K}  \, .
\end{equation}
We assume that, for some $ s_0 >  (d+r)/2$, the nonlinearity
$ f \in C^2(B^{s_0}_1, H^{s_0})$ (where  $B^{s_0}_1$ denotes the unit ball in $H^{s_0}$) 
and   the following ``tame" properties hold: 
given $ S' > s_0 $, for all $s \in[s_0, S')$ there exists a constant $C(s)$  such that for any 
$\|u\|_{s_0} \le2 $, 
\begin{itemize}

\item[(f1)] $\| d f(u)[h]\|_{s}\le C(s) \big( \|u\|_{s} \|h\|_{s_0} +  \|h\|_{s} \big) $,

\item[(f2)] $\| d^{2} f(u)[h,v]\|_{s} \le C(s)\Big(
\|u\|_{s}\|h\|_{s_0}\|v\|_{s_0}+\|h\|_{s}\|v\|_{s_0}
+\|h\|_{s_0}\|v\|_{s}\Big)$
\end{itemize}
hold.
Our goal is to find $u=u_\e (\la)\in H^{s}$ for suitable $s$ which solves the equation
$ F(\e,\la,u_\e (\la)) = 0 $ at least for ``some'' values of $\la \in \mathcal I$.

Then we assume further properties on  the linearised operator
\begin{equation}\label{operatore}
L=L(\e,\la,u)=D(\la)+\e T(u), \quad D(\la)=\diag(D_k(\la)\uno_{d_j})_{k\in \mathfrak K}. 
\end{equation}
where $T(u)$  is the matrix which represents the bounded linear operator $df(u) $, see \eqref{lematrici}.

\begin{hyp}\label{matricione}  
Let $\ol{\om} \in \RRR^d $ satisfy \eqref{dioph}. There exists 
a function $ \gotD: \SSSS\times\gotA\times \RRR\to \CCC$ and $\nu_0>0$  such that 
\begin{subequations}
\begin{align}
 & \text{\bf (Covariance)} \ \qquad
 D_{(l,j,\se)}(\la)=\gotD_{j,\se}(\la\ol{\om}\cdot l)\, , \forall \la\in \mathcal I  
\label{hyp:diag}\\
& \text{\bf (T\"oplitz in time)} \ \qquad
T\in \mathcal M_{{\mathfrak K}}^{{\mathfrak K}}\;: T_{(l,j,\se)}^{(l',j',\se')}= T_{(j,\se)}^{(j',\se')}(l-l')
\label{hyp.topliz}\\
&  \text{\bf (Off-diagonal  decay)} \ \qquad 
 \bvert T(u)\bvert_{s-\nu_0}\le C(s)(1+ \|u\|_{s}) \, ,
\label{hyp.normaT} \\
&  \text{\bf (Lipschitz)} \ \qquad
 \bvert T(u)-T(u')\bvert_{s-\nu_0} \le C(s)(\|u-u'\|_s +(\|u\|_s + \|u'\|_s)\|u-u'\|_{s_0}) \, ,  \label{lip}
\end{align}
for all $ \|u\|_{s_0}, \|u'\|_{s_0}\le 2$ and $ s_0+\nu_0 < s < S' $. 
\end{subequations}
\end{hyp}
For any $\theta\in\RRR$ we set
\begin{subequations}
\begin{align}
&D(\la,\theta)= {\rm Diag}(D_k(\la,\theta) \uno_{d_j})\,,\quad D_k(\la,\theta):= \gotD_{j,\se}(\la\ol{\om}\cdot l+\theta)\\
&L(\e,\la,\theta,u):= D(\la,\theta)+\e T(u)\, .  
\end{align} \label{innomin} 
\end{subequations}
We need the following information about the unperturbed small divisors. 
\begin{hyp}\label{hyp.misura.diag} {\bf (Initialisation)}
There are $ \gotn $ such that for all $\tau_1 > 1$, $N > 1 $, $\lambda\in \mathcal I$, 
$  l \in \ZZZ^d $, $  j \in \Lambda_+ $, $ \gota \in \gotA $, the set
\begin{equation}\label{meas.diag}
\big\{\theta \in \RRR\;:\; |D_{(l,j,\se)}(\lambda,\theta)|\le N^{-\tau_{1}}  \big\} 
\subseteq \bigcup_{q=1}^\gotn I_{q} \quad \mbox{intervals with }  \meas(I_q) \le  N^{-\tau_1} \, . 
\end{equation}
\end{hyp}
We now distinguish which unperturbed small divisors are actually small or not. 
\begin{defi}\label{regular} {\bf (Regular/singular sites)}
We say that the index $k\in{\mathfrak K}$ is
\emph{regular} for a matrix $ D := {\rm diag}( D_k \uno_{d_j}) $, $ D_k \in \CCC $,  if $|D_{k} |\ge 1$, otherwise 
we say that $k$ is \emph{singular}.
\end{defi}
We need an assumption which provides separation properties of clusters of singular sites. 
 
For any $\Sigma\subset \gotK$ and $\widetilde{\jmath}\in \SSSS$ we denote 
the section of $\Sigma$ at fixed $\widetilde{\jmath}$ by
$$
\Sigma^{(\widetilde{\jmath})}:=\{k=(l,\widetilde{\jmath},\se)\in \Sigma\} \,.
$$
\begin{defi}\label{fibre}
Let $ \theta, \la $ be fixed and $ K > 1 $. We denote by $\Sigma_K$ any subset of singular sites of $D(\la,\theta)$
in $\gotK$ such that, for all $\widetilde{\jmath}\in \SSSS $, the cardinality of the section $  \Sigma^{(\widetilde{\jmath})}_K $
satisfies  $  \# \Sigma^{(\widetilde{\jmath})}_K \le K$.
\end{defi}

\begin{defi}\label{gammachain} {\bf  ($ \Gamma$-Chain)} Let 
$\Gamma\ge 2 $. 
A sequence $k_{0},\ldots,k_{\ell}\in\gotK$ with
$k_{p}\ne k_{q}$ for $0\le p\ne q\le \ell$ such that
\begin{equation}\label{Gamma-chain}
\dist( k_{q+1},k_{q})\le \Gamma,\qquad
\mbox{ for all }q=0,\ldots,\ell-1,
\end{equation}
is called a $\Gamma$-chain of length $\ell$.
\end{defi}

\begin{hyp}\label{hyp.separacatene} {\bf (Separation of singular sites)}
There exists a constant $ {\mathtt s}$ and, for any $ N_0 \ge 2 $, 
a set  $ \tilde{\mathcal I} = \tilde{\mathcal I}(N_0)$ such that, 
for all $ \la\in \tilde{\mathcal I} $,  $ \theta \in \RRR $, and for all 
$ K, \Gamma $ with $ K \Gamma \geq N_0 $, any $\Gamma$-chain of
singular sites in $\Sigma_K$ as in Definition \ref{fibre},  has length
$ \ell\le (\Gamma K)^{{\mathtt s}} $.
\end{hyp}

In order to perform the multiscale analysis we need finite dimensional truncations
of the matrices. 
Given a parameter family of matrices $L(\theta)$ with $\theta\in\RRR$ and $N>1$ 
for any $k=(i,\gota)=(l,j,\gota)\in\gotK$ we denote by $L_{N,i}(\theta)$ (or equivalently $L_{N,l,j}(\theta)$) the 
sub-matrix of $L(\theta)$ centered at $ i $, i.e.
\begin{equation}\label{troncato}
L_{N,i}(\theta):=L(\theta)_{F}^F\,,\quad F:=\{k'\in \gotK\,: \; \dist(k,k')\leq N\}.
\end{equation}
If $l=0$, instead of the notation \eqref{troncato}  we shall use the notation
$$
L_{N,j}(\theta):=L_{N,0,j}(\theta) \, ,
$$
if also $j=0$ we write
$$
L_{N}(\theta):=L_{N,0}(\theta) ,
$$
and for $\theta=0$ we denote
$ L_{N,j}:=L_{N,j}(0) $. 

By hypothesis \ref{matricione}, the matrix $L=L(\e,\la,\theta,u)$ has the following covariance property in time
\begin{equation}\label{topliz}
L_{N,l,j}(\e,\la,\theta,u)=L_{N,j}(\e,\la,\theta+\la\ol{\om}\cdot l,u).
\end{equation}
For $ \tau_0 > 0 $, $N_0\ge 1$  we define the set
\begin{equation}\label{lambdabarra}
\ol{\mathcal I} := \ol{\mathcal I}(N_0,\tau_0):= \Big\{ \la\in\mathcal I\;:\; |D_{k}(\la)|
\ge N^{-\tau_{0}}_{0} \mbox{ for all }  k= (i,\gota)\in \gotK: \;|i|\le N_{0} \Big\}.
\end{equation}
\begin{theo}\label{principe}
Let $ \gote > d + r + 1  $. 
Assume that  $ F $ in \eqref{eq.vera}  satisfies \eqref{Hdiag}-\eqref{Ddiag}, (f1)--(f2) and
Hypotheses \ref{matricione}, \ref{hyp.misura.diag}, \ref{hyp.separacatene} 
with $ S' $ large enough, depending on $ \gote $.
Then, there are  $\tau_1>1$, $\ol{ N}_0 \in \NNN $,  
 $ s_1 $, $ S \in (s_0 + \nu_0, S'-\nu_0) $  with $s_1<S$ (all depending on $\gote$)  and  
 $ c(S) >  0 $ such that for all $N_0\ge \ol{N}_0 $, if the smallness condition 
\begin{equation}\label{epN0}
 \e  N_0^{S} < c(S)
 \end{equation} 
holds, then  there exists   
a function $ u_\e \in C^{1}(\mathcal I, H^{s_{1}+\nu})$  with
$u_0(\la)=0$, which solves
\begin{equation}\label{equzero}
F(\e,\la,u_\e (\la))=0 
\end{equation}
 for all $ \la\in \CCCC_{\e} \subset {\mathcal I}  $ defined in \eqref{Cantorfinale} below. 
The set $\CCCC_\e$ is defined  in terms of the {\em ``solution"} $u_\e (\la)$, as
\begin{equation}\label{Cantorfinale}
 \CCCC_\e:= \bigcap_{n \geq 0}\bar \calG^0_{N_0^{2^n}} \cap\bar \gotG_{N_0^{2^n}}\cap \tilde{\mathcal I} 
 \cap \ol{\mathcal I}
\end{equation}
where $  \tilde{\mathcal I}=\tilde{\mathcal I}(N_0)$
 is defined in Hypothesis \ref{hyp.separacatene}, $ \ol{\mathcal I} $
in \eqref{lambdabarra}, 
and, 
for all $ N \in \NNN $, 
 \begin{equation}\label{buoninormal2finale}
\bar \gotG_{N}:= \Big\{\la\in\mathcal I\;:\; \|L^{-1}_{N}(\e,\la,u_\e (\lambda))\|_{0}\le
 N^{\tau_1}/2 \Big\} \, , 
 \end{equation}
 \begin{equation}\label{buoniautovaloriinf}
 \begin{aligned}
 \bar \calG^0_{N}:=\Big\{& \la\in\mathcal I\;:\;
 \forall\; j_{0}\in \SSSS  \mbox{ there is a covering } \\
 \bar B^0_{N}&(j_{0},\e,\la)\subset \bigcup_{q=1}^{N^{\gote}}I_{q},
 \mbox{ with }I_{q}=I_{q}(j_{0})\mbox{ intervals with }
 \meas(I_{q})\le  N^{-\tau_1}
 \Big\}
  \end{aligned}
 \end{equation}
 with 
  \begin{equation}\label{tetacattiviautovalorifinali}
\bar B^0_{N}(j_{0},\e,\la):= \Big\{\theta\in\RRR\;:\;
\|L_{N,j_{0}}^{-1}(\e,\la,\theta, u_\e(\lambda))\|_0>N^{\tau_1}/2\Big\} \, . 
 \end{equation}
Finally,  if the tame estimates (f1)-(f2), \eqref{hyp.normaT}, \eqref{lip}
hold up to $ S' = + \infty $ then 
$ u_\e (\la) \in \cap_{s \geq 0} H^s $. 
\end{theo}
 
 In  applications, it is often useful to work in appropriate closed subspaces 
 $ \widehat H^s(\gotK) \subset H^s(\gotK)$ which are invariant under the action of $F$. 
 The following corollary holds:
 \begin{coro}\label{coromerd}
 Assume, in addition to the hypotheses of Theorem \ref{principe}, 
 that 
 $ F(\e,\la, \cdot ): \widehat H^{s+\nu}(\gotK) \to \widehat H^s(\gotK)$, $ \forall s > s_0 $. 
 Then the function $u_\e$ provided by Theorem \ref{principe} belongs to $C^1(\calI,\widehat H^{s_1+\nu}(\gotK))$.
\end{coro}

 In  Theorem \ref{principe}  
 the Cantor like  $\CCCC_\e$ defined in \eqref{Cantorfinale} may be empty.
 In order to prove that  it has asymptotically full  measure we need  more informations. 
We fix $ N_0 = [\e^{-1/(S+1)}] $ so that the smallness condition \eqref{epN0} is satisfied for $ \e $ small enough.

\begin{theo}\label{mamma}
Let $ N_0 = [\e^{-1/(S+1)}] $ with $ \e $ small enough so that \eqref{epN0} holds. 
Assume, in addition to the hypotheses of Theorem \ref{principe}, that
for all $N\ge N_0 $, 
\begin{equation}\label{meas.bad1}
{\rm meas}({\mathcal I} \setminus {\bar \calG}_{N}^{0} ), {\rm meas}({\mathcal I} \setminus {\bar \gotG}_{N} ) 
= O( N^{-1} )  \, , \quad 
 {\rm meas}({\mathcal I} \setminus (\ol{\mathcal I}\cap\tilde{\mathcal I})) = O (N_0^{- 1 }) \, . 
\end{equation}
 Then 
 $\mathcal C_\e$ satisfies, for some $ \mathtt K > 0 $,   
\begin{equation}\label{misurafin}
{\meas}({\mathcal I} \setminus \mathcal C_\e ) \le \mathtt K \e^{1/(S+1)}.
\end{equation}
 \end{theo}
 
\prova
Let us denote $N_n=N_0^{2^n}$.
By the explicit expression \eqref{Cantorfinale} we have
\begin{equation}\label{asintotica.nls}
\begin{aligned}
\meas(\mathcal I\setminus\CCCC_{\e})& = \meas\Big(
\bigcup_{n\ge 0}(\bar \calG^{0}_{N_{n}})^{c} \bigcup_{n \ge 0}(\bar \gotG_{N_{n}} )^{c}
\cup\tilde{\mathcal I}^{c}\cup\ol{\mathcal I}^c\Big)\\
&\le \sum_{n\ge0}\meas ({\mathcal I} \setminus \bar \calG^{0}_{N_{n}})  +\sum_{n\ge 0}\meas(  {\mathcal I} \setminus  
\bar \gotG_{N_{n}})
+\meas( {\mathcal I} \setminus (\ol{\mathcal I} \cap \tilde {\mathcal I}))\\
&\stackrel{\eqref{meas.bad1}}{\le}
C_0\sum_{n\ge0}N_n^{- 1}+  C_1 N_{0}^{- 1} 
\le C' N_{0}^{-1}  {\le} \mathtt K \e^{1/(S+1)}
\end{aligned}
\end{equation}
which proves \eqref{misurafin}.
\EP

In the applications to NLW and NLS the conditions \eqref{meas.bad1} will be verified taking $ \tau_0 $, $ \tau_1 $ large,  
with a suitable $ \gote $, see Proposition \ref{measure.nlw}. 

\zerarcounters 
\section{Applications to PDEs}\label{section 3}

Now we  apply Theorems \ref{principe}-\ref{mamma} to the NLW and NLS equations
\eqref{nlw}-\eqref{nls}.  To be precise, when $\mathtt M$ is a  manifold which is homogeneous with respect to 
a compact Lie group, we rely on Corollary \ref{coromerd}.

We briefly  recall  the relevant properties of harmonic analysis on compact Lie groups that we need, 
referring to  \cite{Pro} (and \cite{BP}) 
for precise statements and proofs.

A compact manifold  $ \mathtt M $ which is homogeneous with respect to a compact Lie
group is, up to an isomorphism,  diffeomorphic to 
\begin{equation}\label{MGN}
\mathtt M =  G/ N \, , \quad G := \Gacca\times\TTT^{r_2} \, , 
\end{equation}
where $ \Gacca $ is a simply connected compact Lie group, $ \TTT^{r_2} $ is a torus,  and $ N $ is a closed 
subgroup of $  G  $. Then, a function on $ \mathtt M  $ can be seen 
as a function defined on $ G $ which is invariant under the action of $ N $, and  
the space $H^s(\mathtt M,\CCC)$ (or $H^s(\mathtt M,\RRR)$) can be identified 
with the subspace 
\begin{equation}\label{subsp}
\widehat H^s := 
\widehat H^s(G ,\CCC):= 
\Big\{ u \in H^s(G)\,:\; u(x)= u(x g )\,,\ \ \forall  x \in G = \Gacca\times\TTT^{r_2}, g \in N \Big\}. 
\end{equation}
Moreover, the Laplace-Beltrami operator on  $ \mathtt M $ can be identified with the Laplace-Beltrami operator 
on the Lie group $ G  $, acting on functions invariant under $ N $ (see Theorem 2.7, \cite{BP}). 
Then we ``lift" the equations  \eqref{nlw}-\eqref{nls} on  $ G $ and we use harmonic analysis on Lie groups.

\subsection{Analysis on  Lie groups}\label{liegroups}

Any simply connected compact Lie group $\Gacca$ is the product of a finite number of
simply connected 
Lie groups of simple type (which are classified and come in a finite number of families).

Let $\Gacca$ be of simple type, with dimension $\gotd$ and rank $r$.
Denote by ${\mathtt w}_{1},\ldots,{\mathtt w}_{r}\in \RRR^r$ the fundamental weights of $\Gacca$
and consider the cone of dominant  weights
$$
\SSSSGa:=\Big\{j=\sum_{p=1}^{r}j_{p}{\mathtt w}_{p}\,:\,
j_{p}\in\NNN \Big\} \subset \Lambda := \Big\{j=\sum_{p=1}^{r}j_{p}{\mathtt w}_{p}\,:\,
j_{p}\in \ZZZ \Big\} \, . 
$$
Note that $\SSSSGa$ satisfies \eqref{nonsodimeglio} and indexes the finite dimensional irreducible representations of $\Gacca$.

The eigenvalues and the eigenfunctions of the
Laplace-Beltrami operator $\Delta$ on $\Gacca$ are
\begin{equation}\label{autov.lapla}
\mu_{j}:=-|j+\rho|_2^{2}+|\rho|_2^{2},
\qquad
\ff_{j,\s}(x),
\quad
x\in \Gacca,
\quad
j\in\SSSSGa,
\quad
\s=1,\ldots,d_{j},
\end{equation}
where 
$\rho:=\sum_{i=1}^{r}{\mathtt w}_{i} $,  $ |\cdot|_2$ denotes the euclidean norm on $\RRR^{r}$, 
and $\ff_{j}(x)$ is the (unitary) matrix associated to  
an irreducible unitary representation  $(R_{V_{j}},V_{j})$  of $\Gacca $, precisely 
\begin{equation}\nonumber
(\ff_{j}(x))_{h,k}=\langle R_{V_{j}}(x)v_{h},v_{k}\rangle,
\qquad v_{h},v_{k} \in V_{j} \, , 
\end{equation}
where $ (v_h)_{h = 1, \ldots, {\rm dim} V_j} $
is an orthonormal basis of the finite dimensional euclidean space $ V_j $ with scalar product
$ \langle \cdot, \cdot \rangle $.
We denote by  $\mathcal N_j$  the corresponding eigenspace of $ \Delta $.
The \emph{degeneracy of the
eigenvalue} $\mu_{j}$ satisfies
\begin{equation}\label{pol-bound}
d_{j}\le |j+\rho|_2^{\gotd-r} \, . 
\end{equation}

The Peter-Weyl theorem implies the orthogonal decomposition
$$
L^{2} (\Gacca)=\bigoplus_{j\in \SSSSGa} 
\NN_{j} \, . 
$$
Many informations on the eigenvalues $ \mu_j $ are known. 
There exists an integer $\gotZ\in\NNN $ such that  (see \cite{BP}-Lemma 2.6)  
the fundamental weights satisfy 
\begin{equation}\label{Cartan}
w_i  \cdot w_p \in \gotZ^{-1} \ZZZ \, , \quad \forall i,p =1, \ldots, r \, , 
\end{equation}
so that, in particular,  
\begin{equation}\label{stok}
\mu_{j}:=-|j+\rho|_2^{2}+|\rho|_2^{2} \, , \quad  
 |j|_2^2 \, , \quad \rho\cdot j,\quad |\rho|^2_2 \in \gotZ^{-1}\ZZZ,\quad \forall j \in\SSSSGa \, .
\end{equation}
For a product group, $ L^2 (\Gacca_1 \times \Gacca_2) = L^2 (\Gacca_1) \otimes L^2 ( \Gacca_2) $ and 
all the irreducible representations are obtained by the tensor product of the irreducible representations of 
$ \Gacca_1 $ and $ \Gacca_2 $. Hence we extend all the above properties to any compact
Lie group $ \Gacca $. 
For simplicity we still denote the dimension of the group as $\gotd$ and the rank as $r$.
In particular   $\Lambda_+(G)=\Lambda_+(\Gacca)\times\ZZZ^{r_2}$ (see \eqref{MGN}) is the index set 
for the irreducible representations of $ G $,  
with indices $j=(j_1,j_2)$, $ j_1 \in \Lambda_+(\Gacca) $, $ j_2 \in \ZZZ^{r_2} $,  
and $\rho\rightsquigarrow (\rho,0)$.  

We denote the indices $i=(l,j)\in\ZZZ^{d}\times \Lambda_+ (G)$, 
so that $L^2(\TTT^d\times\Gacca\times\TTT^{r_2})$ naturally decomposes as product of subspaces $\NN_k$
of the form
$$
\NN_{k}:=\av{{\rm e}^{\ii \f\cdot l}}\otimes\NN_{j}=\av{{\rm e}^{\ii \f\cdot l}}\otimes\NN_{j_1}\otimes
\av{{\rm e}^{\ii x_2\cdot j_2}}.
$$
We also set
$$
|i|:=\max(|l|,|j|) \,,\quad |l|:=|l|_\io\,,\quad |j|:=|j|_\infty= \max_{i} |j_i|\,,\quad \laa i\raa := \max(1,|i|) \, .
$$
The Sobolev spaces $H^s(\TTT^d\times G) $ and $ H^s(\TTT^d\times G)\times H^s(\TTT^d\times G)$,
 for a Lie group $G$, 
 can be now identified with sequence spaces introduced  in Section \ref{hilbertspaces}. 
\begin{example}\label{zoppa}
Let $\gotA := \{ 1 \} $, $\SSSS :=  \SSSSG $ be the cone
of fundamental weights and,  for $k=(l,j)\in \gotK = \ZZZ^d \times \Lambda_+ $, let  $ w_k  := \sqrt{ |l|_2^2 + |j+\rho|_2^2}$. 
Then we may identify $H^s(\gotK)$ with the Sobolev space $H^s(\TTT^d\times G)$.
\end{example}

\begin{example}\label{zoppa2}
Let $\gotA := \{1,-1\}$, $\SSSS := \SSSSG $ be the cone of fundamental
weights and,  for $k=(l,j,\gota)\in \gotK := \ZZZ^d \times \Lambda_+ \times \gotA $, let  
$ w_k := \sqrt{ |l|_2^2 + |j+\rho|_2^2}$. 
Then we may identify $H^s(\gotK)$ with the Sobolev space $H^s(\TTT^d\times 
G)\times H^s(\TTT^d\times G)$.
\end{example}

The final fundamental property that we exploit concerns the off-diagonal decay of the 
block matrix which represents the multiplication operator, see \eqref{lematrici}.  
The block structure of this matrix takes into account the (large) multiplicity of the degenerate eigenvalues of $ \Delta $. 
We remark that several blocks could correspond to the same eigenvalue (as in the case of the torus).
The next lemmas, proved in \cite{BP}, 
are ultimately connected to the fact that 
the product of two eigenfunctions of the Laplace operator on 
a Lie group is a finite sum of eigenfunctions, see \cite{BP}-Theorem 2.10.

The forthcoming Lemmas are a reformulation of Proposition 2.19 and Lemma 7.1 in \cite{BP}
respectively, and they require the polynomial bound \eqref{pol-bound}.

\begin{lemma}\label{mult.oper} {(\cite{BP}-Proposition 2.19)}
Let $\mathfrak K$ be as in Example \ref{zoppa} 
and $b\in H^{s}(\TTT^d\times G)$ be real valued. Then the multiplication operator
$ B:u(\f,x)\mapsto b(\f,x)u(\f,x) $
is self-adjont in $L^{2}$ and, for any $s>(d+\gotd)/2$, 
$$
\|B_{k}^{k'}\|_{0}\le\frac{C(s)\|b\|_{s}}{\av{k-k'}^{s-(d+\gotd)/2}}, \quad \forall k , k' \in \ZZZ^d \times \Lambda_+ \, , 
$$
where $ B_{k}^{k'} \in {\rm Mat}(d_{j}\times d_{j'},\CCC) $, see \eqref{lematrici}.
\end{lemma}

\begin{lemma}{(\cite{BP}-Lemma 7.1)}\label{moltiplicazione.matrici}
Let $\mathfrak K$ be as in Example \ref{zoppa2}. 
Consider $a,b,c\in H^{s}(\TTT^d\times G) $ with $a,b$ real valued. Then the multiplication operator with matrix
$$
B= \begin{pmatrix}
a(\f,x) & c(\f,x) \\ \bar c(\f,x)& b(\f,x)
\end{pmatrix}
$$
is self-adjont in $L^{2}$ and, for any $s>(d+\gotd)/2$, 
$$
\|B_{\{i\}}^{\{i'\}}\|_0 \le C(s)\frac{\max(\|a\|_{s},\|b\|_{s},\|c\|_{s})}{\av{i-i'}^{s-(d+\gotd)/2}} \, , 
\quad \forall i, i' \in \ZZZ^d \times \Lambda_+ \, . 
$$
\end{lemma}

\begin{coro}\label{nuzero}
Let $B$ be a linear operator as in the previous two Lemmas. Then, for all $ s>(d+\gotd)/2  $, 
$$
\bvert B\bvert_s\leq C(s) \max(\|a\|_{s+\nu_0},\|b\|_{s+\nu_0},\|c\|_{s+\nu_0})\,,\qquad \nu_0:= (2d+\gotd +r+1)/2 \, . 
$$ 
\end{coro}

\subsection{Proof of Theorem \ref{teoremone} for NLW} 
\label{sec:nlw} 

We apply Theorems \ref{principe}-\ref{mamma} to the operator
$$
\begin{aligned}
 F(\e,\la,\cdot):  H^{s+2}(\TTT^d\times\mathtt M,\RRR)  &\longrightarrow \, H^{s}(\TTT^d\times\mathtt M,\RRR) \\  
 u \; &\longmapsto (\la\ol\om\cdot\del_\f)^2u-\Delta u + m u - \e f(\f,x,u) 
 \end{aligned}
$$ 
which can be extended to  
$H^{s+2}(\TTT^d\times\Gacca\times\TTT^{r_2},\RRR)  \to \, H^{s}(\TTT^d\times\Gacca\times\TTT^{r_2},\RRR)$ such that for all 
$u\in  H^{s+2} (\TTT^d) \otimes \widehat H^{s+2}  $ one has $F(\e,\lambda, u)\in H^s (\TTT^d)  \otimes \widehat H^{s}$
where $ \widehat H^{s}$ is defined in \eqref{subsp}.

Setting $\gotA :=\{1\}$, $\SSSS:=\SSSSG$, $ G := \Gacca\times\TTT^{r_2} $, 
we are in the functional setting  of Example \ref{zoppa}.
The Hypothesis  \eqref{Hdiag}-\eqref{Ddiag} holds with  $\nu=2$ and the interpolation estimates
(f1)-(f2) are verified provided that $f(\f,x,u) \in C^q $ for $ q $  large enough  and $s_0 > (d + \gotd)/2\ge(d+r)/2$.

\begin{rmk}\label{dimensioni}
We require $s_0> (d+\gotd)/2$ in view of the embedding $H^{s_0}(\TTT^d\times G)
\hookrightarrow L^{\io}(\TTT^d\times G)$ which, in turn, implies the algebra and interpolation
properties of the spaces $H^s(\TTT^d\times G)$, $ s \geq s_0 $. The weaker bound $s_0>(d+r)/2$ is
sufficient
in order to prove the algebra and interpolation properties of the decay norm
$\bvert\cdot \bvert_s$ (see Section \ref{sub.linearop}), which hold with no
constraint on the multiplicity $d_j$.
\end{rmk}

The linearised operator
$$
D (\la) -\e g(\f,x), \quad 
D (\la) :=(\la \bar  \oo\cdot\partial_{\f})^{2}-\Delta+m,
\quad
g(\f,x):=(\partial_{u}f)(\f,x,u(\f,x)) \, , 
$$
is represented, in the Fourier basis ${\rm e}^{\ii l\cdot\f}\ff_{j}(x)$,  as in 
\eqref{operatore}
with
\begin{equation}\label{DkNLW}
D_k (\lambda) := D_{(l,j)} (\lambda) = -(\la\ol{\oo}\cdot l)^{2}+ m - \mu_{j}
\end{equation}
and $T(u) $ is the matrix associated to the multiplication operator by $-g(\f,x)$.
Corollary \ref{nuzero} implies that 
$T(u) \in \mathcal M^{s-\nu_0}$ for
all $u\in H^s(\TTT^d\times G)$ and the estimates \eqref{hyp.normaT}, \eqref{lip} hold
by interpolation.
Hypothesis \ref{matricione} holds with $\gotD_{j}(y)=  -y^2 + m -\mu_j $.

Also Hypothesis \ref{hyp.misura.diag} holds:
a direct computation shows that 
$$
\{\theta \in \RRR\;:\; |D_{(l,j)} (\la, \theta) |\le N^{-\tau_{1}} \} \subseteq \bigcup_{q=1,2} I_q \, , \ \  
I_q \mbox{ intervals with } {\rm meas} (I_q) \leq \frac{4N^{-\tau_{1}} }{\sqrt{m-\mu_j}}+ O(N^{-2\tau_{1}}),
$$
and Hypothesis \ref{hyp.misura.diag} holds with $ \gotn = 16 /\sqrt{m} $. 

\smallskip

Hypothesis \ref{hyp.separacatene}
about the length of  chains of  singular sites follows as 
in \cite{BB2} because the eigenvalues of the Laplace-Beltrami operator 
are very  similar to those on a torus, see  \eqref{stok}. 
For $\g > 0$ let 
\begin{equation}\label{lambdatilte}
\begin{aligned}
\tilde{\mathcal I} := \tilde{\mathcal I}(\g)
:= &\Big\{\lambda\in[1/2,3/2]\;:\;   |P(\la\ol{\oo})|\ge\frac{\g}{1+|p|^{d(d+1)}}\; ,  \forall \, \mbox{non zero  polynomial} \\
&\; P(X)\in
\ZZZ[X_{1},\ldots,X_{d}]   \mbox{ of the form }
P(X)=p_{0}+\sum_{1\le i_{1},i_{2}\le d}p_{i_{1},i_{2}}X_{i_{1}}X_{i_{2}} \Big\} \, . 
\end{aligned}
\end{equation}

\begin{lemma}\label{catenecattive}
For all $N_0\ge2 $, Hypothesis \ref{hyp.separacatene} is satisfied with $\tilde{\mathcal I}$
defined in \eqref{lambdatilte} and $\g=N_0^{-1}$.
\end{lemma}

\prova
The proof follows  Lemma 4.2 of \cite{BB2}. 
First of all, it is sufficient to bound the length of a
$\Gamma$-chain of singular sites for $D(\la,0)$. 
Then we consider the quadratic form 
\begin{equation}\label{quadratic-form}
Q:\RRR\times\RRR^{r}\to\RRR \, , \quad  Q(x,j):=-x^{2}+|j|^{2}_2,
\end{equation}
and the associated bilinear form  $ \Phi=- \Phi_{1}+\Phi_{2} $ where
\begin{equation}\label{split.b}
\Phi_{1}((x,j),(x',j')):=xx', \qquad
\Phi_{2}((x,j),(x',j')):=j\cdot j' \, .
\end{equation}
For  a $\Gamma$-chain of sites $\{k_{q}=(l_{q},j_{q})\}_{q=0, \ldots, \ell}$ which are singular
for $D(\la,0)$ (Definition \ref{regular})
we have, recalling  \eqref{DkNLW}, \eqref{stok}, and setting $x_{q} := \oo\cdot l_{q}$,   
$$
|Q(x_{q},j_{q}+\rho)|<2+|m-|\rho|^{2}_2|, \qquad \forall q = 0,\ldots,\ell \, .
$$
Moreover, by \eqref{quadratic-form}, \eqref{Gamma-chain},  we derive
$|Q(x_{q}-x_{q_{0}},j_{q}-j_{q_0})|\le C|q-q_{0}|^{2}\Gamma^{2} $, $ \forall 0 \le q, q_0 \le \ell $, and so
\begin{equation}\label{stimameglio}
|\Phi((x_{q_{0}},j_{q_{0}}+\rho),(x_{q}-x_{q_{0}},j_{q}-j_{q_{0}}))|
\le C' |q-q_{0}|^{2}\Gamma^{2} \, . 
\end{equation}
Now 
we introduce the subspace of
$\RRR^{1+r}$ given by
$$
{\mathcal S}:={\rm Span}_{\RRR}\{(x_{q}-x_{q_{0}},j_{q}-j_{q_{0}})\;:
\;q=0,\ldots,\ell\}
$$
and denote by $\gots\le r+1$ the dimension of ${\mathcal S}$.
Let $\de >  0 $ be  a small parameter specified later on. We distinguish two cases. 

\noindent
{\bf Case 1.} {\it For all $q_{0}=0,\ldots,\ell$ one has
\begin{equation}\label{sottospazio.nlw-caso1}
{\rm Span}_{\RRR}\{(x_{q}-x_{q_{0}},j_{q}-j_{q_{0}})\;:
\;|q-q_{0}|\le \ell^{\de},\;q=0,\ldots,\ell\}={\mathcal S}.
\end{equation}}
In such a case, we select a basis 
$ f_b :=(x_{q_{b}}-x_{q_{0}},j_{q_{b}}-j_{q_{0}})=(\oo\cdot\Delta l_{q_{b}},
\Delta j_{q_{b}})$, $b=1,\ldots,\gots$ of ${\mathcal S}$, where 
$\Delta k_{q_{b}}=(\Delta l_{q_{b}}, \Delta j_{q_{b}})$ satisfies
$\bvert \Delta k_{q_{b}}\bvert \le C\Gamma|q_{b}-q_{0}|\le C\Gamma \ell^{\de}$. Hence
we have the bound
\begin{equation}\label{normabase}
\bvert f_{q_{b}}\bvert\le C\Gamma \ell^{\de}, \qquad
b=1,\ldots,\gots.
\end{equation}
Introduce also the matrix $\Omega=(\Omega^{b'}_{b})_{b,b'=1}^{\gots}$ with 
$\Omega^{b'}_{b}:=\Phi(f_{b'},f_{b})$, 
that, according to \eqref{split.b}, we  write
\begin{equation}\label{chieomegone}
\Omega=\Bigl(- \Phi_{1}(f_{b'},f_{b})+\Phi_{2}(f_{b'},f_{b})
\Bigr)_{b,b'=1}^{\gots}=-X+Y,
\end{equation}
where $X^{b'}_{b} :=(\oo\cdot\Delta l_{q_{b'}})(\oo\cdot\Delta l_{q_{b}})$ and
$Y^{b'}_{b} :=(\Delta j_{q_b'})\cdot (\Delta {j_{q_b}})$.
%
By \eqref{Cartan} the matrix $Y$ has entries in $\gotZ^{-1}\ZZZ$  and  the matrix
$X$ has rank $1$ since each  column is
\begin{equation}\nonumber
X^{b}=(\oo\cdot\Delta l_{q_{b}})
\begin{pmatrix}
\oo\cdot \Delta l_{q_{1}} \cr
\vdots \cr
\oo\cdot \Delta l_{q_{\gots}}\end{pmatrix}, \quad
b=1,\ldots,\gots.
\end{equation}
Then, 
since the determinant of a matrix with two collinear columns $X^{b},X^{b'}$, $b\ne b' $, is zero, we get
$$
\begin{aligned}
P(\oo):&=\gotZ^{r+1}{\rm det}(\Omega)=\gotZ^{r+1}{\rm det}(-X+Y)\\
&=\gotZ^{r+1}(\det(Y)-\det(X^{1},Y^{2},\ldots,Y^{\gots})-\ldots-
\det(Y^{1},\ldots,Y^{\gots-1},X^{\gots}))
\end{aligned}
$$
which is a quadratic polinomial as in \eqref{lambdatilte} with coefficients
$\le C(\Gamma \ell^{\de})^{2(r+1)}$. Note that $P\not\equiv0$.
Indeed, if $P\equiv0$ then
\begin{equation}\nonumber
0=P(\ii\oo)=\gotZ^{r+1}\det(X+Y)=\gotZ^{r+1}\det(f_{b}\cdot f_{b'})_{b,b'=1,\ldots,\gots} > 0 
\end{equation}
because $\{f_{b}\}_{b=1}^{\gots}$ is a basis of ${\mathcal S}$. This contradiction proves that $ P \neq 0 $. 
But then,  by \eqref{lambdatilte}, 
$$
\gotZ^{r+1}|\det(\Omega)|=|P(\oo)|\ge\frac{\g}{1+|p|^{d(d+1)}}\ge
\frac{\g}{(\Gamma \ell^{\de})^{C(d,r)}} \, ,
$$
 the matrix $\Omega$ is invertible and 
\begin{equation}\label{inversaomegone}
|(\Omega^{-1})^{b'}_{b}|\le C\gamma^{-1}(\Gamma \ell^{\de})^{C'(d,r)}.
\end{equation}
Now let 
$  {\mathcal S}^{\perp}:= {\mathcal S}^{\perp\Phi} := \{v\in\RRR^{r+1}\;:\;\Phi(v,f)=0,\;\forall\,f\in
{\mathcal S}\}$. Since $\Omega$ is invertible, the quadratic form
$\Phi_{{\mathcal S}}$ is non-degenerate and so
$\RRR^{r+1}={\mathcal S}\oplus{\mathcal S}^{\perp}$.
We denote  $\Pi_{{\mathcal S}}:\RRR^{r+1}\to{\mathcal S}$ the 
projector onto ${\mathcal S}$. Writing 
\begin{equation} \label{proj.esse}
\Pi_{{\mathcal S}}(x_{q_{0}},j_{q_{0}}+\rho)=\sum_{b'=1}^{r+1}a_{b'}f_{b'} \, ,
\end{equation}
and since 
$f_{b}\in{\mathcal S}$, $ \forall b = 1,\ldots,\gots$, we get 
\begin{equation}\nonumber
w_b:=\Phi\big((x_{q_{0}},j_{q_{0}}+\rho),f_{b}\big)=
\sum_{b'=1}^{\gots}a_{b'}\Phi(f_{b'},f_{b})=\sum_{b'=1}^\gots\Omega_b^{b'} a_{b'}
\end{equation}
where $\Omega$ is defined in \eqref{chieomegone}.
The definition
of $f_{b}$, the bound \eqref{stimameglio} and \eqref{sottospazio.nlw-caso1} imply
$|w|\le C(\Gamma \ell^{\de})^{2}$. Hence, by \eqref{inversaomegone},
we deduce $|a|=|\Omega^{-1}w|\le C' \gamma^{-1}(\Gamma \ell^{\de})^{C(r,d)+2}$, whence, by
\eqref{proj.esse} and \eqref{normabase}, 
$$
|\Pi_{{\mathcal S}}(x_{q_{0}},j_{q_{0}}+\rho)|\le \gamma^{-1}(\Gamma \ell^{\de})^{C'(r,d)}.
$$
Therefore, for any $q_{1},q_{2}=0,\ldots,\ell $, one has
$$
|(x_{q_{1}},j_{q_{1}})-(x_{q_{2}},j_{q_{2}})|
=|\Pi_{{\mathcal S}}(x_{q_{1}},j_{q_{1}}+\rho)-
\Pi_{{\mathcal S}}(x_{q_{2}},j_{q_{2}}+\rho)|\le \gamma^{-1}(\Gamma \ell^{\de})^{C_{1}(r,d)},
$$
which in turn implies $|j_{q_{1}}-j_{q_{2}}|\le \gamma^{-1}(\Gamma \ell^{\de})^{C_{1}(r,d)}$
for all $q_{1},q_{2}=0,\ldots,\ell$. Since all the $j_{q}$ have $r$ components
(being elements of $\SSSSG$) they are at most
$C\gamma^{-r}(\Gamma \ell^{\de})^{C_{1}(r,d)r}$. We are considering a $ \Gamma $-chain 
in $ \Sigma_K $ (see Definition \ref{fibre}) and so,
for each $q_{0} $, the number of $q \in \{0, \ldots, \ell \} $ such that $j_{q}=j_{q_{0}}$
is at most $K$ and hence
$$
\ell\le \gamma^{-r}(\Gamma \ell^{\de})^{C_{2}(r,d)}K \leq (\Gamma K)^r  (\Gamma \ell^{\de})^{C_{2}(r,d)}K \leq
 \ell^{\de C_{2}(r,d)} (\Gamma K)^{r +C_{2}(r,d)} 
$$
because of the condition $\Gamma K \geq N_0 $ (Hypothesis \ref{hyp.separacatene}) and
$ N_0 = \gamma^{-1} $.
Choosing $\de<1/(2C_{2}(r,d))$ we  get $\ell \leq (\Gamma K)^{2(r +C_{2}(r,d))}$. 

\noindent
{\bf Case 2.} There is $q_{0}=0,\ldots,\ell$ such that
\begin{equation}\nonumber
{\rm dim}(
{\rm Span}_{\RRR}\{(x_{q}-x_{q_{0}},j_{q}-j_{q_{0}})\;:
\;|q-q_{0}|\le \ell^{\de},\;q=0,\ldots,\ell\})\le \gots-1.
\end{equation}
Then we  repeat the argument of Case 1 for the sub-chain
$\{(l_{q},j_{q})\;:\:|q-q_{0}|\le \ell^{\de}\}$ and obtain a bound for $\ell^\delta $.
Since this procedure should be applied at most $r+1$ times, at the end
we get a bound like $\ell\le(\Gamma K)^{C_{3}(r,d)}$.
\EP

We have verified the hypotheses of  Theorem \ref{principe} and Corollary \ref{coromerd}.  
The next proposition proves that also the assumptions \eqref{meas.bad1} in Theorem \ref{mamma} hold.

\begin{prop}\label{measure.nlw}
Fix  $ \tau_0 > r + 3d + 1$, $\tau_1 > d + \gotd + 2 $ and $\gote = d+\gotd+r+4$. There exists $ N_0 \in \NNN $ (possibly larger than
the $N_0 $ found in Theorem \ref{principe}) such that  \eqref{meas.bad1} holds. 
\end{prop}

The proof of Proposition \ref{measure.nlw} -which will continue until the end of this section- follows by  basic
properties of the eigenvalues of a self-adjoint operator, which are a
consequence of their variational characterisation. 
Proposition \ref{measure.nlw} is indeed a reformulation of Proposition  5.1 of \cite{BB2}.
With respect to \cite{BB2}, the eigenvalues of the Laplacian in \eqref{stok} are different and
the index set $\Lambda$  is not an orthonormal lattice.

 \begin{rmk}\label{pippon}
  There are two positive constants $\pippa< \Pippa $ such that 
 $\pippa |j|\leq |j|_2 \leq \Pippa |j|   $. Hence if $|j|> \alpha \pippa^{-1}N$, $ N > 2 |\rho|_2 $, then the eigenvalues  $ \mu_j  $ in 
\eqref{stok}   satisfy 
 $
 -\mu_j > \alpha(\alpha-1)N^2.
 $ On the other hand if $|j|\leq \alpha \pippa^{-1}N$ then $-\mu_j\leq \alpha(\alpha+1)(\Pippa/\pippa)^{2}N^2 $. 
 \end{rmk}
 
 Recall that if $A,A'$ are self-adjoint matrices, then their eigenvalues $ \mu_p (A) $, $ \mu_p (A') $ (ranked in nondecreasing order) 
satisfy
\begin{equation}\label{eq:7.10}
|\m_p (A) - \m_p (A') |\le \|A-A'\|_{0} \, .
\end{equation}
We study finite dimensional restrictions of the the self-adjoint operator $L(\e,\la)=L(\e,\la,u_\e(\la)) = D(\lambda) + \e  T(\e, \lambda) $.

One proceeds differently for $|j_{0}| \ge ( \pippa + 5 ) \pippa^{-1}N$ and $|j_{0}|< ( \pippa + 5 ) \pippa^{-1}N$. 
We assume $N\ge N_0>0$ large enough and $\e \|T\|_0\le 1$.

\begin{lemma}\label{primapartizione} 
For all $j_{0}\in\SSSSG$, $| j_{0}|\ge ( \pippa + 5 ) \pippa^{-1}N$, 
and for all $\la\in[1/2,3/2]$ one has
$$
B_{N}^{0}(j_{0},\e,\la)\subset\bigcup_{q=1}^{N^{d+\gotd+2}}I_{q}\, , 
\  \mbox { with } I_{q} = I_q(j_0)  \mbox{ intervals with } \meas(I_q) \le N^{-\tau_1} \, .
$$
\end{lemma}

\prova
We first show that $B_{N}^{0}(j_{0},\e,\la)\subset\RRR\setminus[-2N,2N]$. Indeed by
\eqref{eq:7.10} all the eigenvalues $\mu_{l,j, \sigma}(\theta)$, $ \sigma = 1, \ldots, d_j $, of $L_{N,j_{0}}(\e,\la,\theta)$ (recall
that $ d_j $ denotes the degeneracy of the eigenvalues $ \mu_j $ in \eqref{autov.lapla}), 
are of the form 
\begin{equation}\label{diff auto}
\mu_{l,j, \sigma}(\theta)
=\de_{l,j}(\theta)+O(\e\|T\|_{0}), 
\quad
\de_{l,j}(\theta) :=-(\oo\cdot l+\theta)^{2}+m-\mu_{j} \, .
\end{equation}
Since $|\oo|_1=\la|\ol{\oo}|_1\le3/2$, $|j-j_{0}|\le N$, $|l|\le N$, one has, by Remark \ref{pippon},  
$$
\de_{l,j}(\theta)\ge-\Big(\frac{3}{2}N+|\theta|\Big)^{2}+ 20 N^2 > 7 N^2 \, \, , \quad \forall |\theta | < 2 N \, . 
$$
By \eqref{diff auto} we deduce $ \mu_{l,j, \sigma}(\theta)  \geq 6 N^2  $
and this implies $B_{N}^{0}(j_{0},\e,\la)\cap[-2N,2N]=\emptyset$.
Now set
$ B_{N}^{0,+}:=B_{N}^{0}(j_{0},\e,\la)\cap(2N,+\io) $,  $ B_{N}^{0,-}:=B_{N}^{0}(j_{0},\e,\la)\cap(-\io, - 2 N) $.
Since
$$
\partial_{\theta}L_{N,j_{0}}(\e,\la,\theta)=
\diag_{{\substack|l|\le N,\\ |j-j_{0}|\le N}}-2(\oo\cdot l+\theta) \uno_{d_j} \geq  N\uno,
$$
we apply Lemma 5.1 of \cite{BB2} 
with $\alpha= N^{-\tau_1}$, $\beta= N$ and $|E|\leq C N^{\gotd+d}$ (this is due to the bound  $ d_j \leq |j + \rho |_2^{\gotd-r} $) and obtain
\begin{equation}\nonumber
B_{N}^{0,-}\subset\bigcup_{q=1}^{N^{d+\gotd+1}}I_{q}^{-}\ ,  \quad  I_{q}^{-}=I_{q}^{-}(j_{0}) \mbox{ intervals with }
\meas (I_q) \le  N^{-\tau_1} \, . 
\end{equation}
We can reason in the same way for $B_{N}^{0,+}$ and the lemma
follows.
\EP

Consider now $|j_{0}|< ( \pippa + 5 ) \pippa^{-1}N$. We  obtain a complexity
estimate for $B_{N}^{0}(j_{0},\e,\la)$ by knowing the measure of the set
$$
B^{0}_{2,N}(j_{0},\e,\la):=\left\{\theta\in\RRR\;:\;\|L_{N,j_{0}}^{-1}(\la,
\e,\theta)\|_{0}>N^{\tau_1}/2\right\}.
$$

\begin{lemma}\label{lem.controllobdoppio}
For all $|j_{0}|< ( \pippa + 5 ) \pippa^{-1} N$ and all $\la\in[1/2,3/2]$ one has
$$
B^{0}_{2,N}(j_{0},\e,\la)\subset I_{N}:=[-  g N,  g N] 	, 	\quad  g := (2 \pippa + 8 )\Pippa \pippa^{-1} \, .
$$
\end{lemma}

\prova
If $|\theta|> g  N$ one has
$|\oo\cdot l+\theta|\ge|\theta|-|\oo\cdot l|> (g - (3/2) )N > ( 2 \pippa+6) \Pippa \pippa^{-1} N $.
Using Remark \ref{pippon},  all the eigenvalues 
$$
\mu_{l,j, \sigma}(\theta)=-(\oo\cdot l+\theta)^{2}+m-\mu_{j}+O(\e\|T\|_{0}) \le - (\Pippa \pippa^{-1}N)^{2} \, , \quad \forall |\theta| > g N \, ,
$$
 proving the lemma.   
\EP

\begin{lemma}\label{secondocontrollo}
For all $|j_{0}|\le ( \pippa + 5 ) \pippa^{-1} N$
and all $\la\in[1/2,3/2]$ one has
$$
B^{0}_{N}(j_{0},\e,\la)\subset\bigcup_{q=1}^{\hat{C}\gotM N^{\tau_1+1}}I_{q} \, , 
\ I_{q}=I_{q}(j_{0}) \mbox{ intervals with } \meas (I_q) \le N^{-\tau_1} 
$$
where $\gotM:=\meas(B^{0}_{2,N}(j_{0},\e,\la))$ and $\hat{C}=\hat{C}(r)$.
\end{lemma}

\prova
This is Lemma 5.5 of \cite{BB2}, where our  exponent $\tau_1$ is denoted by $\tau$.
\EP

Lemmas \ref{lem.controllobdoppio} and \ref{secondocontrollo} imply that for all $\la
\in[1/2,3/2]$ the set $B^{0}_{N}(j_{0},\e,\la)$ can be covered by 
$\sim N^{\tau_1+2}$ intervals of length $\le N^{-\tau_1}$. This estimate is not enough. 
Now we prove that for ``most" $\la$ the number of such intervals does not depend on $\tau_1$, showing
that 
$ \meas (B^{0}_{2,N}(j_{0},\e,\la))  = O( N^{\gote - \tau_1}) $ where $ \gote = d+ \gotd + r + 4  $ has been fixed in 
Proposition \ref{measure.nlw}.
To this purpose first  we  provide an estimate for the set
$$
\BBB^{0}_{2,N}(j_{0},\e):=\left\{(\la,\theta)\in[1/2,3/2]\times \RRR\;:\;
\|L_{N,j_{0}}^{-1}(\e,\la,\theta)\|_{0}>N^{\tau_1}/2\right\} \, . 
$$
Then in Lemma \ref{usofubini} we  use Fubini Theorem to obtain the desired bound for $\meas
(B^{0}_{2,N}(j_{0},\e,\la))$.

\begin{lemma}\label{lem.misurabione.nlw}
For all $|j_{0}|< ( \pippa + 5 ) \pippa^{-1} N$ one has
$ {\rm meas}(\BBB^{0}_{2,N}(j_{0},\e))\le C N^{-\tau_1+\gotd+d+1}  $
for some $C >  0 $.
\end{lemma}

\prova
Let us introduce the variables
\begin{equation}\label{variabilinuove}
\pD=\frac{1}{\la^{2}},  \ \h=\frac{\theta}{\la}, \qquad
(\pD,\h)\in[4/9,4]\times [- 2 g N, 2 g N]=: [4/9,4] \times J_N,
\end{equation}
and set
$$
L(\pD,\h):= {\la^{-2}}L_{N,j_{0}}(\e,\la,\theta)=
\diag_{|l|\le N,|j-j_{0}|\le N}\Big(\big(-(\ol{\oo}\cdot l+\h)^{2} +
\pD(-\mu_{j}+m)\big) \uno_{d_j} \Big)+ \e\pD T(\e,1/\sqrt{\pD}).
$$
Note that
\begin{equation}\label{minla}
\min_{j\in\SSSSG} -\mu_{j}+m \ge m.
\end{equation}
Then, except for $(\pD,\h)$ in a set of measure
$O(N^{-\tau_1+\gotd+d+1})$ one has
\begin{equation}\label{stimotta}
\|L(\pD,\h)^{-1}\|_{0}\le N^{\tau_1}/8.
\end{equation}
Indeed
$$
\begin{aligned}
\partial_{\pD}L(\pD,\h)&=\diag_{|l|\le N,|j-j_{0}|\le N}\left((-\mu_{j}+m) \uno_{d_j}\right)
+ \e T(\e,1/\sqrt{\pD}) - \frac{\e}{2} \pD^{-1/2}\partial_{\la}T
\stackrel{\eqref{minla}}  \ge \frac{m}{2}\uno,
\end{aligned}
$$
for $ \e $ small (we used that $\pD\in[4/9,4]$).
Therefore Lemma 5.1 of \cite{BB2} implies that for each $\h$,
the set of $\pD$ such that at least one eigenvalue of $L(\pD,\h)$ has
modulus $\le8 N^{-\tau_1}$, is contained in the union of $O(N^{d+\gotd})$ intervals
with length $O(N^{-\tau_1})$ and hence has measure $\le O(N^{-\tau_1+d+\gotd})$. Integrating
in $\h\in J_{N}$ we obtain \eqref{stimotta} except in a set with measure
$O(N^{-\tau_1+d+\gotd+1})$. 
The same measure estimates hold in the original variables $ (\la, \theta) $ in \eqref{variabilinuove}. 
Finally \eqref{stimotta} implies
$$
\|L^{-1}_{N,j_{0}}(\e,\la,\theta)\|_{0}\le\la^{-2}N^{\tau_1}/8\le N^{\tau_1}/2,
$$
for all $(\la,\theta)\in[1/2,2/3]\times\RRR$ except in a set with
measure $\le O(N^{-\tau_1+d+\gotd+1})$.
\EP

The same argument implies that 
\begin{equation}\label{cattivinormalta}
\meas([1/2,3/2]\setminus \bar \gotG_{N})\le N^{-\tau_1+d+\gotd +1} 
\end{equation}
where $ \bar \gotG_N$ is defined in \eqref{buoninormal2finale}.

Define the set
\begin{equation}\label{FNj0}
\calF_{N}(j_{0}):=\left\{\la\in[1/2,3/2]\;:\;\meas(B_{2,N}^{0}(j_{0},\e,\la))
\ge\hat{C}N^{-\tau_1+d+\gotd +r+2}\right\}
\end{equation}
where $\hat{C}$ is the constant appearing in Lemma \ref{secondocontrollo}.

\begin{lemma}\label{usofubini}
For all $|j_{0}|\le  ( \pippa + 5 ) \pippa^{-1} N$ one has
$ \meas(\calF_{N}(j_{0})) = O ( N^{-r-1}) $. 
\end{lemma}

\prova
By Fubini Theorem we have
$$
\meas(\BBB_{2,N}^{0}(j_{0},\e))=\int_{1/2}^{3/2}{\rm d}\la\,
\meas(B_{2,N}^{0}(j_{0},\e,\la)).
$$
Now, for any $\beta > 0 $, using Lemma \ref{lem.misurabione.nlw} we have
$$
\begin{aligned}
CN^{-\tau_1+d+\gotd+1}&\ge\int_{1/2}^{3/2}{\rm d}\la\,
\meas(B_{2,N}^{0}(j_{0},\e,\la))\\
&\ge\beta \meas(\{\la\in[1/2,3/2]\;:\;\meas(B_{2,N}^{0}(j_{0},\e,\la))\ge\beta \})
\end{aligned} 
$$
and for $\beta = \hat{C}N^{-\tau_1+r+\gotd +d+2}$  we prove the lemma (recall \eqref{FNj0}).
\EP

\begin{lemma}\label{misuralambdabarra} If $ \tau_0 > r+ 3d+ 1 $ then 
$ \meas([1/2,3/2]\setminus\ol{\mathcal I})= O(N_{0}^{-1} )$
where $\ol{\mathcal I}$ is defined in \eqref{lambdabarra}.
\end{lemma}

\prova 
Let us write
$$
[1/2,3/2]\setminus\ol{\mathcal I}=\bigcup_{|l|,|j|\le N_{0}}
\RR_{l,j},\qquad \RR_{l,j}:= \Big\{\la\in\Lambda\;:\; |(\la\ol{\oo}\cdot l)^{2}+\mu_{j}-m|\le N_{0}^{-\tau_{0}} 
\Big\}.
$$
Since $-\mu_{j}+m\ge m>0$, then $\RR_{0,j}=\emptyset$ if $N_{0}>m^{-1/\tau_{0}}$. For $ l \neq 0 $,
 using the Diophantine condition \eqref{dioph}, we get  $ \meas(\RR_{l,j})\le C N_{0}^{-\tau_{0}+2d}, $
so that
$$
\meas([1/2,3/2]\setminus\ol{\mathcal I})\le \sum_{|l|,|j|\le N_{0}}\meas(\RR_{l,j})\le
C N_{0}^{-\tau_{0}+r+3d} = O( N_0^{-1})
$$
because $\tau_{0}-r-3d>1$. 
\EP

The measure of the set $ \tilde{\mathcal I} $ in \eqref{lambdatilte} is estimated in \cite{BB2}-Lemma 6.3 
(where $\tilde{\mathcal I}$ is denoted by $\tilde{\mathcal G}$).

\begin{lemma}[]\label{lem.lambdatilde}
If $\g<\min(1/4,\g_{0}/4)$ (where $ \g_0 $ is that in \eqref{diophquad}) then 
$ \meas([1/2, 3/2] \setminus\tilde{\mathcal I})=O(\g)$.
\end{lemma}

\noindent
{\it Proof of Proposition \ref{measure.nlw} completed.}
By the definition in \eqref{FNj0} 
 for all $\la\not\in\calF_{N}(j_{0})$
one has $\meas(B_{2,N}^{0}(j_{0},\e,\la))<O(N^{-\tau_1+r+\gotd+d+2})$. Thus for any
$\la\not\in\calF_{N}(j_{0})$, applying Lemma \ref{secondocontrollo} we have
\begin{equation}\nonumber 
B_{N}^{0}(j_{0},\e,\la)\subset\bigcup_{q=1}^{N^{r+\gotd+d+4}}I_{q} \, , \quad
I_{q} \mbox{ intervals with } \meas(I_q) \le N^{-\tau_1} \, .
\end{equation}
But then, using also Lemma \ref{primapartizione}, we have that (recall \eqref{buoniautovaloriinf}
with $ \gote = r + \gotd + d + 4 $)
\begin{equation}\nonumber 
[1/2, 3/2]  \setminus \bar  \calG_{N}^{0}
\subset\bigcup_{|j_{0}|\le ( \pippa + 5 ) \pippa^{-1}N}\calF_{N}(j_{0}) \, .
\end{equation}
Hence,  using Lemma \ref{usofubini},   
\begin{equation}\nonumber 
{\rm meas} ({\mathcal I} \setminus \bar  \calG_{N}^{0} )\le\sum_{|j_{0}|\le ( \pippa + 5 ) \pippa^{-1}N}
\meas(\calF_{N}(j_{0}))\le O(N^{-1}) 
\end{equation}
which is the first bound in \eqref{meas.bad1}.
The second bound follows by \eqref{cattivinormalta} with  $ \tau_1 > d+ \gotd + 2 $.  
Finally, Lemmas \ref{misuralambdabarra} and \ref{lem.lambdatilde}   with  $ \g = N_0^{-1} $ implies
the third estimate in \eqref{meas.bad1}.
\EP

\subsection{Proof of Theorem \ref{teoremone} for NLS} 
\label{sec:nls} 

In order to apply Theorems \ref{principe}-\ref{mamma} to the Hamiltonian NLS,
we start by defining two extensions $\gotF(u,v),\gotH(u,v) $ of class $ C^q (\TTT^d \times {\mathtt M \times \CCC^2; \CCC}) $ (in the real sense)  of ${\mathtt f}(u) $
in such a way that $\gotF(u,\ol{u})=\gotH(u,\ol{u})={\mathtt f}(u)$ and
$\del_u \gotF(u,\ol{u})=\del_v \gotH(u,\ol{u})\in\RRR$,
$\del_{\ol u} \gotF(u,\ol{u})=\del_{\ol u}\gotH(u,\ol{u})=
\del_{\ol v} \gotF(u,\ol{u})=\del_{\ol v}\gotH(u,\ol{u})=0$ and
$\del_v \gotF(u,\ol{u})=\ol{\del_u\gotH(u,\ol{u})}$\footnote{where, for $ u = r + \ii s $, $ v = a + \ii b $ we set
$ \partial_u := (\partial_r - \ii \partial_s)/2 $, $ \partial_{\ol u} := (\partial_r + \ii \partial_s)/2 $, 
$ \partial_v := (\partial_a - \ii \partial_b)/2 $, $ \partial_{\ol v} := (\partial_a + \ii \partial_b)/2 $. A possible extension is the following: writing
 ${\mathtt f}(\f,x,r+\ii s)=\gotf_1(r,s)+\ii\gotf_2(r,s)$ we consider
%
 \begin{equation}\nonumber
\begin{aligned}
\gotF(\f,x,r+\ii s,a+\ii b)&:= (1+i)\gotf_1(\frac{r+a}{2},r-a+s)-i \gotf_1(\frac{r+a}{2}- \frac{s+b}{2},r-a+s) \\
&+ \gotf_2(2a-r-(s+b),\frac{r-a}{2}+\frac{s-b}{2})+ (-1+i)\gotf_2(a,\frac{s-b}{2})\\
\gotH(\f,x,r+\ii s,a+\ii b)&:=(1+i)f_1(\frac{r+a}{2},\frac{a-r}{2}+\frac{s-b}{2})-i f_1(\frac{r+a}{2}- \frac{s+b}{2},\frac{r-a}{2}+\frac{s-b}{2}) \\
&+f_2(\frac{r+a}{2}+\frac{s+b}{2},\frac{a-r}{2}+\frac{s-b}{2})- (1+i)f_2(a,\frac{s-b}{2}).
\end{aligned}
\end{equation}
}.
Then we   ``double'' the NLS equation, namely we look for a zero of the 
\emph{vector} NLS operator
\begin{equation}\label{vnls}
F(\e, \la, u^+,u^-):=\left\{
\begin{aligned}
\ii \la \overline \oo\cdot\del_{\f}u^{+}-\Delta u^{+}+mu^+-\e \gotF(\f,x,u^{+},u^{-})\\
-\ii \la \overline \oo\cdot\del_{\f}u^{-}-\Delta u^{-}+mu^--\e {\gotH}(\f,x,u^{+},u^{-})
\end{aligned}
\right.
\end{equation}
on the space $ H^s(\TTT^d\times G)\times H^s(\TTT^d\times G)$. 
Note that \eqref{vnls} reduces to \eqref{nls} on the invariant subspace
$$
{\cal U} := \{u=(u^{+},u^{-})\in H^{s}\times H^{s}\;:\; u^{-}=\ol{u^{+}}\}.
$$
Setting $\gotA=\{1,-1\}$, $\SSSS:=\SSSSG $,  
 we are in the functional setting  of Example \ref{zoppa2}, namely 
$ H^s(\TTT^d\times G)\times H^s(\TTT^d\times G)\equiv H^s(\gotK) $
where $ \gotK= \ZZZ^d\times \SSSSG\times \{1,-1\} $. Then equation \eqref{vnls} is of the form \eqref{eq.vera} with 
$$
D (\lambda) :=\begin{pmatrix}
\ii\la\ol{\oo}\cdot\del_{\f}-\Delta+m & 0 \\
0 & -\ii\la\ol{\oo}\cdot\del_{\f}-\Delta +m
\end{pmatrix},
\quad
f(u):= - \begin{pmatrix}
\gotF(\f,x,u^{+},u^{-})\\
{\gotH}(\f,x,u^{+},u^{-})
\end{pmatrix},
$$
and \eqref{Hdiag} holds with $\nu=2$.
Again the interpolation estimates (f1)--(f2) are verified if $s_0\ge (\gotd+d)/2$.

In the Fourier basis ${\rm e}^{\ii l\cdot\f}\ff_{j}(x)$ the operator
$D(\la) $ is represented by an infinite dimensional matrix
as in \eqref{operatore}, with $ D_{(l,j,\gota)}(\lambda)=  - \gota \lambda \ol{\oo}\cdot l-\mu_{j} +m $
and
$$
T^{(l,j)}_{(l',j')}:=\begin{pmatrix}
P^j_{j'}{(l-l')} & Q_{j'}^j{(l-l')} \\
\ol{Q}_{j'}^j{(l-l')} & P_{j'}^j{(l-l')}
\end{pmatrix},
$$
where $P_{j'}^{j}{(l)}, Q_{j'}^j{(l)}$ are the matrix representation in the Fourier basis of
the multiplication operators
$$
P(\f,x):= - \del_{u^+ } \gotF(\f,x,u^{+}(\f,x),u^{-}(\f,x)) ,
\quad
Q(\f,x):= - \del_{u^-} \gotH(\f,x,u^{+}(\f,x),u^{-}(\f,x))  \, .
$$
By the Hamiltonian assumption \eqref{HamNLS}, the constraints on
$\gotF,\gotH$ and $ ( u^-, u^+) \in {\cal U} $,   it results 
 $ P(\f,x) \in \RRR $ and   $(T^{(l',j')}_{(l,j)})^{\dagger}=T^{(l,j)}_{(l',j')}$. 
Hypothesis \ref{matricione} is then verified because
\eqref{hyp:diag} holds with $\gotD_{j,\gota}(y)= - \gota y -\mu_j+m $, 
$$
T^{(l',j',\gota)}_{(l,j,\gota)}=P_{j'}^j{(l-l')}\,, \ \gota = \pm 1 \, , 
\quad T^{(l',j',-1)}_{(l,j,1)}=Q_{j'}^j{(l-l')}\,,
$$
 Corollary \ref{nuzero} implies that $T\in \mathcal M^{s-\nu_0} $
and the estimates \eqref{hyp.normaT}, \eqref{lip} hold
by interpolation. 

We introduce the additional parameter $\theta$ and following  \eqref{innomin} we define the matrices 
$$
L(\e,\la,\theta,u):=D(\la,\theta)+\e T(\la,u), \quad 
D(\la,\theta):=D(\la)+\theta Y,\quad
Y:=\diag_{(l,j)\in\ZZZ^{d}\times\SSSSG}\begin{pmatrix}
- \uno_{d_j} & 0 \\
0 & \uno_{d_j}
\end{pmatrix} 
$$
so that  Hypothesis \ref{hyp.misura.diag} holds with  $ \gotn = 2 $.

Hypothesis \ref{hyp.separacatene} and 
the measure estimates needed in Theorem \ref{mamma} are obtained 
as in the case of NLW, following \cite{BB1} instead of \cite{BB2}.

\zerarcounters 
\section{An abstract Nash-Moser theorem}\label{see:nm1}

Let us consider a scale of Banach spaces $(X_s, \| \ \|_s)_{s \geq 0}$, such that
$$
\forall  s \leq s', \  \ X_{s'} \subseteq X_s \ \ {\rm and} \  \ \|u\|_s \leq \|u \|_{s'} \, , \ \forall u \in X_{s'} \, , 
$$
and define
$ X := \cap_{s\geq 0} X_s $.

We assume that there is a non-decreasing family $(E^{(N)})_{N \geq 0}$ of subspaces of $X$ such that
$\cup_{N \geq 0} E^{(N)}$ is dense in $X_s$ for any $s\geq 0$, and that there are projectors
$$
\Pi^{(N)}: X_0 \to E^{(N)}
$$
satisfying:  for any $s\ge0$ and any $\nu\ge0$ there is a positive
constant $ C :=C(s, \nu)$ such that

\begin{itemize}

\item[(P1)] $\|\Pi^{(N)}u\|_{s+\nu}\le CN^{\nu}\|u\|_{s}$
for all $u\in X_{s}$

\item[(P2)] $\|(\uno - \Pi^{(N)})u\|_{s}\le C
N^{-\nu}\|u\|_{s+\nu}$ for all $u\in X_{s+\nu}$.
\end{itemize}
In every Banach scale with smoothing operators satisfying (P1)-(P2) as above, the following interpolation inequality
holds (see Lemma 1.1 in \cite{BBP}): for all $ s_1 < s_2 $, $ t \in [0,1] $, 
\begin{equation}\label{interpolation}
\| u \|_{t s_1 + (1-t)s_2} \leq K(s_1, s_2) \| u \|_{s_1}^t  \| u \|_{s_2}^{1-t}  \, . 
\end{equation}

\begin{rmk}\label{pantheon}
The sequence spaces $H^s(\gotK)$ defined in \eqref{spazi}
admit spaces $ E^{(N)} := \{ u=\{u_k\}_{k\in \gotK}:\; u_k=0 \mbox{ for } |k|>N \} $ 
whose  corresponding projectors $\Pi^{(N)}$ satisfy $(P1)-(P2)$. 
\end{rmk}
Let us consider a  parameter family of $ C^2 $ maps
$ F:[0,\e_{0})\times\mathcal I\times X_{s_0 + \nu} \to X_{s_0} $
for some $s_0 \ge 0 $, $ \nu > 0 $, $ \e_{0}>0$ and $\mathcal I$
an interval in $\RRR$.  We assume
\begin{itemize}
\item[(F0)] $ F(0,\la,0)=0$ for any $\la\in\mathcal I$,
\end{itemize}
and  the following tame properties: given $ S' > s_0 $, $\forall s \in[s_0, S')$, for all 
$\|u\|_{s_0} \le1$, $(\e,\la)\in[0,\e_{0})\times\mathcal I $, 
\begin{itemize}
\item[(F1)] $\|\partial_{\la}F(\e,\la,u)\|_{s}\le
C(s)(1+\|u\|_{s+\nu})$,

\item[(F2)] $\|D_{u}F(\e,\la,0)[h]\|_{s}\le C(s)\|h\|_{s+\nu}$,

\item[(F3)] $\|D^{2}_{u}F(\e,\la,u)[h,v]\|_{s} \le C(s)\big(
\|u\|_{s+\nu}\|h\|_{s_0}\|v\|_{s_0}+\|h\|_{s+\nu}\|v\|_{s_0}
+\|h\|_{s_0}\|v\|_{s+\nu}\big)$,

\item[(F4)]  $\|\partial_{\la}D_{u}F(\e,\la,u)[h]\|_{s}\le
C(s)(\|h\|_{s+\nu}+\|u\|_{s+\nu}\|h\|_{s_0})$.
\end{itemize}
In application the following assumption is often verified 
\begin{equation}\label{Fsimme}
\exists \, \widehat X_s \subset X_s \ 
\mbox{closed subspaces  of } X_s, \, s \geq 0,  \, \mbox{such that} \ 
F:\widehat X_{s+\nu}\to \widehat X_s  \, . 
\end{equation} 
In order to prove the existence of a zero for $F$ we shall follow a Nash-Moser approach whose main assumption concerns the 
invertibility of the linearised
operators
$$
L^{(N)}(\e,\la,u):= \Pi^{(N)}D_{u}F(\e,\la,u)|_{E^{(N)}} 
$$
in a neighborhood of $u=0$. 

We introduce parameters $ \sigma $, $\tau > 0  $, $ s_1 > s_0 $, $\delta \in (0,1) $,
$S\in(s_0,S')$ satisfying
\begin{equation}\label{exponents}
\s\ge \max\{ 2(\tau+\de s_{1})+3\nu+2, 4(\tau + \delta s_1 +\nu)\}, \quad
2(2(\tau+\de s_{1})+\nu+3+\s)\le S-s_{1}\le 4(\s+1) \, . 
\end{equation}
Define the sets
\begin{equation}\label{buonini}
\begin{aligned}
J^{(N)}_{\tau,\de}&:=\Big\{
(\la,u)\in\mathcal I\times E^{(N)}\;:\; \|u\|_{s_1}\leq 1\,,\;
L^{(N)}(\e,\la,u)\mbox{ is invertible and, } \forall s\in[s_{1},S],\\ &\e\in[0,\e_{0}), \
\|L^{(N)}(\e,\la,u)^{-1}[h]\|_{s}\le { C(s) N^{\mu}}
(\|h\|_{s}+N^{\de({s}-s_{1})}\|u\|_{s}\|h\|_{s_{1}}), \mu := \tau + \de s_1 
\Big\} \, . 
\end{aligned}
\end{equation}
For $ K > 0 $ and  
 $ u\in C^{1}( \mathcal I,E^{(N)}) $ satisfying $ \|u\|_{s_{1}}\le1 $, $ \|\partial_{\la}u\|_{s_{1}}\le K $,  we set 
\begin{equation}\label{bonazzi}
{\mathtt G}^{(N)}_{\tau,\de}(u):=\{\la\in\mathcal I
\;:\; (\la,u(\la))\in J^{(N)}_{\tau,\de}\}.
\end{equation}

Given $N_{0}\in\NNN$ set $N_{n}:=N_{0}^{2^{n}}$
and denote with $E_{n},\Pi_{n},J^{n}_{\tau,\de}$ the subspace $E^{(N_{n})}$, the
projector $\Pi^{(N_{n})}$ and the set $J^{(N_{n})}_{\tau,\de}$ respectively.
Given any set $A$ and a positive real number $\h$ we denote by
$\NN(A,\h)$ the open neighborhood of $A$ with width $\h$.

\begin{theorem}\label{thm:nm1} {\bf (Nash-Moser)}
Assume (F0)-(F4). Then, for all $ \tau > 0 $, $\de\in(0,1/4)$, $ \s, s_1 >  s_0 $, $ S < S'$, satisfying \eqref{exponents}, 
there are $ c $, $ \ol{N}_{0}$, $ K_{0} > 0 $, such that, for all $N_0\ge \ol{N}_0$ and $ \e_0 $ small enough
such that 
\begin{equation}\label{piccoep} 
\e_0 N_0^S \leq c \, , 
\end{equation}
and, for all $\e\in[0,\e_0)$ a sequence
$\{u_n=u_{n}(\e,\cdot)\}_{n\ge0}\subset C^{1}(\mathcal I, X_{s_{1}+\nu})$ such
that

\begin{itemize}

\item[(S1)$_{n}$] $u_{n}(\e,\la)\in E_{n}$, $u_{n}(0,\la)=0$, $\|u_{n}\|_{s_{1}}
\le1$ and $\|\partial_{\la}u_{n}\|_{s_{1}}
\le K_{0}N_0^{\s/2}$. 
\item 
[($\widehat{S1}$)$_{n}$]
If \eqref{Fsimme} holds then $ u_n (\e, \la)\in 
E_{n} \cap \widehat{X}_s $.
\item[(S2)$_{n}$] For all $1\le i \le n$ one has
$\|u_{i}-u_{i-1}\|_{s_{1}}\le N_{i}^{-\s-1}$ and
$\|\partial_{\la}(u_{i}-u_{i-1})\|_{s_{1}}\le N_{i}^{-1-\nu}$.

\item[(S3)$_{n}$] Set $u_{-1}:=0$ and define 
\begin{equation}\label{defAn}
A_{n}:=\bigcap_{i=0}^{n}
{\mathtt G}^{(N_{i})}_{\tau,\de}(u_{i-1}) \, . 
\end{equation} 
For $\la\in\NN(A_{n},N_n^{-\s/2})$
the function $u_{n}(\e,\la)$ solves the equation
$ \Pi_{n}F(\e,\la,u)=0 $.
\item[(S4)$_{n}$] Setting $B_{n}:=1+\|u_{n}\|_{S}$ and
$B'_{n}:=\|\partial_{\la}u_{n}\|_{S}$ one has
\begin{subequations}
\begin{align}
&B_{n}\le2N^{p}_{n+1}\,,\qquad p:=\m+\frac{\nu}{2}+1 \, , 
\label{approx.normalta.a} \\
&B'_{n}\le 2N^{q}_{n+1}\,,\qquad q:= 2\m+\nu+2+ ( \s \slash 2) \, . 
\label{approx.normalta.b}
\end{align}
\label{approx.normalta}
\end{subequations}
\end{itemize}
\vskip-10mm
As a consequence, for all $\e\in[0,\e_0)$, 
 the sequence $\{u_{n}(\e,\cdot)\}_{n\ge0}$ converges uniformly in
$ C^{1}(\mathcal I, X_{s_{1}+\nu})$ to $u_\e$ with
$ u_0(\la)\equiv 0$, at a superexponential rate
\begin{equation}\label{exponentialrate}
\| u_\e (\la)  - u_n (\la) \|_{s_1} \leq N_{n+1}^{-\sigma-1} \, , \quad \forall \la \in {\cal I} \, ,
\end{equation}
and for all $ \la\in A_{\io}:=\bigcap_{n\ge0}A_{n}$ one has $ F(\e,\la,u_\e (\la))=0 $. 

Finally, if \eqref{Fsimme} holds then $ u_\e (\la) \in {\widehat X}_{s_{1}+\nu} $. 
\end{theorem}

\subsection{Proof of Theorem \ref{thm:nm1}} 
\label{sec:preliminari} 

Taylor formula and (F0)--(F4) imply the following 
tame properties: for any $s\in[s_0,S]$ there
is $C=C(s)$ such that for any $u,h\in X_{s}$ with  $\|u\|_{s_0}\le2$ and
$\|h\|_{s_0} \le1$ one has

\begin{itemize}

\item[(F5)] $\|F(\e,\la,u)\|_{s}\le C(s)(\e+\|u\|_{s+\nu})$,

\item[(F6)] $\|D_{u}F(\e,\la,u)[h]\|_{s}\le C(s)(\|u\|_{s+\nu}\|h\|_{s_0}+
\|h\|_{s+\nu})$

\item[(F7)] $\|F(\e,\la,u+h)-F(\e,\la,u)-D_{u}F(\e,\la,u)[h]\|_{s}
\le C(s)(\|u\|_{s+\nu}\|h\|_{s_0}^{2}+\|h\|_{s+\nu}\|h\|_{s_0})$.

\end{itemize}

The following Lemma follows as  in \cite{BBP},  Lemma 2.2.

\begin{lemma}\label{lem.bbp2}
Let $(\la,u)\in J^{(N)}_{\tau,\de}$ with $\|u\|_{s_{1}}\le 1$.
For $\e$ small enough, there exists $c=c(S)>0$ such that,  if
$|\la'-\la|+\|h\|_{s_{1}}\le c N^{-(\m+\nu)}$, 
$h\in E^{(N)}$, then $L^{(N)}(\e,\la',u+h)$ is invertible and
\begin{subequations}
\begin{align}
& \Big\|L^{(N)}(\e,\la',u+h)^{-1}[v] \Big\|_{s_{1}}  \le
2{N^{\m}}\|v\|_{s_{1}}, \quad \forall v\in E^{(N)} \, , 
\label{inv.a} \\
& \Big\|L^{(N)}(\e,\la',u+h)^{-1}[v]\Big\|_{S} \le
2{N^{\m}}\|v\|_{S} +KN^{\m}
\big({N^{\m+\nu}}(\|u\|_{S}+\|h\|_{S})
+{N^{\de(S-s_{1})}}
\|u\|_{S}\big)\|v\|_{s_{1}}.
\label{inv.b}
\end{align}
\label{inv}
\end{subequations}
\end{lemma}

\subsection{Initialisation of the Nash-Moser scheme} 
\label{inizio.nm1} 

Set $A_{0}:={\mathtt G}^{(N_{0})}_{\tau,\de}(0)$. By \eqref{bonazzi} we have that  $\la\in A_{0}$ if and only if
$(\la,0)\in J^{(N_{0})}_{\tau,\de}$. Therefore, if $N_{0}$ is large enough
Lemma \ref{lem.bbp2}
ensures that for all $\la\in \NN(A_{0},2N_{0}^{-\s/2})$ the operator $L^{(N_{0})}(\e,\la,0)$
is invertible for $\e$ small enough and 
\begin{equation}\label{initialbound}
\|L^{(N_{0})}(\e,\la,0)^{-1}\|_{s_{1}}\le 2 N_{0}^{\mu} ,\qquad
\|L^{(N_{0})}(\e,\la,0)^{-1}\|_{S}\le 2 N_{0}^{\mu} .
\end{equation}
Let us denote
\begin{equation*}
\begin{aligned}
&L_{0}:= L^{(N_{0})}(\e,\la,0), \, , \quad r_{-1}:=\Pi_{0}F(\e,\la,0) \\
&R_{-1}(u):=\Pi_{0}\left(F(\e,\la,u)-F(\e,\la,0)-D_{u}F(\e,\la,0)[u]\right).
\end{aligned}
\end{equation*}
We look for a solution of the equation $
\Pi_{0}F(\e,\la,u)=0, $ as a fixed point of the map
$$
\calH_{0}:E_{0} \longrightarrow E_{0} \, , \quad 
u  \longmapsto \calH_{0}(u):=-L_{0}^{-1}(r_{-1}+R_{-1}(u)) \, . 
$$
Let us  show that $\calH_{0}$ is a contraction in the set
$ \BBB_{\!\rho_0}= \{
u\in E_{0} :  \|u\|_{s_{1}}\le \rho_{0}:=c_{0}N_{0}^{\m}\e \}, $
for all $\e\in[0,\e_0(N_{0})]$ and some $c_{0}=c_{0}(s_{1})$.
We  bound
\begin{equation}\label{rinitial}
\|r_{-1}\|_{s_{1}}\le \| F(\e,\la,0)\|_{s_1}\stackrel{(F5)}{\le} C(s_{1})\e, \quad
\|R_{-1}\|_{s_{1}}\stackrel{(F7)}{\le} C(s_{1})\|\Pi_{0}u\|_{s_{1}+\nu}^2\stackrel{(P1)}{\le}
C(s_{1})N_{0}^{\nu}\|u\|_{s_{1}}^2,
\end{equation}
so that for any $u\in\BBB_{\!\rho_0}$ one has
$$
\begin{aligned}
\|\calH_{0}(u)\|_{s_{1}} & 
\stackrel{\eqref{initialbound}, \eqref{rinitial}}{\le} 2N_{0}^{\m}C(s_{1})(\e+N_{0}^{\nu}\|u\|_{s_{1}}^{2}) 
\stackrel{def}{\le}2C(s_{1})N_{0}^{\mu}\e+2C(s_{1})N_{0}^{\mu+\nu}\rho_{0}^{2} \le \rho_{0}
\end{aligned}
$$
where we have set $c_{0}=4 C(s_{1})$ and using \eqref{piccoep}.
This means that $\calH_{0}$ maps $\BBB_{\!\rho_0}$ into
iteself. In the same way (using (F3)) we obtain $\|D\calH_{0}(u)[h]\|_{s_{1}}\le\|h\|_{s_{1}}/2$
so that $\calH_{0}$ is a contraction on $(\BBB_{\!\rho_0},\|\cdot\|_{s_{1}})$ and hence it admits
a unique fixed point $\tilde{u}_{0}(\e,\la)$ for all $\la\in\NN(A_{0},2 N_{0}^{-\s/2})$.

Now, for $\la\in\NN(A_{0},2 N_{0}^{-\s/2})$ one has,
by $(F0)$, that $\tilde{u}_{0}(0,\la)=0$. The Implicit Function Theorem ensures
that $\tilde{u}_{0}(\e,\cdot)\in C^{1}(\NN(A_{0},2 N_{0}^{-\s/2});\BBB_{\!\rho_0})$ and
$\del_{\la}\tilde{u}_{0}=
L^{(N_{0})}(\e,\la,\tilde{u}_{0})^{-1}[\Pi_{0}\del_{\la}F(\e,\la,\tilde{u}_{0})]$. Hence
$$
\begin{aligned}
\|\del_{\la}\tilde{u}_{0}\|_{s_{1}}&\stackrel{\eqref{inv.a}}{\le}
2N_{0}^{\m}\|\Pi_{0}\del_{\la}F(\e,\la,\tilde{u}_{0})\|_{s_{1}}
\stackrel{(F1), (P1)}{\le}  2N_{0}^{\m}C(s_{1})(1+N_{0}^{\m+\nu}c_{0}\e) \le C(s_1) N_{0}^{\m}
\end{aligned}
$$
for $ \e $ small.  

We now define $u_{0}:=\psi_{0}\tilde{u}_{0}:[0,\e_0]\times\mathcal I\to E_{0}$ where
$\psi_{0}$ is a $C^{1}(\mathcal I, \RRR )$ cut-off function such that $0\le\psi_{0}\le1$ and 
\begin{itemize}
\item $\psi_{0}(\la)=1$ for $\la\in\NN(A_{0},N_{0}^{-\s/2})$
and $\psi_{0}(\la)=0$ for $\la\notin\NN(A_{0},2 N_{0}^{-\s/2})$,
\item $|\del_{\la}\psi_{0}|\le CN_{0}^{\s/2}$.
\end{itemize}
Of course $u_{0}$ satisfies (S3)$_{0}$. Moreover $u_{0}(0,\la)=0$ and satisfies also (S1)$_{0}$,
since
\begin{equation}\label{questa}
\|u_{0}\|_{s_{1}}\le 1 / 2,\qquad
\|\del_{\la}u_{0}\|_{s_1}\le CN_{0}^{\s/2}+ C(s_1) N_{0}^{\m}
\le C(s_1) N_{0}^{\s/2}
\end{equation}
because $ \s > 2 \m $. Finally, the bounds (S4)$_{0}$ follow in the same way.

\subsection{Iterative step} 
\label{sec.iterazione} 

Suppose inductively that we 
have defined $u_{n}\in C^{1}(\mathcal I,E_{n})$ such that properties (S1)$_{n}$ --
(S4)$_{n}$ hold. We define $u_{n+1}$ as follows. For $h\in E_{n+1}$ let us write
$$
\Pi_{n+1}F(\e,\la,u_{n}(\e,\la)+h)=r_{n}+ L_{n+1}[h]+R_{n}(h),
$$
where
$$
\begin{aligned}
r_{n}:=&\Pi_{n+1}F(\e,\la,u_{n}),\qquad
L_{n+1}:=L_{n+1}(\e,\la):=L^{(N_{n+1})}(\e,\la,u_{n}(\e,\la)), \\
&R_{n}(h):=\Pi_{n+1}\big(F(\e,\la,u_{n}(\e,\la)+h)-F(\e,\la,u_{n})
-D_{u}F(\e,\la,u_{n})[h]\big) \,. 
\end{aligned}
$$
Note that  (F7) and $ (S1)_n $ imply 
\begin{equation}\label{stimadelresto}
\|R_{n}(h)\|_{s}\le C(s)(\|u_{n}\|_{s+\nu}\|h\|_{s_{1}}^{2}+\|h\|_{s+\nu}
\|h\|_{s_{1}}),
\end{equation}
and for $ \la \in \NN(A_{n}, N^{-\s/2}_{n})$, we  use (S3)$_{n}$
to obtain
\begin{equation}\label{riscrivoilresto}
r_{n}=\Pi_{n+1}F(\e,\la,u_{n})-\Pi_{n}F(\e,\la,u_{n})=
\Pi_{n+1} ( \uno-\Pi_{n})F(\e,\la,u_{n}) \, .
\end{equation}
If $A_{n+1}=\emptyset$ we define $u_{i} := u_{n}$ for all $ i >n$, otherwise we proceed as follows.
By definition \eqref{defAn},  for all $ \la \in A_{n+1} $, the operator $L_{n+1}(\e,\la) $ is invertible. 
We also note that for $N_{0}$ large enough, one has
\begin{equation}\label{insiemotti}
\NN(A_{n+1},2 N_{n+1}^{-\s/2})\subset \NN(A_{n}, N_{n}^{-\s/2}) \, . 
\end{equation}

\begin{lemma}\label{invertibilitadecente}
For $\e$ small enough,  
$ \forall \la\in\NN(A_{n+1},2 N_{n+1}^{-\s/2}) $, the operator
$L_{n+1}(\e,\la)$ is invertible and  
\begin{subequations}
\begin{align}
\|L_{n+1}^{-1}[v]\|_{s_{1}}&\le 2{N^{\m}_{n+1}}\|v\|_{s_{1}},
\label{inversaL.a}\\
\|L_{n+1}^{-1}[v]\|_{S}&\leq KN_{n+1}^{\mu}\left(\|v\|_{S}+
\left(
N_{n+1}^{\x}B_n+ N_{n+1}^{(\m+\nu)-\s/2}B_n'\right)\|v\|_{s_{1}}\right)
\label{inversaL.b}
\end{align}
\label{inversaL}
\end{subequations}
where $ \x := \max\{\mu+\nu,\delta(S-s_{1}) \} $.
\end{lemma}

\prova
We apply Lemma \ref{lem.bbp2}.  
For $\la\in\NN(A_{n+1},2 N_{n+1}^{-\s/2})$ there is
$\la'\in A_{n+1}$ such that, setting $h(\la):=u_{n}(\e,\la)-u_{n}(\e,\la')$
one has (use $(S1)_{n}$)
$$
|\la-\la'|+\|h\|_{s_{1}} {\le} 3 N_{n+1}^{-\s/2}(1+
KN_{0}^{\s/2})\le c N_{n+1}^{-(\m+\nu)}
$$
using   \eqref{exponents} 
and $N_{0}$ large.
\EP

For all $\la\in\NN(A_{n+1},2 N_{n+1}^{-\s/2}) $, let consider the map
\begin{equation}\label{contraction}
\calH_{n+1}:E_{n+1} \longrightarrow E_{n+1} \, , \quad 
h \longmapsto\calH_{n+1}(h):=-L_{n+1}^{-1}[r_{n}+R_{n}] \, . 
\end{equation}
\begin{lemma}\label{puntofisso}
For $N_{0}$ large enough $\calH_{n+1}$ has a unique
fixed point $ {\tilde h}_{n+1}=\calH({\tilde h}_{n+1})$ in 
$ \BBB_{\!\rho_{n+1}}:=\{h\in E_{n+1}\;:\; \|h\|_{s_{1}}\le \rho_{n+1}:= N_{n+1}^{-\s-1}\} $.
Moreover, the following estimate holds 
\begin{equation}\label{poverinoi}
\| {\tilde h}_{n+1} \|_{s_1} \leq 2 C(s) N_{n}^{ - (s-s_1) + 2 \tau+ 2 \delta s_1 +  \nu } (1+ \| u_n \|_s) \, , \quad \forall s \in (s_1, S') \, .   
\end{equation}
\end{lemma}

\prova
Let us prove that $\calH_{n+1}$ is a contraction
in $\BBB_{\!\rho_{n+1}} $. 
For $\la\in\NN(A_{n+1},2 N_{n+1}^{-\s/2})$ we  use
\eqref{inversaL.a} and the definition of $\calH_{n+1} $ in order to bound
$ \|\calH_{n+1}(h)\|_{s_{1}}\le 2 N_{n+1}^{\m} \big(
\|r_{n}\|_{s_{1}}+\|R_{n}(h)\|_{s_{1}}\big) $.
Now we have 
\begin{equation}\label{rnbassa}
\|r_{n}\|_{s_{1}} \stackrel{\eqref{riscrivoilresto}, (P2), (F5)}{\le} C(S) N_{n+1}^{-(S-s_{1})/2}(\e+\|u_{n}\|_{S+\nu})
\stackrel{(P1)}{\le} C(S) N_{n+1}^{-(S-s_{1}-\nu)/2}B_{n}
\end{equation}
where  $ B_n := 1+ \| u_n \|_S $, see $ (S4)_n $. 
On the other hand
\begin{equation}\label{Rnbassa}
\begin{aligned}
\|R_{n}(h)\|_{s_{1}}&\stackrel{\eqref{stimadelresto},(P1)}{\le}
C(s_1)N_{n+1}^{\nu}\|h\|_{s_{1}}^{2}\le C(s_1)N_{n+1}^{\nu}\rho_{n+1}^{2}.
\end{aligned}
\end{equation}
Therefore
\begin{equation}\label{piccorR}
 2 N_{n+1}^{\mu}(\|r_{n}\|_{s_{1}}+\|R_{n}(h)\|_{s_{1}})\le \rho_{n+1} 
\end{equation}
  is implied by 
$ N_{n+1}^{-(S-s_{1}-\nu)/2}B_{n}+N_{n+1}^{\nu}\rho_{n+1}^{2}
\le \rho_{n+1}N_{n+1}^{-\mu-1} $, 
which in turn follows by
\eqref{approx.normalta.a},  \eqref{exponents}.   
So 
$\calH_{n+1}(\BBB_{\!\rho_{n+1}})\subseteq\BBB_{\!\rho_{n+1}}$. 
The derivative
$ D_{h}\calH_{n+1}(h)[v]=-L_{n+1}^{-1}\Pi_{n+1} (
D_{u}F(\e,\la,u_{n}+h)[v]-D_{u}F(\e,\la,u_{n})[v]),$
satisfies 
$$ 
\begin{aligned}
\left\|D_{h}\calH_{n+1}(h)[v]\right\|_{s_{1}}
&\stackrel{\eqref{inversaL.a}}{\le}
2 N_{n+1}^{\m}\left\|D_{u}F(\e,\la,u_{n}+h)[v]-D_{u}F(\e,\la,u_{n})[v]
\right\|_{s_{1}}\\
&\stackrel{(F3)}{\le} 4N_{n+1}^{\m}\Big(
\|u_n\|_{s_{1}+\nu}\|h\|_{s_{1}}\|v\|_{s_{1}}+\|h\|_{s_{1}+\nu}\|v\|_{s_{1}}
+\|h\|_{s_{1}}\|v\|_{s_{1}+\nu}\Big)\\
&\stackrel{(P1)}{\le}12N_{n+1}^{\m+\nu}\rho_{n+1}\|v\|_{s_{1}}
\stackrel{{\rm def}}{\le}12N_{n+1}^{\m+\nu-\s-1}\|v\|_{s_{1}}
\le \|v\|_{s_{1}} \slash 2 
\end{aligned}
$$
using  \eqref{exponents} and $N_{0}$ is large.
Then the Contraction Lemma implies the existence of a unique fixed point 
$ {\tilde h}_{n+1} = {\cal H}_{n+1} ({\tilde h}_{n+1}) $.
Now, by 
 \eqref{inversaL.a},  and the first inequality in \eqref{Rnbassa}, we get
 $$
\|\tilde{h}_{n+1}\|_{s_1} \le 2 N_{n+1}^\mu( \|r_n\|_{s_1}+C(s_1) N_{n+1}^\nu \|\tilde h_{n+1}\|_{s_1}^2)
\leq 2 N_{n+1}^\mu \|r_n\|_{s_1}+ 2 N_{n+1}^{\mu + \nu}C(s_1) \rho_{n+1} \|\tilde h_{n+1}\|_{s_1} \, .
$$
Using \eqref{rnbassa} with $S\rightsquigarrow s$, $B_n\rightsquigarrow  (1 + \| u_n \|_s )$, 
and  $ 2 N_{n+1}^{\mu + \nu}C(s_1) \rho_{n+1} < 1 /2 $, the bound \eqref{poverinoi} follows.
\EP

For $\la \in \NN(A_{n+1},2 N_{n+1}^{-\s/2})$, let $\tilde{h}_{n+1}(\e, \la) $ 
be the unique solution of $h=\calH_{n+1}(h)$.
It results $ \tilde{h}_{n+1} (0, \la) = 0 $. 

\begin{rmk}\label{cacchiocacchio} 
If \eqref{Fsimme} holds then 
$L_{n+1}, \mathcal H_{n+1}: E_{n+1}\cap \widehat{X}_s \to E_{n+1}\cap \widehat{X}_s $ 
and so $\tilde h_{n+1} (\e, \la) \in E_{n+1}\cap \widehat X_s $. 
\end{rmk}
 
\begin{lemma}\label{lem:normalta.sol}
$ \|\tilde{h}_{n+1}\|_{S}\le N_{n+1}^{2p}$ where $ 2p = 2 (\tau + \delta s_1 ) + \nu + 2 $, see \eqref{approx.normalta.a}.
\end{lemma}

\prova
Using \eqref{inversaL.b}, \eqref{piccorR}, $ \xi \geq \mu + \nu $ (see Lemma \ref{invertibilitadecente}), we estimate 
\begin{equation*}
\begin{aligned}
\|{\tilde h}_{n+1}\|_{S} & {\leq} K
N_{n+1}^\mu\left( \| r_{n}\|_{S}+
\|R_{n}( {\tilde h}_{n+1})\|_{S}+N_{n+1}^{\x-\mu-1}\rho_{n+1}\big(B_{n}+
N_{n+1}^{-\s/2}B_{n}' \big)\right)\\
&\le K N_{n+1}^\mu \Big(
N_{n+1}^{\nu/2}B_{n}(1+N_{n+1}^{\nu/2}\rho_{n+1}^{2})+N_{n+1}^{\nu}\rho_{n+1}\|{\tilde h}_{n+1}\|_{S}
+ N_{n+1}^{\x-\mu-1}\rho_{n+1}\big(B_{n}+N_{n+1}^{-\s/2}B_{n}' \big) \Big)
\end{aligned}
\end{equation*}
where in the last bound we used that, , $ \forall\; \|h\|_{s_{1}}\le \rho_{n+1}$, 
\begin{equation} \label{633b}
\|r_{n}\|_{S} \stackrel{(P1), (F5)} \le C N_{n+1}^{\nu/2}B_{n}, \quad
\|R_{n}(h)\|_{S} \stackrel{(P1), (F7)} \le K (N_{n+1}^{\nu}B_{n}\rho_{n+1}^{2}+
N_{n+1}^{\nu}\rho_{n+1} \|h\|_{S}) \, . 
\end{equation}
Now, since $N_{n+1}^{\nu+\mu}\rho_{n+1}\le 1/2$
by \eqref{exponents}, we  shift $ \| {\tilde h}_{n+1}\|_{S}$ on the
l.h.s. and obtain
$$
\| {\tilde h}_{n+1}\|_{S}\le K\big(
(N_{n+1}^{(\nu/2)+\mu}+N_{n+1}^{\x-\s-2})B_{n}+N_{n+1}^{\x- (3\s/2)-2}B_{n}' \big) 
$$
and the lemma follows by  \eqref{approx.normalta}, \eqref{exponents}.
\EP

\begin{lemma}
The map $ {\tilde h}_{n+1}$ is in $C^{1}(\NN(A_{n+1},2 N_{n+1}^{-\s/2});\BBB_{\!\rho_{n+1}})$
and 
\begin{equation}\label{derivhbassa} 
\|\del_{\la} {\tilde h}_{n+1} \|_{s_{1}} \le N_{n+1}^{-\nu-1}/ 2 , \quad \|\del_{\la} {\tilde h}_{n+1} \|_{S} \le N_{n+1}^{2q} \, . 
\end{equation} 
\end{lemma}

\prova
For $\la\in \NN(A_{n+1},2 N_{n+1}^{-\s/2})$ we have $U_{n+1}(\la, {\tilde h}_{n+1}(\la))\equiv 0$
where we have set $ U_{n+1}(\la,h):=\Pi_{n+1}F(\e,\la,u_{n}(\la)+h) $.
Hence 
\begin{equation}
0=\frac{{\rm d}}{{\rm d}\la}\big( U_{n+1}(\la, {\tilde h}_{n+1}(\la))\big)
=(\del_{\la}U_{n+1})(\la, {\tilde h}_{n+1}(\la))+
(D_{h}U_{n+1})(\la, {\tilde h}_{n+1}(\la))\del_{\la}{\tilde h}_{n+1}(\la).
\end{equation}
On the other hand, since $\| {\tilde h}_{n+1} \|_{s_{1}}\le N_{n+1}^{-\s-1}\ll c
N_{n+1}^{-(\mu+\nu )}$ for $N_{0}$
large enough (recall that $\mu+\nu<\s$), so that we  apply Lemma
\ref{lem.bbp2} and obtain 
\begin{subequations}
\begin{equation}
\|\big(D_{h}U_{n+1}(\la, {\tilde h}_{n+1})\big)^{-1}[v] \|_{s_{1}}
\le 4N_{n+1}^{\mu}\|v\|_{s_{1}},
\label{normabassaU}
\end{equation}
and, using also Lemma \ref{lem:normalta.sol}, 
\begin{equation}
\Big\| \Big(D_{h}U_{n+1} 
(\la, {\tilde h}_{n+1}) 
\Big)^{-1}[v] \Big\|_{S}
\le 4N_{n+1}^{\mu}\|v\|_{S}+ KN_{n+1}^{\mu}
\Big(B_{n}\left(N_{n+1}^{\nu+\m}+ 
N_{n+1}^{\de(S-s_{1})}\right) + N_{n+1}^{\mu+\nu+2p}
\Big)\|v\|_{s_{1}}.
\label{normaltaU}
\end{equation}
\end{subequations}
Therefore by the Implicit Function Theorem we have
$ {\tilde h}_{n+1} \in C^{1}(\NN(A_{n+1},2 N_{n+1}^{-\s/2});\BBB_{\!\rho_{n+1}})$
and
$$
\del_{\la} {\tilde h}_{n+1} =-\big((D_{h}U_{n+1})(\la, {\tilde h}_{n+1} (\la) )\big)^{-1}
(\del_{\la}U_{n+1})(\la, {\tilde h}_{n+1}(\la)).
$$
Now, by \eqref{insiemotti} we  use (S3)$_{n}$ to deduce
$$
\begin{aligned}
\del_{\la}U_{n+1}(\la,h)&= 
\Pi_{n+1}((\del_{\la}F)(\e,\la,u_{n}+h)-(\del_{\la}F)(\e,\la,u_{n}))\\
&\qquad+\Pi_{n+1}\left((D_{u}F)(\e,\la,u_{n}+h)-(D_{u}F)(\e,\la,u_{n})
\right)[\del_{\la}u_{n}])\\
&\qquad+\Pi_{n+1}(\uno-\Pi_{n})\left((\del_{\la}F)(\e,\la,u_{n})+(D_{u}F)(\e,\la,u_{n})
\right)[\del_{\la}u_{n}].
\end{aligned}
$$
Now, by (F1)-(F4), (F6), (P1)-(P2), (S1)$_{n}$, \eqref{normabassaU}, we  get
$$
\|\del_{\la} {\tilde h}_{n+1}\|_{s_{1}}\le
C N_{n+1}^{\mu}(N_{n+1}^{\nu-\s-1}+N_{n+1}^{- \frac{S-s_{1}}{2} + (\nu/2)}
(N_0^{\s/2}  B_{n}+B'_{n}))  \leq N_{n+1}^{-\nu-1}/2
$$
by \eqref{approx.normalta} and \eqref{exponents}. 
Now to get the estimate for the $S$-norm we use \eqref{normaltaU}
and obtain
$$ 
\begin{aligned}
\Big\| \Big(D_{h}U_{n+1} 
&(\e,\la,\tilde h) 
\Big)^{-1}[\del_{\la}U_{n+1}(\e,\la, {\tilde h}_{n+1} )] \Big\|_{S} {\le} 4N_{n+1}^{\mu}\|\del_{\la}U_{n+1}
(\e,\la, {\tilde h}_{n+1} )\|_{S}+
KN_{n+1}^{\mu}\\
&\qquad\qquad\qquad\times
\Big(B_{n} \Big(N_{n+1}^{\nu+\m}+ N_{n+1}^{\de(S-s_{1})}\Big) + 
N_{n+1}^{\mu+\nu+2p}\Big)\|\del_{\la}U_{n+1}(\e,\la,{\tilde h}_{n+1})\|_{s_{1}}\\
& {\le}
4 N_{n+1}^{\mu}( N_0^{\s/2} \|u_{n}\|_{S+\nu}+\| {\tilde h}_{n+1} \|_{S+\nu}+
\|\del_{\la}u_{n}\|_{S+\nu})+ KN_{n+1}^{\mu}\\
&\qquad\qquad\qquad\times
\Big(B_{n}\left(N_{n+1}^{\nu+\m}+ 
N_{n+1}^{\de(S-s_{1})}\right) + 
N_{n+1}^{\mu+\nu+2p}\Big)N_{n+1}^{-\m-\nu-1} \leq N_{n+1}^{2q}
\end{aligned}
$$
by \eqref{approx.normalta} and \eqref{exponents}. 
\EP

Let us define
\begin{equation}\label{extension}
h_{n+1}(\e,\la):=
\psi_{n+1}(\la) {\tilde h}_{n+1}(\e,\la) 
\end{equation}
where $\psi_{n+1}$ is a $C^{1}$ cut-off function such that $0\le\psi_{n+1}\le1 $ and
\begin{itemize}
\item $\psi_{n+1}(\la)=1$ for $\la\in\NN(A_{n+1},N_{n+1}^{-\s/2})$
and $\psi_{n+1}(\la)=0$ for $\la\notin\NN(A_{n+1},2 N_{n+1}^{-\s/2})$,

\item $|\del_{\la}\psi_{n+1}|\le N_{n+1}^{\s/2}$.

\end{itemize}

Then, by Lemma \ref{puntofisso}, \eqref{derivhbassa}, we get

\begin{lemma}\label{lem:extension}
One has $h_{n+1} \in C^{1}(\mathcal I;\BBB_{\!\rho_{n+1}})$ and 
\begin{equation}\label{stimafin+1}
h_{n+1}(0,\la)=0,\qquad
\|h_{n+1}\|_{s_{1}}\le N_{n+1}^{-\s-1},\qquad
\|\del_{\la}h_{n+1}\|_{s_{1}}\le N_{n+1}^{-\nu-1}.
\end{equation}
\end{lemma}

We now conclude the proof of Theorem \ref{thm:nm1}. Let 
\begin{equation}\label{un+1}
u_{n+1} := u_n+ h_{n+1} \, . 
\end{equation}
We want to show that (S1)$_{n+1}$-(S4)$_{n+1}$ are satisfied.
Property (S1)$_{n+1}$ follows by \eqref{questa} and \eqref{stimafin+1}.
Moreover Remark \ref{cacchiocacchio} implies that $ u_{n+1 } \in E_{n+1} \cap \widehat X_s $, i.e. 
($\widehat{S1}$)$_{n+1} $ holds.
Property (S2)$_{n+1}$ is  \eqref{stimafin+1}. 
Property (S3)$_{n+1}$ follows by the definition \eqref{extension} and since
$ {\widetilde h}_{n+1}(\e,\la)  $ solves $ \Pi_{n+1} F(\e,\la, u_{n} (\e,\la) + h ) = 0 $, for all
$ \la \in A_{n+1} $.
Finally
$$
B_{n+1}\le B_{n}+\|h_{n+1}\|_{S}\stackrel{(S4)_n, {\rm Lem.}\ref{lem:normalta.sol}}{\le}
2N_{n+1}^{p}+N_{n+1}^{2p}\le 2N_{n+2}^{p} 
$$
and $ (S4)_n $, 
Lemma \ref{lem:normalta.sol}, \eqref{derivhbassa}, imply 
$$
\begin{aligned}
B_{n+1}' &\le B_{n}' +\|\del_{\la}h_{n+1}\|_{S}  \le 2N_{n+1}^{q}+ N_{n+1}^{(\s/2)+2p}+N_{n+1}^{2q}
\le 2N_{n+2}^{q}
\end{aligned}
$$
because $ q = 2 p + (\s / 2) $.
Hence also (S4)$_{n+1}$ follows.
\EP

Note that so far the set $A_{\io} := \cap_{n\geq 0} A_n $ where $ u_\varepsilon $ is a solution of
$ F(\e, \la, u_\e (\la) ) = 0 $ (Theorem \ref{thm:nm1})
 could have zero measure or even be the
empty set. The goal of the next section is to show that, under 
further  assumptions on $ F $ (i.e. those of Theorem \ref{principe}), the set $ {\cal C}_\e $ 
 in \eqref{Cantorfinale}  (which is defined in a non inductive way)
is contained in $A_{\io}$. 

\zerarcounters 
\section{Proof of Theorem \ref{principe}} \label{section 5}
\label{sec:nm2} 

We now specialise   the abstract Nash-Moser Theorem \ref{thm:nm1} 
to the scale of  sequence spaces $ X^s = H^s(\gotK) $  and to operators $ F $ of the form \eqref{eq.vera}
satisfying the assumptions of Theorem \ref{principe}. 
In particular \eqref{Hdiag}-\eqref{Ddiag}, (f1)--(f2), imply the assumptions (F0)-(F4) of  Theorem \ref{thm:nm1}.

In addition to the parameters $ 	\tau > 0 $, $ \delta \in (0, 1/4 ) $, $ \sigma $,  $ s_1 > s_0 $, $  S<S'  $ satisfying
\eqref{exponents} needed in Theorem \ref{thm:nm1}, we now introduce other parameters 
 $\tau_{1} $,  $ \chi_0 $, $ \tau_0 $, $ C_1 $ and add the following constraints
\begin{equation}\label{esponenti}
s_1>s_0+\nu_0\,,\quad S< S'-\nu_0\,,\quad  
\tau>\tau_0 \, ,\quad \tau_{1}>2\chi_0d \, , 
\quad
\tau>2\tau_{1}+d+r+1,\quad C_{1}\ge 2 \, , 
\end{equation}
then, setting $\ka:=\tau+d+r+s_{0}$, $s_2:=s_1-\nu_0$
\begin{subequations}
\begin{align}
& \chi_0 (\tau-2\tau_{1}-d-r) > 3(\ka+(s_{0}+d+r)C_{1}),
\qquad \chi_0\de>C_{1},
\label{esponenti1.a} \\
&  s_{2} >3\ka+2\chi_0(\tau_{1}+d+r)+C_{1}s_{0}, 
\quad \ 2\de s_1>\nu_0 \, . 
\label{esponenti1.b}
\end{align}
\label{esponenti1}
\end{subequations}
Note that no restrictions from above on $S'$ are required, i.e. it  could be $S'=+\io $.

\begin{rmk}
In the applications, the constants $ \tau_0 $, $ \tau_1 $ have to be taken large enough, 
in order to verify condition \eqref{meas.bad1}. 
Nevertheless, all the constraints \eqref{exponents}, \eqref{esponenti}, \eqref{esponenti1} may be verified.
\end{rmk}

Given $\Omega,\Omega'\subset{\mathfrak K} $, we define
$$
{\rm diam}(\Omega):=\sup_{k,k'\in \Omega}\dist(k,k'),
\qquad
{\rm dist}(\Omega,\Omega'):=\inf_{\substack{k\in\Omega, k'\in\Omega'}}
\dist(k,k') \, , 
$$
where $\dist(\cdot,\cdot)$ is defined in \eqref{modulo}.

\begin{defi}\label{Ngood} {\bf ($N$-good/$ N$-bad matrices).}
Let $F\subset{\mathfrak K}  $ be such that ${\rm diam}(F)\le4N$
for some $N\in \NNN$.
We say that a matrix $A\in\MM_{F}^{F}$ is $N$\emph{-good}
if $A$ is invertible and for all $s\in[s_{0},s_{2}]$ one has
$$
\bvert A^{-1}\bvert_{s}\le N^{\tau+\de s}.
$$
Otherwise we say that $A$ is $N$\emph{-bad}.
\end{defi}

\begin{defi}\label{AN-reg} {\bf ($(A,N)$-regular, good, bad sites).}
For any finite  $E\subset{\mathfrak K} $, let $ A = D + \e T \in \MM^E_E $ with 
 $ D := {\rm diag}( D_k \uno_{d_j}) $, $ D_k \in \CCC $. 
An index $k\in E$ is
\begin{itemize}

\item $(A,N)$-\emph{regular} if there
exists $F\subseteq E$ such that ${\rm diam}(F)\le4N$,
${\rm dist}(\{k\},E\setminus F)\ge N$ and the matrix $A^{F}_{F}$ is $N$-good.

\item $(A,N)$-\emph{good} if either it is regular for $D$ (Definition \ref{regular}) or it is
$(A,N)$-regular. Otherwise $k$ is $(A,N)$-\emph{bad}.
\end{itemize}
\end{defi}
The above definition could be extended to infinite $E$. 

\begin{defi}\label{param.N-buoni}
{\bf ($N$-good/$ N$-bad parameters).}
For $ L $ as in \eqref{innomin}, we denote
\begin{equation}\label{bad}
B_{N}(j_0,\e,\la):= \Big\{\theta\in\RRR\,:\, L_{N,j_0}(\e,\la,\theta,u) \mbox{ is }N\mbox{--bad}\, \Big\}.
\end{equation}
A parameter $\la\in\mathcal I$ is $N$--good for $L $ if for any $j_0 \in \SSSS$ one has
\begin{equation}\label{partizione}
B_{N}(j_0,\e,\la)\subseteq \bigcup_{q=1}^{N^{\gote}} I_{q} \, , \quad 
I_{q} \mbox{ intervals with } \meas(I_q) \le N^{-\tau_{1}} 
\end{equation}
where  $ \gote $ is the parameter introduced in Theorem \ref{principe}.
Otherwise we say that $\la$ is
$N$--bad. We denote the set of $N$--good parameters as
\begin{equation}\label{param.good}
\calG_{N}=\calG_{N}(u):= \Big\{\la\in\mathcal I\,:\, \la\mbox{ is }N\mbox{--good for }L \Big\}.
\end{equation}
\end{defi}

Note that the above definition deals only with  
finite dimensional truncations of $ L $.

The following assumption is needed for the multiscale Proposition \ref{multiscala}.

\begin{hyp}\label{hyp.separazione} {\bf (Separation of bad sites)} 
There exist $ C_{1} = C_{1}(\mathfrak{K})>2$, 
$\hat{N}=\hat{N}(\mathfrak{K},\tau_{0})
\in \NNN $ and $\hat {\mathcal I} \subseteq 
\overline{\mathcal I}$ (see \eqref{lambdabarra})
such that, for all $ N \ge \hat{N} $, and $ \| u \|_{s_1} < 1 $ (with $s_1$ satisfying \eqref{esponenti1.b}), if
$$
\la\in\calG_{N}(u)\cap \hat {\mathcal I},
$$
then for any $\theta\in\RRR$, for all $\chi\in[\chi_0,2\chi_0]$ and all $j_0\in\Lambda_+$
the $(L,N)$-bad sites $k=(l,j, \gota)\in\mathfrak{K}$
of $L=L_{N^{\chi},j_0}(\e,\la,\theta,u)$
admit a partition $\cup_{\al}\Omega_{\al}$ in disjoint
clusters satisfying
\begin{equation}\label{separazione}
{\rm diam}(\Omega_{\al})\le N^{C_{1}},
\qquad
{\rm dist}(\Omega_{\al},\Omega_{\be})\ge N^{2},
\;\mbox{ for all }\al\ne\be.
\end{equation}
\end{hyp}

For $ N > 0 $, we  denote 
\begin{equation}\label{buoniautovalori}
\begin{aligned}
\calG_N^{0}(u):=\Big\{& \la\in\mathcal I\;:\;
\forall\; j_{0}\in \SSSS  \mbox{ there is a covering } \\
B_{N}^{0}&(j_{0},\e,\la)\subset \bigcup_{q=1}^{N^{\gote}}I_{q}, \quad I_{q}=I_{q}(j_{0})\mbox{ 
intervals with }
\meas(I_{q})\le  N^{-\tau_1}
\Big\} 
\end{aligned}
\end{equation}
where
\begin{equation}\label{tetacattiviautovalori}
B_{N}^{0}(j_{0},\e,\la) := B_{N}^{0}(j_{0},\e,\la, u) := \Big\{\theta\in\RRR\;:\;
\|L_{N,j_{0}}^{-1}(\e,\la,\theta, u)\|_0>N^{\tau_1}\Big\} \, . 
\end{equation}
We also set
\begin{equation}\label{buoninormal2}
\gotG_{N}(u):=\Big\{\la\in\mathcal I\;:\; \|L^{-1}_{N}(\e,\la,u)\|_{0}\le
N^{\tau_1} \Big\} \, .
\end{equation}
Under the smallness condition \eqref{piccoep}, Theorem  \ref{thm:nm1} applies, thus 
defining the sequence $ u_ n $ and the sets $ A_n $.
We  now introduce the sets
\begin{equation}\label{Cn}
\CCCC_{0}:=\hat{\mathcal I},\qquad
\CCCC_{n}:=\bigcap_{i=1}^n \calG^{0}_{N_i}(u_{i-1}) \bigcap_{i=1}^n \gotG_{N_i}(u_{i-1})
 \cap \hat {\mathcal I}
\end{equation}
where $\hat {\mathcal I}$ is defined in Hypothesis  \ref{hyp.separazione},  
$\gotG_{N}(u)$  in \eqref{buoninormal2}, and  $\calG^{0}_{N}(u)$ in \eqref{buoniautovalori}. 

\begin{theo}\label{thm:nm2}
Assume that  $ F $ in \eqref{eq.vera} satisfies \eqref{Hdiag}-\eqref{Ddiag}, (f1)-(f2)  and
Hypotheses \ref{matricione}, \ref{hyp.misura.diag} and \ref{hyp.separazione}. 
Assume that the parameters satisfy  \eqref{exponents}, \eqref{esponenti}, \eqref{esponenti1}. 
Then
there exists $ \ol{N}_0 \in \NNN $, such that, for all $ N_0 \geq  \ol{N}_0 $ and $   \e \in [0, \e_0) $
with $ \e_0 $ satisfying \eqref{piccoep}, 
the following inclusions hold:
\begin{align*}
&(S5)_{0}\qquad\quad \;
\|u\|_{s_{1}}\le 1 \quad \Rightarrow\quad
\calG_{N_{0}}(u) ={\mathcal I} \\
&(S6)_0 \qquad\qquad \CCCC_{0}\subseteq A_{0},
\end{align*}
and for all $n\ge1$ (recall the definitions of $ A_n $ in \eqref{defAn}) 
\begin{align*}
&(S5)_{n}\qquad\quad \;
\|u-u_{n-1}\|_{s_{1}}\le N_{n}^{-\s} \quad \Rightarrow\quad
\bigcap_{i=1}^{n}\calG^{0}_{N_i}(u_{i-1})\cap \hat {\mathcal I}\subseteq
\calG_{N_{n}}(u)\cap\hat {\mathcal I}, \\
&(S6)_{n}\qquad\qquad \CCCC_{n}\subseteq A_{n} \, .
\end{align*}
Hence  $ \CCCC_{\io}:=\bigcap_{n\ge0} \CCCC_{n} \subseteq A_\infty :=\bigcap_{n\ge0} A_n $. 
\end{theo}

\subsection{Initialisation}
\label{sec:inizio.nm2} 

Property  $(S5)_0$ follows from the following Lemma.

\begin{lemma}\label{iniz.1}
For all $\|u\|_{s_{1}}\le1$, $N\le N_0 $,   the set 
$\calG_{N}(u)=\mathcal I$.
\end{lemma}

\prova
We claim that, for any $\la\in\mathcal I$ and any $j_{0}\in\SSSS$, if 
\begin{equation}\label{autov.grandi}
|D_{k}(\la,\theta)| > N^{-\tau_{1}},\quad\forall k=(l,j,\gota)\in \gotK \; \mbox{with}\; |(l,j-j_{0})|\le N \, , 
\end{equation}
then $ L_{N,j_{0}}(\e,\la,\theta)\;\mbox{ is }N\mbox{--good} $.
This implies that 
$$
B_{N}(j_{0},\e,\la)\subset \bigcup_{ |(l,j-j_{0})|\le N}\left\{ \theta\in\RRR\;:\; |D_{k}(\la,\theta)| \leq N^{-\tau_{1}}
\right\},
$$
which  in turn, by  Hypothesis \ref{hyp.misura.diag}, implies the thesis, see \eqref{partizione}, \eqref{param.good}, for some
$ \gote \geq d + r + 1 $.
The above claim  follows by a perturbative argument.
Indeed for $\|u\|_{s_{1}} \le1 $, $ s_1 = s_2 + \nu_0 $, we
use \eqref{hyp.normaT} to obtain
$$
\e\bvert (D_{N,j_{0}}^{-1}(\la,\theta))\bvert_{s_{2}}\bvert T_{N,j_{0}}(u)\bvert_{s_{2}} \leq \e C(s_1)
\bvert D_{N,j_{0}}^{-1}(\la,\theta) \bvert_{s_{2}}(1+\|u\|_{s_2 + \nu_0})
\stackrel{\eqref{autov.grandi}}{\le}
\e N^{\tau_{1}}C(s_{1})\stackrel{\eqref{piccoep}}{\le} \frac{1}{2} \, . 
$$
Then we  invert $L_{N,j_{0}}$ by Neumann series and Lemma \ref{lem.inversasinistra}
implies
$$
\bvert L_{N,j_{0}}^{-1} (\e,\la,\theta)\bvert_{s}\le 2\bvert D_{N,j_{0}}^{-1}(\la,\theta) \bvert_{s}
\le 2 N^{\tau_{1}} {\le} N^{\tau+\de s}, \quad \forall s \in[s_{0},s_{2}]  \, , 
$$
by \eqref{esponenti}, which proves the claim.
\EP

\begin{lemma}\label{s60}
Property $(S6)_0$ holds.
\end{lemma}

\prova
Since $\hat{\mathcal I}\subset\ol{\mathcal I}$ it is sufficient to prove that $ \ol{\mathcal I} \subset A_0 $.
By the definition of $ A_0 $ in \eqref{defAn}, \eqref{bonazzi}, \eqref{buonini}, 
we have to prove that
\begin{equation}\label{cio da dim}
\la\in\ol{\mathcal I} \quad \Longrightarrow \quad 
\|L_{N_0}^{-1}(\e,\la,0)[h]\|_{s} \le C(s)N_0^{\tau + \delta s_1 } \|h\|_s \, , \ \forall s \in [s_1, S] \, .
\end{equation}
Indeed, if $\la\in\ol{\mathcal I}$ then $|D_{k}(\la)| \geq N_0^{-\tau_0}$, for all $|k|<N_0$, and so
$ \bvert D(\la)^{-1} \bvert_{s} \le  N_0^{\tau_0} $, $ \forall s $. 
Since $ \e \bvert D(\la)^{-1}T(0)\bvert_{s_1}\le \e N_0^{\tau_0}\bvert T(0)\bvert_{s_1}<1/2$ for
$\e$ small enough,  Lemma \ref{lem.inversasinistra} 
implies
$$
\bvert L_{N_0}^{-1}(\e,\la,0)\bvert_{s_1}\le 2N_{0}^{\tau_0},\qquad
\bvert L^{-1}_{N_0}(\e,\la,0)\bvert_{s}\le C(s)N_{0}^{\tau_0}(1+\e N_{0}^{\tau_0}\bvert T(0)\bvert_s) \, ,\forall s > s_1 \, ,
$$
 and, by \eqref{esponenti}, \eqref{hyp.normaT} and Lemma \ref{lem.questo}, the estimate  \eqref{cio da dim} follows.
\EP

\subsection{Inductive step} \label{sec:ter.nm2} 

By the Nash-Moser Theorem \ref{thm:nm1} we know that (S1)$_{n}$--(S4)$_{n}$ hold for all $n\ge0$.
Assume inductively that (S5)$_i $ and (S6)$_i $ hold for all $i \le n$. In order to prove
(S5)$_{n+1}$, we need the following {\em multiscale Proposition} \ref{multiscala} which
allows to deduce estimates on the $\bvert\cdot \bvert_s$--norm of the
inverse of $L$ from informations on the $L^2$-norm of the inverse $ L^{-1} $, 
the off-diagonal decay of $ L $, and separation properties of the bad sites.

\begin{prop}\label{multiscala} {\bf (Multiscale)} Assume \eqref{esponenti}, \eqref{esponenti1}. For any $ \ol{s} > s_2 $, 
$\Upsilon>0$ there exists $\e_0=\e_0(\Upsilon,s_{2})>0$
and $N_{0}=N_{0}(\Upsilon, \ol{s})\in\NNN$ such that, for all $N\ge N_{0}$, $|\e|< \e_0$, $\chi\in [\chi_0,2\chi_0]$, 
$ E\subset{\mathfrak K}  $ with ${\rm diam}(E)\le 4N^{\chi}$, 
if the matrix $A=D+\e T\in\MM^{E}_{E}$ satisfies

\begin{itemize}
\item[{\rm (H1)}] $\bvert T\bvert_{s_{2}}\le \Upsilon$,
\item[{\rm (H2)}] $\|A^{-1}\|_{0}\le N^{\chi\tau_{1}}$,
\item[{\rm (H3)}] there is a partition $\{\Omega_{\al}\}_{\al}$
of   the $(A,N)$-bad sites (Definition \ref{AN-reg}) such that
$$
{\rm diam}(\Omega_{\al})\le N^{C_{1}},
\qquad
{\rm dist}(\Omega_{\al},\Omega_{\be})\ge N^{2},\;\mbox{ for }\al\ne\be,
$$
\end{itemize}

then the matrix $A$ is $N^{\chi}$-good and
\begin{equation}\label{stimabella}
\bvert A^{-1}\bvert_{s}\le\frac{1}{4}N^{\chi\tau}\left(N^{\chi\de s}+\e
\bvert T\bvert_{s}\right) \, , \quad \forall s \in [ s_0,\ol{s}] \, .  
\end{equation}
\end{prop}
Note that the bound \eqref{stimabella} is much more than requiring that the matrix $A$ is 
$N^{\chi}$--good, since it holds also for $s > s_{2} $.

This Proposition is proved by ``resolvent type arguments'' and it coincides  essentially with \cite{BB1}-Proposition 4.1. The correspondences in the notations of this paper and \cite{BB1}
respectively are the following:
$(\tau,\tau_1,d+r,s_2, \ol{s})\rightsquigarrow (\tau',\tau,b,s_1, S)$, 
and, since we do not have a potential, we can fix  $\Theta=1$ in  Definition 4.2 of \cite{BB1}. 
Our conditions \eqref{esponenti}, \eqref{esponenti1} imply conditions (4.4) and (4.5) of \cite{BB1} 
for all $\chi\in [\chi_0,2\chi_0]$ and our (H1) implies the corresponding 
Hypothesis (H1) of \cite{BB1}  with $\Upsilon\rightsquigarrow 2\Upsilon$. The other hypotheses are the same.   
Although  the $s$--norm  in this paper is different, the proof of  \cite{BB1}-Proposition 4.1 relies only on abstract algebra and 
interpolation properties of the $s$--norm (which indeed hold also in this case -- see 
section \ref{sub.linearop}). Hence it can be repeated verbatim and we report it in the Appendix for completeness.

\smallskip

Now,  we distinguish two cases:

\begin{itemize}

\item[{\bf case 1:}] $2^{n+1} \le \chi_0$. Then there exists $\chi\in[\chi_{0},2\chi_{0}]$ (independent of $n$) such that
\begin{equation}\label{case1}
N_{n+1}=\ol{N}^{\chi},\qquad \ol{N}:=[N_{n+1}^{1/\chi_0}]\in(N_{0}^{1/\chi},N_{0}) \, . 
\end{equation}
 This case may occur only in the first steps.

\item[{\bf case 2:}] $2^{n+1}>\chi_{0}$. Then there exists a unique $p\in[0,n]$ such that
\begin{equation}\label{case2}
N_{n+1}=N_{p}^{\chi},\qquad \chi=2^{n+1-p}\in[\chi_{0},2\chi_{0}) \, . 
\end{equation}
\end{itemize}

Let us start from {\bf case 1} for $n+1=1$; the other (finitely many) steps are identical.

\begin{lemma}\label{buonobarra}
Property (S5)$_{1}$ holds.
\end{lemma}

\prova
We have to prove that $ \calG_{N_1}^0(u_0) \cap \hat {\cal I} \subseteq {\cal G}_{N_1} (u) \cap \hat {\cal I} $
where $\|u-u_{0}\|_{s_{1}}\le N_{1}^{-\s}$.  By Definition \ref{param.N-buoni} and \eqref{buoniautovalori}
it is sufficient to prove that, for all $j_{0}\in\SSSS$, 
$$
B_{N_{1}}(j_{0},\e,\la,u)\subseteq B_{N_{1}}^{0}(j_{0},\e,\la,u_{0}),
$$
where we stress the dependence on $u,u_{0}$  in \eqref{bad}, \eqref{tetacattiviautovalori}. 
By the definitions \eqref{tetacattiviautovalori}, \eqref{bad} this amounts to prove that 
\begin{equation}\label{vweek}
\|L_{N_1,j_0}^{-1}(\e,\la,\theta,u_0)\|_{0}\le N_1^{\tau_1} \quad \Longrightarrow \quad  L_{N_1j_0}(\e,\la,\theta,u)\mbox{ is }N_1-\mbox{good} \, . 
\end{equation}
We first claim that $ \|L_{N_1,j_0}^{-1}(\e,\la,\theta,u_0)\|_{0}\le N_1^{\tau_1} $ implies
\begin{equation}\label{stimapotenteuno}
\bvert L^{-1}_{N_{1},j_{0}}(\e,\la,\theta,u_{0})\bvert_{s}\le
\frac{1}{4} N_{1}^{\tau}\left(N_{1}^{\de s}+\e
\bvert T(u_{0})\bvert_{s}\right) \, , \quad \forall  s\in[s_0,S] \, .
\end{equation}
Indeed we may apply Proposition \ref{multiscala} to the matrix $A=L_{N_{1},j_{0}}(\e,\la,\theta,u_{0})$ with $\ol{s}=S$,
$ N = \ol{N} $, $ N_1 = \ol{N}^\chi $ and  $E=\{ |l|\le N_1,|j-j_0|\le N_1\}\times\gotA$. 
Hypothesis (H1)  with $\Upsilon = 3 C(s_1) $ 
 follows by  \eqref{hyp.normaT} and $ \| u_0 \|_{s_1 } \leq 1 $. Moreover (H2) is
$ \| L_{N_1,j_0}^{-1}(\e,\la,\theta,u_0)\|_{0}\le N_1^{\tau_1} $ . 
Finally (H3) is  implied by  Hypothesis \ref{hyp.separazione} 
provided we take $N_0^{1/\chi_0}> \hat{N}(\gotK, \tau_0)$ (recall \eqref{case1})
and noting that  $ \la \in \calG_{\overline N}(u_0) \cap \hat {\cal I} $
by Lemma \ref{iniz.1} (since $ \overline{N} \leq N_0 $ then $ \calG_{\overline N}(u_0) = \mathcal I  $). 
 Hence \eqref{stimabella} implies  \eqref{stimapotenteuno}. 
 
 We now prove \eqref{vweek} by a perturbative argument. 
Since  $\|u-u_{0}\|_{s_{1}}\le N_{1}^{-\s}$ (recall that $\|u_0\|_{s_1}\le 1$ so $\|u\|_{s_1}\le 2$)
then, for $\nu_1 := \max(\nu,\nu_0)$,  
$$
\bvert L_{N_{1},j_{0}}(\e,\la,\theta,u_{0})-L_{N_{1},j_{0}}(\e,\la,\theta,u)\bvert_{s_{2}}
\leq \bvert L_{N_{1},j_{0}}(\e,\la,\theta,u_{0})-L_{N_{1},j_{0}}(\e,\la,\theta,u)
\bvert_{s_{1}-\nu_1}
$$
$$ \stackrel{\eqref{Hdiag}, \eqref{lip}} \le C 
\|u-u_{0}\|_{s_{1}}\le C N_{1}^{-\s} < 1 / 2 \, . 
$$
By Neumann series (see Lemma \ref{lem.inversasinistra}) 
and \eqref{stimapotenteuno} one has
$ \bvert L^{-1}_{N_{1},j_{0}}(\e,\la,\theta,u)\bvert_{s}\le 
N_{1}^{\tau+\de s} $
for all $s\in [s_{0},s_{2}]$, namely  $ L_{N_{1},j_{0}}(\e,\la,\theta,u)$ is $ N_1 $-good.
\EP

\begin{lemma}\label{buonobarraS6}
Property (S6)$_{1}$ holds.
\end{lemma}

\prova
Let $ \la \in \CCCC_1 := \calG^0_{N_{1}}(u_0)\cap \gotG_{N_1}(u_0)  \cap \hat {\mathcal I} $, 
see \eqref{Cn}. By the definitions \eqref{defAn}, \eqref{bonazzi}, and (S6)$_0 $, 
in order to prove that $ \la \in A_1 $, 
it is sufficient to prove that $(\la,u_{0}(\la)) \in J^{(N_{1})}_{\tau,\de}$. Since 
$ \la \in 
\gotG_{N_1}(u_0)  $ the matrix $ \| L_{N_1}^{-1}(\e,\la,u_0) \|_0 \leq N_1^{\tau_1} $ (see \eqref{buoninormal2}) and so
\eqref{stimapotenteuno} holds with $ j_0 = 0 $, $ \theta = 0  $. 
Hence, using  \eqref{hyp.normaT} and  $u_0\in E_0 $,
we have
$$
\bvert L_{N_1}^{-1} (\e,\la,u_0)\bvert_s\le \frac{1}{4}N_{1}^{\tau}(N_1^{\de s}+\e C(s)
(1+N_1^{\nu_0/2}\|u_0\|_s)) \, , \quad \forall s \in [s_0, S] \, ,
$$
that, by Lemma \ref{lem.questo}, $ \de s_1>\nu_0/2 $, (P1), implies, $ (S1)_n $, $ \forall s \in [s_1, S]  $, 
$$
\|L_{N_1}^{-1}(\e,\la,u_0) [h]\|_{s}\le { C(s) N_1^{\tau + \de s_1}}
(\|h\|_{s}+N_1^{\de({s}-s_{1})}\|u_0 \|_{s}\|h\|_{s_{1}}),  \quad \forall h \in E_1 \, ,
$$
which is the inequality in \eqref{buonini} with $ N = N_1 $, $ u = u_0 $. 
Hence  $(\la,u_{0}(\la)) \in J^{(N_{1})}_{\tau,\de}$.  
\EP

Now we consider {\bf case 2}.

\begin{lemma}\label{lem:boh}
$
\bigcap_{i=1}^{n+1}\calG_{N_{i}}^{0}(u_{i-1})\cap \hat{\mathcal I} \subseteq \calG_{N_p}(u_n)\cap \hat{\mathcal I} $.
\end{lemma}

\prova
By $(S2)_{n} $ of Theorem \ref{thm:nm1} we get
$ \|u_{n}-u_{p-1}\|_{s_{1}}\le $ $ \sum_{i=p}^{n}\|u_{i}-u_{i-1}\|_{s_{1}} \le $ $ \sum_{i=p}^{n}N_{i}^{-\s-1}{\le} $ $ N_{p}^{-\s}\sum_{i=p}^{n}N_{i}^{-1}\le N_{p}^{-\s} $. Hence  $ (S5)_{p} $ ($ p \le n $) 
implies 
$$
\bigcap_{i=1}^{n+1}\calG_{N_{i}}^{0}(u_{i-1})\cap \hat{\mathcal I}\subseteq
\bigcap_{i=1}^{p}\calG_{N_{i}}^{0}(u_{i-1})\cap \hat{\mathcal I}\stackrel{(S5)_{p}}{\subseteq} \calG_{N_{p}}(u_{n})\cap \hat{\mathcal I} 
$$
proving the lemma.
\EP

\begin{lemma}\label{induz:s5}
Property (S5)$_{n+1}$ holds.
\end{lemma}

\prova
Fix $\la\in \bigcap_{i=1}^{n+1}\calG_{N_{i}}^{0}(u_{i-1})\cap \hat{\mathcal I}                                         $.
Reasoning as in the proof of Lemma \ref{buonobarra}, it is sufficient to prove that, 
for all $ j_0 \in\SSSS $, $\|u-u_{n}\|_{s_{1}}\le N_{n+1}^{-\s}$, 
one has
\begin{equation}\label{cisiamo}
\|L^{-1}_{N_{n+1},j_{0}}(\e,\la,\theta,u_{n})\|_{0}\le N_{n+1}^{\tau_{1}} \quad \Longrightarrow \quad 
L_{N_{n+1},j_{0}}(\e,\la,\theta,u)\mbox{ is }N_{n+1}\mbox{--good} \, .
\end{equation}
We apply the multiscale Proposition \ref{multiscala} 
to the matrix $A=L_{N_{n+1},j_{0}}(\e,\la,\theta,u_{n})$ with $ N^\chi =  N_{n+1} $ and $ N = N_p $, see \eqref{case2}.
Assumption (H1) holds and (H2) is $ \|L^{-1}_{N_{n+1},j_{0}}(\e,\la,\theta,u_{n})\|_{0}\le N_{n+1}^{\tau_{1}} $. 
Lemma \ref{lem:boh}  implies that $ \la \in \calG_{N_p}(u_n)\cap \hat{\mathcal I} $ and therefore 
also (H3) is satisfied
since we are assuming Hypothesis \ref{hyp.separazione}. 
But then Proposition \ref{multiscala} implies
\begin{equation}\label{Ln+1ind}
\bvert L^{-1}_{N_{n+1},j_{0}}(\e,\la,\theta,u_{n})\bvert_{s}\le
\frac{1}{4} N_{n+1}^{\tau}\left(N_{n+1}^{\de s}+\e
\bvert T(u_{n})\bvert_{s}\right), \quad \forall s\in[s_{0},S] \, .
\end{equation}
Finally, for $\|u-u_{n}\|_{s_{1}}\le N_{n+1}^{-\s}$ (recall that $\|u_n\|_{s_1}\le 1$ so $\|u\|_{s_1}\le 2$) one has
$$
\bvert L_{N_{n+1},j_{0}}(\e,\la,\theta,u_{n})-L_{N_{n+1},j_{0}}(\e,\la,\theta,u)\bvert_{s_{2}}
= \bvert L_{N_{n+1},j_{0}}(\e,\la,\theta,u_{n})-L_{N_{n+1},j_{0}}(\e,\la,\theta,u)\bvert_{s_{1}-\nu_1}
$$
$$\le C 
\|u-u_{n}\|_{s_{1}}\le C N_{n+1}^{-\s},
$$
where the second bound follows by \eqref{Hdiag} and \eqref{lip} with $\nu_1=\max(\nu,\nu_0)$.
Hence \eqref{Ln+1ind} and Lemma \ref{lem.inversasinistra} imply 
$ \bvert L_{N_{n+1},j_{0}}(\e,\la,\theta,u)^{-1}\bvert_{s}\le N_{n+1}^{\tau+\de s} $
for all $s\in [s_{0},s_{2}]$, proving \eqref{cisiamo}. 
\EP

\begin{lemma}\label{lem:finale}
Property (S6)$_{n+1}$ holds.
\end{lemma}

\prova
Follow word by word the proof of Lemma \ref{buonobarraS6} 
with $N_{n+1}$ instead of $N_1$, and $u_n$ instead of $u_0$. 
Since $ \la \in \gotG_{N_{n+1}} (u_n) $  (see \eqref{buoninormal2})
the bound  \eqref{Ln+1ind} holds with $ j_0 = 0  $, $ \theta = 0 $, and so
\begin{equation}\label{Ln+1j0}
\bvert L^{-1}_{N_{n+1}}(\e,\la, u_{n})\bvert_{s}\le
\frac{1}{4} N_{n+1}^{\tau}\left(N_{n+1}^{\de s}+\e
\bvert T(u_{n})\bvert_{s}\right), \quad \forall s\in[s_{0},S] \, .
\end{equation}
\EP

\subsection{Separation properties}\label{separazioneastratta}

In order to complete the proof of Theorem \ref{principe}  we show that Hypothesis \ref{hyp.separacatene}
implies Hypothesis \ref{hyp.separazione}.

\begin{prop}\label{separageo}
Hypothesis \ref{hyp.separacatene} implies Hypothesis \ref{hyp.separazione} with $C_1=(\gote+d+r+3){\mathtt s}+3$,
$ \hat{\mathcal I} := \tilde{\mathcal I}\cap\ol{\mathcal I}$ and $\hat{N}$ large enough.
\end{prop}

We split the proof of Proposition \ref{separageo} in several Lemmas. For $ \| u \|_{s_1} < 1 $,
 we consider $ L := L(\la, \theta, u) = D (\la, \theta) +\e T(u)  $ defined in \eqref{innomin}.

\begin{defi}\label{stronglygood}
A site $k=(i,\gota)\in{\mathfrak K}  $ is

\begin{itemize}

\item $(L,N)$-strongly-regular if $L_{N,i}$ is $N$-good,
\item   $(L,N)$-weakly-singular if, otherwise, $L_{N,i}$ is $N$-bad,
\item $(L,N)$-strongly-good if either it is regular for $ D = D(\la, \theta ) $ (recall Definition \ref{regular}) or all the
sites $k'=(i',\se')$ with
$\dist(k,k')\le N$ are $(L,N)$-strongly-regular. Otherwise $k$ is $(L,N)$-weakly-bad.

\end{itemize}
\end{defi}

The above definition differs from that of $(L,N)$-good matrix (Definition \ref{AN-reg}) 
in the following way. Here we do not introduce a finite subset $E$ 
but study the infinite dimensional matrix $L$, 
and require invertibility conditions on the $N$-dimensional submatrices centered at a strongly-regular point $k$
(with respect to \cite{BB1} we use a different notation, see Definition 5.1-\cite{BB1}).

\begin{lemma}\label{stronglygoodgood}
For any $j_0\in\Lambda_+$, $\chi\in [\chi_0,2\chi_0]$  consider $k=(l,j,\se)$ such that 
$|l|,|j-j_0|\le N^{\chi}$,  if 
$k$ is $(L,N)$-strongly-good   then $k$ is  
$(L_{N^\chi,j_0},N)$-good.
\end{lemma}

\prova
Set $N'=N^{\chi}$ and (recall the definition \eqref{cono})
\begin{equation}\nonumber
E= 
\gotT\times\gotS\times\gotA,
\quad
\gotT:=[-N',N']^d\cap\ZZZ^d,\
\gotS:=\Big(j_0+\Big\{ \sum_{p=1}^r \al_p{\mathtt w}_p\;:\; \al_p\in[-N',N']\Big\}\Big)
\cap\SSSS.
\end{equation}

If $k\in E$ is regular then it is $(L,N)$-good.
If $k=(l,j,\se)\in E$ is singular but $(L,N)$-strongly-regular, we  define the neighborhood $F_{N}=F_{N}(k)$ as
\begin{equation}\nonumber
F_{N}:=\gotT_{N}\times\gotS_{N}\times\gotA,
\qquad
\gotT_{N}:=\Big(\prod_{q=1}^{d}I_{q}\Big)\cap\ZZZ^d,\quad
\gotS_{N}:=\Big\{ \sum_{i=1}^r \beta_p {\mathtt w}_p \,:\;\beta_p\in J_p\Big\}\cap\SSSS,
\end{equation}
where the intervals $J_p\in \RRR$ are defined as follows (we set $a_{p}:=(j_{0})_{p}-N'$, $
b_{p}:=(j_{0})_{p}+N'$):
\begin{equation}\nonumber
\begin{aligned} 
j_{p}-a_{p}>N, \ b_{p}-j_{p}>N  \qquad & \Rightarrow \qquad J_{p}:=  [j_{p}-N,j_{p}+N], \\
j_{p}-a_{p}\le N, \ b_{p}-j_{p} > N  \qquad & \Rightarrow \qquad J_{p}:= [a_{p},a_{p}+2N] \\
j_{p}-a_{p}> N,  b_{p}-l_{p}\le N  \qquad & \Rightarrow \qquad  J_{p}:=[b_{p}-2N,b_{p}]\, , 
\end{aligned}
\end{equation}
same for $I_q$.
By construction ${\rm dist}(k,E\setminus F_{N})>N$ and ${\rm diam}(F_{N}) \le 2N
<4N$.
Moreover, by \eqref{nonsodimeglio}
there exists $\ol{k}\in E$ with ${\rm dist}(\ol{k},k)<N$ such that
$$
F_{N}=\Big((\ol{l}+[-N,N]^{d})\times\Big(\ol{\jmath}+\Big\{\sum_{p=1}^r\al_p
{\mathtt w}_p\;:\; \al_p\in[-N,N]\Big\}
\Big)\times\gotA\Big)\cap \gotK.
$$ 
Then, since $k$ is $(L,N)$-strongly-regular, the
$ \bvert (L^{F_{N}}_{F_{N}})^{-1}\bvert_{s}=
\bvert (L_{N,\ol{k}})^{-1}\bvert_{s}\le N^{\tau+\de s} $ proving the lemma.
\EP

\begin{rmk}
The assumption that $\Lambda_+$ has the ``product structure'' \eqref{nonsodimeglio} has been used only
in Lemma \ref{stronglygoodgood} above.
\end{rmk}

The Diophantine condition \eqref{dioph} implies (since $\la\ge1/2$) that 
\begin{equation}\label{dioph.meglio}
|\oo\cdot l|=\la|\ol{\oo}\cdot l |  \ge \g_{0} {|l|^{-d}} \, , \quad \forall l\in\ZZZ^{d}\setminus\{0\} \, .
\end{equation}

\begin{lemma}\label{contasingolari.nlw} 
Let $\la$ be $N$-good for $ L $, see 
Definition \ref{param.N-buoni}.  
Then, for any
$j_{0}\in\SSSS$, $\chi\in [\chi_0,2\chi_0]$,  the cardinality
\begin{equation}\label{numero.nlw}
 \# \Big\{k=(l,j_{0},\se)\in\gotK\;:\;
k\mbox{ is }(L,N)-\mbox{weakly-singular},
\;|l|\le2N^{\chi}\Big\} \le |\gotA| N^{\gote}.
\end{equation}
\end{lemma}

\prova
By Definition \ref{stronglygood}, if $(l,j_{0},\se)$ is $(L,N)$-weakly-singular, then $L_{N,l,j_{0}}(\e,\la,
\theta)$ is $N$-bad. By  \eqref{topliz} 
this means that $L_{N,j_{0}}(\e,\la,\theta+\la\ol{\oo}\cdot l)$ is $N$-bad, i.e.
$\theta+\la\ol{\oo}\cdot l\in B_{N}(j_{0},\e,\la)$, see \eqref{bad}. By assumption $\la$ is
$N$-good for $L $ and hence \eqref{partizione} holds. We claim that in each interval
$I_{q}$ there is at most one element $\theta+\oo\cdot l'$ with
$\oo=\la\ol{\oo}$ and $|l'|\le2N^{\chi}$. This, of course, imply
\eqref{numero.nlw}. Indeed, if there are $l'\ne l''$ with $|l'|,|l''|\le
2N^{\chi}$ such that $ \theta + \om \cdot l'  $, $ \theta + \om \cdot l''\in I_{q}$, then
\begin{equation}\label{dallalto}
|\oo\cdot(l'-l'')| = |(\oo\cdot l'+\theta)-(\oo\cdot l''+\theta)|
\le |I_{q}|\le N^{-\tau_{1}}.
\end{equation}
On the other hand \eqref{dioph.meglio} implies
\begin{equation}\label{dalbasso}
|\oo\cdot(l'-l'')|\ge \frac{\g_{0}}{|l'-l''|^{d}}\ge
\frac{\g_{0}}{(4N^{\chi})^{d}}=4^{-d}\g_{0}N^{-\chi d}.
\end{equation}
Clearly \eqref{dallalto} and \eqref{dalbasso} are in contradiction
   for $N\ge N_{0}$ large, because $ \tau_1>2\chi_0 d $, see \eqref{esponenti}.
\EP

\begin{coro}\label{coro1.nlw}
Let $\la$ be $N$-good for $L$. Then, for all $j_0\in\SSSS$, the number of
$(L,N)$--weakly-bad sites $(l,j_0 ,\se)\in\gotK$ with $|l|\le N^{\chi}$ is
bounded from above by $N^{\gote+r+d+1}$. 
Hence the  $(L,N)$--weakly-bad sites are a set $\Sigma_K$ as in Definition \ref{fibre} with $K= N^{\gote+d+r+1}.$
\end{coro}

\prova
By Lemma \ref{contasingolari.nlw} above, the set of $(L,N)$-weakly-singular sites $(l,j,\se)$ with
$|l|\le N^{\chi}+N$, $|j-j_0 |\le N$ has cardinality at most
$C N^{\gote}\times N^{r}$. Each $(L,N)$--weakly-bad site
$(l,j_0,\gota )$ with $|l|\le N^{\chi}$ is included in some $N$--ball
centered at an $(L,N)$-weakly-singular site and each of these balls contains
at most $CN^{d}$ sites with $j=j_0 $. Therefore there are at
most $C N^{\gote+r}\times N^{d}$ of such bad sites.
\EP

\begin{defi}\label{equivalence}
Given two sites $k,k'\in\gotK$ we say that $k\cong k'$ if there exists a $\Gamma$-chain
$\{k_{q}\}_{q=0}^{\ell}$  (Definition \ref{gammachain})
of $(L,N)$-weakly-bad sites connecting $k$ to $k'$, namely $k_0=k$ and $k_{\ell}=k'$.
\end{defi}

\noindent
{\it Proof of Proposition \ref{separageo} completed.} Let 
$ \la \in {\cal G}_N  \cap \widetilde {\cal I} $, see Definition \ref{param.N-buoni} and recall that 
$ \widetilde {\cal I} $ is introduced in Hypothesis \ref{hyp.separacatene}. 
Set $\Gamma=N^{2}$.
A $N^{2}$-chain of $(L,N)$-weakly-bad sites is formed by sites which are singular for $ D(\la,\theta)$,
see Definition \ref{stronglygood}.
Corollary \ref{coro1.nlw} implies that $  \# \Sigma_{K}^{(\tilde{\jmath})} \le K =N^{\gote+d+r+1} $, $ \forall \tilde{\jmath} \in \Lambda_+ $, 
so that Hypothesis \ref{hyp.separacatene} (for $ \la \in \widetilde {\cal I}(N_0)$ and since
$\Gamma K=N^{\gote+d+r+3}>N_0$) implies
\begin{equation}\label{ellN}
\ell\le (N^2 N^{\gote+d+r+1})^{\mathtt s} =N^{(\gote+d+r+3) \mathtt s }.
\end{equation}
The equivalence relation introduced in Definition \ref{equivalence}
induces a partition of the $(L,N)$--weakly-bad sites in disjoint equivalence classes $\Omega_{\al}$ 
with (recall  $ \Gamma = N^2 $)
$$
{\rm dist}(\Omega_{\al},\Omega_{\be})>N^{2},\qquad
{\rm diam}(\Omega_{\al}) {\le} N^{2}\ell \stackrel{\eqref{ellN}} \leq  N^{2+ (\gote+d+r+3) \mathtt s} \leq N^{C_{1}} \, , 
$$
which is \eqref{separazione}. We have verified 
 Hypothesis \ref{hyp.separazione} 
 with $ C_1 = 3 + (\gote+d+r+3)   \mathtt s  $ and 
 $ \hat{\mathcal I} := \tilde{\mathcal I}\cap\ol{\mathcal I} $.
\EP

The assumptions of Theorem \ref{principe} imply those of  Theorems \ref{thm:nm1} and \ref{thm:nm2}, by Proposition \ref{separageo}.
We fix the parameters  to satisfy  \eqref{exponents}, \eqref{esponenti}, \eqref{esponenti1}. 
Note that, by Proposition \ref{separageo},  the constant $ C_1 $  is large with $ \gote $ and so, by \eqref{esponenti1}, 
the constant  $ S'  $ has to be large with $ \gote $. 
Then,  Theorem \ref{thm:nm1} implies the existence of  
a solution $ u_\e (\la) $ of $ F(\e,\la,u_\e (\la))=0 $ for all  $ \la \in A_{\io}:=\bigcap_{n\ge0}A_n $ and   
Theorem \ref{thm:nm2} implies that $ \CCCC_\io \subseteq A_{\io} $.

\begin{lemma}\label{CinftyCep}
 $ \CCCC_\e\subseteq \CCCC_\io $ where the set $\CCCC_{\e}$
 is defined in \eqref{Cantorfinale}. 
\end{lemma}

\prova
We claim that, for all $ n \geq 0 $,  the sets $\bar\calG^0_{N_n}\subseteq \calG^0_{N_n}(u_{n-1})$
and $ \bar \gotG_{N_n} \subseteq \gotG_{N_n}(u_{n-1}) $. 
These inclusions are a consequence of the super-exponential convergence \eqref{exponentialrate} of $ u_n $ to $ u_\e $. 
In view of the definitions \eqref{buoniautovaloriinf} and \eqref{buoniautovalori}, it is sufficient to prove that
$ B^0_{N_n}(j_0,\e,\la, u_{n-1})\subseteq \bar B_{N_n}^0(j_0,\e,\la)$, $ \forall  j_0 $. 
Equivalently, if  $\theta \notin \bar B_{N_n}^0(j_0,\e,\la)$ then 
$\|L_{N_n,j_0}^{-1}(\theta, u_{n-1})\|_0 \leq N_n^{\tau_1}$, namely $ \theta \notin B^0_{N_n}(j_0,\e,\la, u_{n-1}) $ (recall 
\eqref{tetacattiviautovalori}). Indeed,  $ \|L_{N_n,j_0}^{-1}(\e, \la, \theta, u_{\e})\|_0 \leq N_{n}^{\tau_1} /2 $ by \eqref{tetacattiviautovalorifinali},
and so
 \begin{equation}\nonumber
\begin{aligned}
\|L_{N_n,j_0}^{-1}(\theta, u_{n-1})\|_0&\le \|L_{N_n,j_0}^{-1}(\theta, u_{\e})\|_0 \
\Big\| \Big(\uno+L_{N_n,j_0}^{-1}(\theta, u_{\e})(L_{N_n,j_0}(\theta, u_{n-1})- L_{N_n,j_0}(\theta, u_{\e}) )\Big)^{-1}\Big\|_0 \\
 &\le (N_{n}^{\tau_1} /2) \, 2 = N_{n}^{\tau_1}
\end{aligned}
\end{equation}
by Neumann series expansions and using Lemma \ref{controlloautov}, \eqref{lip}, 	\eqref{exponentialrate},
and \eqref{exponents}. The inclusions $ \bar \gotG_{N_n} \subseteq \gotG_{N_n}(u_{n-1}) $
follow similarly. 
\EP

The last conclusion of Theorem \ref{principe} is proved in the next section.

\subsection{Regularity}\label{sec.reg}

We now consider the case
 $ S' = + \infty $. 
The key estimate is the following upper bound for the divergence of the high Sobolev norm of the approximate solutions
 $ u_n $, which extends $(S4)_n$.
 It requires only a small modification of Lemma 
 \ref{lem:normalta.sol}.  

\begin{lemma}\label{normalitissima}
For all $s\ge S$ one has  $B_n(s) := 1+ \|u_n\|_s \le C(s)N_n^{2p} $ where $ 2p := 2(\tau+\de s_{1})+\nu+2 $.
\end{lemma}

\prova
Given $ s $, we take $ n_0(s) \in \NNN $ such that $ N_{n_0(s)} > N_0 (\Upsilon, s) $ where  $ N_0 (\Upsilon, s) $
is introduced in the multiscale Proposition \ref{multiscala}.  
For all $ n\le n_0(s)$,  the required 
bound $ B_n(s) \le C(s)N_{n}^{2(\tau+\de s_{1})+\nu+2} $ holds taking $ C(s) $ large enough.
We now prove the same bound for $ n > n_0 (s) $.
In this case, Proposition \ref{multiscala}  implies the estimate  \eqref{Ln+1j0} also for $s> S$, see 
\eqref{stimabella}. Then
\eqref{hyp.normaT}  and (P1) imply
\begin{equation}\label{fighetta}
\bvert L^{-1}_{N_{n+1}}(\e,\la,u_n)\bvert_{s}\le C(s) N_{n+1}^\tau 
\big(N_{n+1}^{\de s}+  N_{n+1}^{\nu_0/2}\|u_n\|_{s} \big) \, . 
\end{equation}
Hence, for all  $h \in E_{n+1}$, using \eqref{soboh}, \eqref{fighetta}, $ (S1)_n $, $ \nu_0/2 < \de s_1 $, 
\begin{equation} \label{quasifinale}
\|L_{N_{n+1}}^{-1}(\e,\la,u_n) [h]\|_{s} 
\leq C'(s)N_{n+1}^{\tau+\de s_1}(\|h\|_s+  (N_{n+1}^{\de(s-s_1)}+B_n(s)) \|h\|_{s_1}) \, . 
\end{equation}
Now, as in \eqref{rnbassa}, \eqref{Rnbassa}, we get
\begin{equation}\label{rs1Rs1}
\begin{aligned} 
\| r_n + R_n ( \tilde{h}_{n+1} ) \|_{s_1} & \leq  C(s) N_{n+1}^{-(s-s_{1}-\nu)/2}B_{n} (s) +  
C(s_1)N_{n+1}^{\nu}\| {\tilde h}_{n+1} \|_{s_{1}}^{2} \\
& 
 \leq C'(s) N_{n+1}^{- ( s - s_1 - \nu)/2 } B_n (s)   
\end{aligned}
\end{equation}
using \eqref{poverinoi} and 
 $ N_{n+1}^{\nu + \tau + \delta s_1} N_{n+1}^{- \sigma -1} \leq 1 $ (by  \eqref{exponents}).
Since $ \tilde{h}_{n+1} $ defined in Lemma \ref{puntofisso}
is the fixed point of  $ \calH_{n+1} $  in \eqref{contraction},  we have, using 
\eqref{quasifinale}, \eqref{633b} (with $ S \rightsquigarrow s $), \eqref{exponents}, $ \delta \in (0, 1/4) $, \eqref{rs1Rs1}
\begin{equation}\label{allafine}
\|\tilde{h}_{n+1}\|_s  
\le C'(s)N_{n+1}^{\tau+\de s_1+ (\nu/2)} 
B_n(s) \big (1+ N_{n+1}^{-(s-s_1)/2} B_n(s) \big) +C'(s)N_{n+1}^{\tau+\de s_1+\nu}\rho_{n+1}\|\tilde h_{n+1}\|_s \, . 
\end{equation}
For all $ n > n_0(s) $ (possibly larger)  
we have $  C'(s) N_{n+1}^{\tau+\de s_1+ \nu -\s-1} < 1 / 2 $ (see \eqref{exponents}). 
Moreover $N_{n+1}^{-(s-s_1)/2} B_n(s)\le 2 $ (by (P1) and $ (S1)_n $) and \eqref{allafine} implies 
\begin{equation}\label{accatilde}
\|\tilde{h}_{n+1}\|_s \le C(s)N_{n+1}^{\tau+\de s_1 + (\nu/2)}B_n(s) \, .  
\end{equation}
Therefore (recall \eqref{un+1} and \eqref{extension})
$ B_{n+1}(s) \le B_n(s) + \| \tilde{h}_{n+1}\|_s\le C'(s)N_{n+1}^{\tau+\de s_1+(\nu/2)} B_n(s) $, and so 
the sequence $B_n(s)N_{n}^{-2(\tau+\de s_{1})-\nu-2}$ is bounded.
\EP

By \eqref{extension},  Lemma \ref{normalitissima} and the estimate \eqref{accatilde} imply that (use also Lemma \ref{lem:normalta.sol} for 
$ s_1 < s \leq S $)
\begin{equation}\label{htildealtissima}
\| h_{n}\|_s  \le 
\|\tilde{h}_{n}\|_s\le C(s)N_n^{2(\tau+\de s_1)+\nu+1} \leq C(s) N_n^{2p}  \, , \quad \forall  s > s_1 \, . 
\end{equation}
Now for all $s>s_1$ let $s'=2s-s_1>s$. The interpolation inequality \eqref{interpolation} implies 
$$
\| h_{n}\|_{s} \le C(s)\| h_n\|_{s_1}^{1/2}\| h_n\|_{s'}^{1/2} \stackrel{(S2)_n, \eqref{htildealtissima}} \le C(s)
N_{n}^{- \frac{\s+1}{2} + p}\stackrel{\eqref{exponents}}{\le}
C'(s)N_n^{-1}
$$
which implies $ \|u_\e \|_{s}\le\sum_{n\ge0}\| h_{n}\|_{s}<\io $,  i.e. $u_{\e} \in H^{s}$ for all $s$.
\EP

\noindent
{\it Proof of Corollary \ref{coromerd}.} Since 
\eqref{Fsimme} holds with $\widehat{X}_s=\widehat{H}^s(\gotK)$ by assumption,  
the solution $ u_\e \in \widehat{H}^{s_1+\nu}(\gotK)$ 
by the last sentence of the Nash-Moser Theorem \ref{thm:nm1}.
\EP

\appendix\zerarcounters 


\section{Proof of the multiscale Proposition \ref{multiscala}} 
\label{prova.multiscala} 


We first prove a lemma about  left invertible block diagonal matrices. 

\begin{lemma}\label{blocchi}
Let $\DD\in\MM^{B}_{C}$ be a left invertible block
diagonal matrix, with $B\subseteq C$ i.e.
\begin{equation}\label{dstorto}
\DD^{k}_{k'}=\left\{
\begin{aligned}
&\DD^{k}_{k'}, & \mbox{if }(k,k')\in\cup_{\al}(\Omega_{\al}\times\Omega'_{\al}),\\
&\;0 & \mbox{if }(k,k')\notin\cup_{\al}(\Omega_{\al}\times\Omega'_{\al}), 
\end{aligned}
\right.
\end{equation}
where $\{\Omega_{\al}\}_{\al} $
 is a partition of $ B $, i.e. 
$ \cup_{\al}\Omega_{\al}=B $, $ \Omega_{\al}\cap\Omega_{\be}=\emptyset $, $ \forall \al \neq \be $, 
and the family $\{\Omega'_{\al}\}_{\al}$ is such that 
$ \cup_{\al}\Omega'_{\al}\subseteq C $,  $ \Omega_{\al}\subseteq\Omega_{\al}' $
and $ \Omega_{\al}' \cap\Omega_{\be}' =\emptyset $, $ \forall \al \neq \be $. 
Then $\DD$ has a block diagonal left inverse
and, given any left-inverse $L \in {\cal M}^C_B $ of $\DD$, its restriction
\begin{equation}\label{restrizione}
(R)^{k'}_{k}=\left\{
\begin{aligned}
&L^{k'}_{k}, & \mbox{if }(k,k')\in\cup_{\al}(\Omega_{\al}\times\Omega'_{\al}),\\
&\;0 & \mbox{if }(k,k')\notin\cup_{\al}(\Omega_{\al}\times\Omega'_{\al}),
\end{aligned}
\right.
\end{equation}
is a left inverse of $\DD$.
\end{lemma}

\prova
$\DD$ is left invertible and block diagonal, hence each block is left
invertible. This produces a block diagonal left inverse of $\DD$.
In order to prove that $ R $  in \eqref{restrizione} is a left inverse of $ {\cal D } $
it is sufficient to show that  
$ (L-R)\DD = 0$.
Indeed, for any $k\in B$ there is a (unique) index $\al $ such that
$k\in\Omega_{\al}$, and for any $k'\in B$ one has
\begin{equation}\label{componenti}
((L-R)\DD)^{k'}_{k}=\sum_{q\notin\Omega'_{\al}}(L-R)^{q}_{k}\DD^{k'}_{q},
\end{equation}
since, by definition \eqref{restrizione}, $(L-R)^{q}_{k}=0$ if $q\in\Omega'_{\al}$.
Now, if $k'\in\Omega_{\al}$ then $\DD^{k'}_{q}=0$ for all $q\notin
\Omega'_{\al}$ (see \eqref{dstorto}) and hence \eqref{componenti} implies $((L-R)\DD)^{k'}_{k}=0$. If, otherwise,
$ k' \in \Omega_\be $ for some $\be\ne\al$, then, by \eqref{componenti} and 
$\DD^{k'}_{q}=0 $ for all $q\notin \Omega'_{\be}$, we have 
\begin{equation}\nonumber
((L-R)\DD)^{k'}_{k}= {\mathop \sum}_{q \in\Omega'_{\be}} (L-R)^{q}_{k}\DD^{k'}_{q} = 
{\mathop \sum}_{q\in\Omega'_{\be}}L^{q}_{k}\DD^{k'}_{q}
=(L\DD)^{k'}_{k}=(\uno_{B})^{k'}_{k}=0
\end{equation}
where in the second equality we used that $ R^q_k = 0, \forall q \in\Omega'_{\be}$ (since 
$ \Omega_\be' \cap \Omega_\al' = \emptyset $), and in the third that $\DD^{k'}_{q}=0 $,  
$ \forall q\notin\Omega'_{\be}$.
\EP

Call $G$ the set of the $(A,N)$-good sites
and $B$ the set of $(A,N)$-bad sites,  see Definition \ref{AN-reg}. Let 
$ \Pi_B, \Pi_G $ be the projectors on the subspaces
$ H_B $, $ H_G $ (see \eqref{HB subspaces}) and decompose 
$ u = u_B + u_G $, $ u_B :=\Pi_B u $, $ u_G:=\Pi_G u $.

\begin{lemma}\label{controlloridotti} {\bf (Semi-reduction on the good sites).}
There exists $N_{1}=N_{1}(\Upsilon)$ such that, for $N\ge N_{1}$,
there exist $\BB\in\MM^{B}_{G}$ and $\calG\in\MM^{E}_{G}$ satisfying (recall that $\kappa= \tau+d+r+s_0$)
\begin{equation}\label{normas0}
\bvert\calG\bvert_{s_{0}}\le cN^{\ka},\qquad
\bvert\BB\bvert_{s_{0}}\le c \e \Upsilon,
\end{equation}
for some $c=c(s_{2})$ and, for all $s\ge s_{0}$, 
\begin{equation}\label{normas}
\bvert\calG\bvert_{s}\le C(s)N^{2\ka}(N^{s-s_{0}}+ \e  N^{-d-r}\bvert T
\bvert_{s+d+r}), \quad \bvert\BB\bvert_{s}\le C(s) \e N^{\ka}(N^{s-s_{0}}+ \e N^{-d-r}\bvert T
\bvert_{s+d+r}),
\end{equation}
such that if $u$ solves $Au=h$ then
\begin{equation}\label{riduco}
u_{G}=\BB u_{B}+\calG h.
\end{equation}
Conversely, if $u_{G}=\BB u_{B}+\calG h $ then, for all $k$ regular, one has
$ (Au)_{k}=h_{k} $. 
\end{lemma}

\prova
We first prove that there exist matrices  $W,R\in\MM^{E}_{G}$, satisfying
\begin{equation}\label{stimebasse}
\bvert W\bvert_{s_{0}}\le N^{\ka},\qquad
\bvert R\bvert_{s_{0}}\le C (s_2)  \e \Upsilon \, , 
\end{equation}
\begin{equation}\label{stimealte}
\bvert W\bvert_{s}\le C(s)N^{\ka+s-s_{0}},\quad
\bvert R\bvert_{s}\le \e C(s)N^{\ka}(N^{s-s_{0}}+N^{-d-r}\bvert T
\bvert_{s+d+r}), \quad \forall s \ge s_{0},
\end{equation}
such that if $u$ solves $Au=h$ then 
\begin{equation}\label{step1M} 
u_G + R u = W h  \, . 
\end{equation}
Indeed, fix $k\in G$. If $k$ is regular set $F=\{k\}$, while, if $k$ is
singular but $(A,N)$-regular, let $F\subset E$ with ${\rm diam}(F)\le 4N$
be such that ${\rm dist}(k,E\setminus F)\ge N$ and $A^{F}_{F}$ is $N$-good.
If $u$ solves $Au=h$ then
$ A^{F}_{F}u_{F}+A^{E\setminus F}_{F}u_{E\setminus F}=h_{F} $,
and hence
\begin{equation}\label{uff}
u_F + Qu_{E\setminus F}=(A^{F}_{F})^{-1}h_{F}, \quad 
Q:=(A^{F}_{F})^{-1}A^{E\setminus F}_{F}=\e (A^{F}_{F})^{-1}T^{E\setminus F}_{F} \, . 
\end{equation}
Lemma \ref{lem.algebra.bvert}, the fact that $ A_F^F $ is $ N $-good (see Definition \ref{Ngood})
 and the Hypothesis (H1) imply 
\begin{equation}\label{stimaQ}
\bvert Q\bvert_{s_{2}}\le C(s_{2}) \e \Upsilon N^{\tau+\de s_{2}}.
\end{equation}
Moreover, since ${\rm diam}(F)\le 4N $, \eqref{nanbnbna.bvert} and the first inequality in \eqref{dentrodiag.a},  
we get
\begin{equation}\label{stimaltaQ}
\bvert Q\bvert_{s+d+r} 
\le \e C(s)N^{(\de-1)s_{0}}(\Upsilon N^{s+d+r+\tau}+
N^{\tau+s_{0}}\bvert T\bvert_{s+d+r}).
\end{equation}
 Projecting \eqref{uff} onto $\{k\}$
we obtain
$ u_{k}+\sum_{k'\in E}R^{k'}_{k}u_{k'}=\sum_{k'\in E}W^{k'}_{k}h_{k'} $
where
\begin{equation}\label{ReW}
R^{k'}_{k} :=\left\{
\begin{aligned}
&Q^{k'}_{k}, & \mbox{if }k'\in E\setminus F,\\
&\;0 & \mbox{if } k'\in F, 
\end{aligned}\right.
\qquad
\mbox{ and }
\qquad
W^{k'}_{k} :=\left\{
\begin{aligned}
&[(A^{F}_{F})^{-1}]^{k'}_{k}, & \mbox{if }k'\in  F,\\
&\;0 & \mbox{if } k'\in E\setminus F 
\end{aligned}\right.
\end{equation}
which is \eqref{step1M}. 

If $k$ is regular  (see Definition \ref{regular}) then $F=\{k\}$, and, for 
$\e$  small, one has 
$ \|(A^{k}_{k})^{-1}\|_0\le 2 $. Then  
the $k$-th line of the matrix $R$ is bounded by
$ \bvert R_{k}\bvert_{s_{0}+d+r}\le $ $ \e \bvert (A^{k}_{k})^{-1}T_{k}
\bvert_{s_{0}+d+r} \le 2\Upsilon \e $ by \eqref{esponenti1.b}.

If $k$ is singular but $(A,N)$-regular then $ R^{k'}_{k}=0$ for $ \dist(k', k) \leq N$ and 
$$
\bvert R_{ k}\bvert_{s_{0}+d+r} \stackrel{\eqref{fuoridiag}} \le N^{-(s_{2}-s_{0}-d-r)}\bvert R_k
\bvert_{s_{2}} \le N^{-(s_{2}-s_{0}-d-r)}\bvert Q\bvert_{s_{2}}
 \stackrel{\eqref{stimaQ}}  \le C(s_{2})\Upsilon \e ,
$$
using also that $\tau+s_{0}+d+r-(1-\de)s_{2} < 0 $, see \eqref{esponenti1.b} and $ \delta \in (0, 1/4)$.
But then,  Lemma \ref{decadonolerighe} implies the second inequality in
\eqref{stimebasse}. The first inequality in
\eqref{stimebasse} follows in the same way. 
The first estimate in \eqref{stimealte} is a consequence of the first
in \eqref{stimebasse} and  \eqref{dentrodiag.a}. The second
estimate in \eqref{stimealte} follows by 
$$
\bvert R \bvert_{s} \stackrel{\eqref{eq.decadonolerighe}} \le K \sup_{k\in G}\bvert R_{k}\bvert_{s+d+r}
\le K\bvert Q\bvert_{s+d+r} \stackrel{\eqref{stimaltaQ}} \le C(s) \e N^{\ka}(N^{s-s_{0}}+N^{-(d+r)}\bvert T\bvert_{s+d+r}) \, . 
$$
We rewrite \eqref{step1M} as $ (\uno_{G}+R^G)u_G =W h-R^{B}u_{B} $ where $ R^G  $ 
denotes the restriction  $ R^G : G \to G $. 
By \eqref{stimebasse}, if $|\e| \leq \e_0 (s_2) $ is small enough, then 
$\bvert R^{G}\bvert_{s_{0}}<1/2$ and, by  Lemma \ref{lem.inversasinistra},
\begin{equation}\label{inversabassa}
\bvert(\uno_{G}+R^{G})^{-1}\bvert_{s_{0}}<2
\end{equation}
\begin{equation}\label{inversalta}
\bvert(\uno_{G}+R^{G})^{-1}\bvert_{s}\le C(s)(1+\bvert R^{G}\bvert_{s})
\le C(s)\big(1+ \e  N^{\ka}(N^{s-s_{0}}+N^{-d-r}\bvert T\bvert_{s+d+r})\big) \, , \ \forall s \ge s_{0} \, . 
\end{equation}
Then we obtain \eqref{riduco} with 
$ \calG:=(\uno_{G}+R^{G})^{-1} W $,  $ \BB:=-(\uno_{G}+R^{G})^{-1}R^{B} $. 
The bounds \eqref{normas0}, \eqref{normas} follow by Lemma
\ref{lem.algebra.bvert} and  \eqref{inversabassa}, \eqref{inversalta},
\eqref{stimebasse}, \eqref{stimealte}.
Finally,  \eqref{riduco} is equivalent to \eqref{step1M} which, 
for $ k $ regular, gives $  u_{k}+(A^{k}_{k})^{-1}\sum_{k'\ne k}A^{k'}_{k}u_{k'}=(A^{k}_{k})^{-1}h_{k} $.
The final assertion follows.
\EP

\begin{lemma}\label{risolvocattivi} {\bf (Reduction on the bad sites).} 
If $ u $ solves $ A u = h $ then 
\begin{equation}\label{L'eZ} 
A'u_{B}=Zh \quad 
\text{where} \quad
A':=A^{B}_{E}+A^{G}_{E}\BB\in \MM^{B}_{E},
\quad
Z:=\uno_{E}-A^{G}_{E}\calG\in\MM^{E}_{E}
\end{equation}
satisfy
\begin{subequations}
\begin{align}
&\bvert A'\bvert_{s_{0}}\le c_0,
\qquad
\bvert A'\bvert_{s}\le C(s)N^{\ka}(N^{s-s_{0}}+ \e N^{-d-r}\bvert T
\bvert_{s+d+r}),
\label{stima.a} \\
&\bvert Z\bvert_{s_{0}}\le cN^{\ka},
\qquad
\bvert Z\bvert_{s}\le C(s)N^{2\ka}(N^{s-s_{0}}+ \e N^{-d-r}\bvert T
\bvert_{s+d+r}).
\label{stima.b}
\end{align}
\label{tiodio2}
\end{subequations}
Moreover $(A^{-1})_{B}^{E}$ is a left inverse of $A'$.
\end{lemma}

\prova
If $u$ solves $ A u = h $ then $ A^{G}u_{G}+A^{B}u_{B}=h $
(we denote $A^G_E= A^G$ same for $A^B$) and, by \eqref{riduco} 
we deduce $ (A^{G}\BB + A^{B})u_{B}=h-A^{G}\calG h $, 
which is \eqref{L'eZ}.  
By the last assertion of Lemma \ref{controlloridotti}, 
for any  $ k $ regular and for all $ h $, we have  
$ \big((A^{G}\BB + A^{B})u_{B} \big)_{k}=\left(h-A^{G}\calG h\right)_{k} $
identically, namely the lines
\begin{equation}\label{suiregolari}
A'_{k}=0\quad\mbox{ and }\quad Z_{k}=0,
\quad \mbox{ for all }k\; \mbox{ regular}.
\end{equation}
That is, denoting  $ R \subset E $ the set of the regular sites in $ E $, we have $ \Pi_R A' = 0 $. 
Then  \eqref{tiodio2} follow by applying the interpolation estimates 
\eqref{nanbnbna.bvert},  \eqref{normas0},
\eqref{normas}, $ \ka > d+r $, and $\bvert D_{E\setminus R}\bvert_{s}\le 1$ for all $s$. 

Finally, 
$(A^{-1})_{B}$ is a left inverse of $A'$ 
because $A^{-1}A'=A^{-1}(A^{B}_{E}+A^{G}_{E}\BB)=\uno^{B}_{E}+\uno^{G}_{E}\BB$
which, in turn, implies $ (A^{-1})_{B}A' = \uno_{B} $. 
\EP

\begin{lemma}\label{invertoL'} {\bf (Left inverse with decay).} 
The matrix $A'$ in \eqref{L'eZ} has a left inverse
$\leftinv{{A'}}$ such that
\begin{equation}\label{sL'}
\bvert\leftinv{{A'}}\bvert_{s}\le C(s)N^{2\chi\tau_{1}+\ka+2(s_{0}+d+r)C_{1}}
(N^{C_{1}s}+ \e \bvert T\bvert_{s+d+r}), \quad \forall s\ge s_{0} \, .
\end{equation}
\end{lemma}

\prova
Let us define the matrix $\DD\in\MM^{B}_{E}$ as
$$
\DD^{k}_{k'}:=\left\{
\begin{aligned}
&(A')^{k}_{k'}, & \mbox{if }(k,k')\in\cup_{\al}(\Omega_{\al}\times\Omega'_{\al}),\\
&\;0 & \mbox{if }(k,k')\notin\cup_{\al}(\Omega_{\al}\times\Omega'_{\al}), 
\end{aligned}
\right.
$$
where the family $\{\Omega_{\al}\}_{\al\in \ZZZ}$ is the one in
Hypothesis (H3) of Proposition \ref{multiscala} and
$ \Omega'_{\al}:=\{k\in E\;:\; {\rm dist}(k,\Omega_{\al})\le N^{2}/4\}. $
First of all we prove that $\DD$ admits a left inverse $ \cal W$ with
$\| {\cal W} \|_{0}\le 2N^{\chi\tau_{1}}$. Indeed, setting $ \RR:=A'-\DD $,  we 
have that $\RR^{k}_{k'}=0$ for all $k,k'$ such that ${\rm dist}(k,k')<N^{2}/4$ and
\begin{equation}\label{controlloresto}
\bvert\RR\bvert_{s_{0}} \stackrel{\eqref{fuoridiag}} 
\le 4^{s_{2}}N^{-2(s_{2}-s_{0}-d-r)}\bvert\RR
\bvert_{s_{2}-d-r}\le 4^{s_{2}}N^{-2(s_{2}-s_{0}-d-r)}\bvert A'
\bvert_{s_{2}-d-r}  \le 
C(s_{2})N^{2\ka-s_{2}} 
\end{equation}
using \eqref{stima.a}, (H1). But then, by Lemma \ref{controlloautov}, 
$$
\|\RR\|_{0}\|(A^{-1})_{B}\|_{0}  \le \bvert\RR\bvert_{s_{0}}\|A^{-1}\|_{0}
\stackrel{\eqref{controlloresto}, (H2)} \le C(s_{2})N^{2\ka-s_{2}+\chi\tau_{1}} \stackrel{ \eqref{esponenti1.b}} \le 1/ 2,
$$
for $N$ large enough. 
Now, since $(A^{-1})_{B}$ is a left inverse of $A'$ (Lemma \ref{risolvocattivi}),
Lemma \ref{lem.inversasinistra} implies that $\DD=A'-\RR$
has a left inverse ${\cal W} $ such that (see \eqref{norma0leftinv})
\begin{equation}\label{elle2W}
\|{\cal W}\|_{0} \le 2\|(A^{-1})_{B}\|_{0}\le 2\|A^{-1}\|_{0}\le2N^{\chi\tau_{1}},
\end{equation}
by hypothesis (H2) of Proposition \ref{multiscala}. 
Now Lemma \ref{blocchi} allows to define a  block diagonal left inverse of $\DD $, denoted by 
$\leftinv{\DD}$,   as the restriction of $ {\cal W} $ as in \eqref{restrizione}. Since 
${\rm diam}(\Omega_{\al})\le  N^{C_{1}}$ (Hypothesis (H3)) then  
${\rm diam}(\Omega'_{\al})\le 2N^{C_{1}}$ and so $\leftinv{\DD}^{k'}_{k}=0$
if $ \dist (k, k' ) >2N^{C_{1}}$. Therefore, for any $s\ge0$ one has
\begin{equation}\label{leftinvds}
\bvert\leftinv{\DD} \bvert_{s} \stackrel{\eqref{dentrodiag.a}} \le C(s)N^{(s+d+r)C_{1}}\|\leftinv{\DD}\|_{0}
\le C(s)N^{(s+d+r)C_{1}+\chi\tau_{1}},
\end{equation}
by \eqref{elle2W} and Lemma \ref{blocchi}.
Finally, $A'=\DD+\RR$ and, \eqref{leftinvds}, \eqref{controlloresto} imply 
\begin{equation}\label{s0leftinvd}
\bvert \leftinv{\DD}\bvert_{s_{0}}\bvert\RR\bvert_{s_{0}}\le
C(s_{2})N^{(s_{0}+d+r)C_{1}+\chi\tau_{1}+2\ka-s_{2}}  \le 1 / 2
\end{equation}
by \eqref{esponenti1} for $N$ large enough. But then, Lemma \ref{lem.inversasinistra} implies
\begin{equation}\label{s0L'}
\bvert\leftinv{{A'}}\bvert_{s_{0}}\le 2\bvert\leftinv{\DD}\bvert_{s_{0}}
\le C(s_{0})N^{(s_{0}+d+r)C_{1}+\chi\tau_{1}},
\end{equation}
and, using 
\eqref{normaltaleftinv}, \eqref{stima.a}, \eqref{leftinvds}, we obtain \eqref{sL'}.
\EP

\noindent
{\it Proof of Proposition \ref{multiscala} completed.}
By Lemmas \ref{controlloridotti}, \ref{risolvocattivi} and \ref{invertoL'}, if $u$
solves $Au=h$ then
$ u_{G} =\calG h+\BB u_{B} $, $ u_{B} =(\leftinv{{A'}})Zh $,
which in turn implies
$$
(A^{-1})_{B}=(\leftinv{{A'}})Z,\qquad
(A^{-1})_{G}=\calG+\BB(\leftinv{{A'}})Z=
\calG+\BB(A^{-1})_{B}.
$$
Then 
\eqref{nanbnbna.bvert},  \eqref{s0L'}, \eqref{sL'}, \eqref{stima.b}, \eqref{dentrodiag.a}, \eqref{normas0}, \eqref{normas} imply, setting
$ \zeta := 2\tau_{1}+d+r+2\chi^{-1}(\ka+C_{1}(s_{0}+d+r)) $, for all $ s \in [s_0, \bar s] $, that 
$$
\bvert(A^{-1})_{B}\bvert_{s} +  \bvert(A^{-1})_{G}
\bvert_{s} \le C(s) N^{\chi\zeta}(N^{s C_{1}}+ \e \bvert T\bvert_{s}) \leq N^{\chi\tau}(N^{\chi\de s}+
\e \bvert T\bvert_{s}) / 4 
$$
by \eqref{esponenti1.a} and $ N  \geq N_0 (\bar s) $ large enough.  Thus \eqref{stimabella} is proved.
\EP

This research was supported by the European Research Council under
FP7
and partially by the PRIN2009
grant ``Critical point theory and perturbative methods for
nonlinear differential equations".

\end{document}